\documentclass[12pt]{amsart}

\usepackage{amssymb,amsmath,amsfonts}
\begin{document}
%\addtolength{\textwidth}{6truecm}
\newtheorem{theorem}{Theorem}[section]
\newtheorem{lemma}{Lemma}[section]
\newtheorem{proposition}{Proposition}[section]
\newtheorem{corollary}{Corollary}[section]
\newtheorem{conjecture}{Conjecture}[section]

\theoremstyle{definition}
\newtheorem{definition}{Definition}[section]
\theoremstyle{remark}
\newtheorem{remark}{Remark}[section]

\numberwithin{equation}{section}

%    Absolute value notation
\newcommand{\abs}[1]{\lvert#1\rvert}

%    Blank box placeholder for figures (to avoid requiring any
%    particular graphics capabilities for printing this document).
\newcommand{\blankbox}[2]{%
  \parbox{\columnwidth}{\centering
%    Set fboxsep to 0 so that the actual size of the box will match the
%    given measurements more closely.
%    \setlength{\fboxsep}{0pt}%
%    \fbox{\raisebox{0pt}[#2]{\hspace{#1}}}%
}%
}
%%%%
%\newcommand{\be}{\begin{eqnarray}}
%\newcommand{\ee}{\end{eqnarray}}
%\newcommand{\bes}{\begin{eqnarray*}}
%\newcommand{\ees}{\end{eqnarray*}}
%\newcommand{\beqn}{\begin{equation}}
%\newcommand{\eeqn}{\end{equation}}

%\renewcommand{\theenumi}{(\roman{enumi})}

\newcommand\Pa{Painlev\'e}
\newcommand{\PI}{${\rm P}_{\textrm{\rm\tiny I}}$}
\newcommand\Pp{\Pi_+}
\newcommand\Pm{\Pi_-}
\title[The First Painlev\'e Equation]{The Dirichlet Boundary Value
Problem for Real Solutions of the First Painlev\'e Equation on
Segments in Non-Positive Semi-Axis}
\author{N. Joshi}
\address{School of Mathematics and Statistics, University of Sydney,
NSW 2006, Australia}
\email{nalini@maths.usyd.edu.au}
\thanks{Research supported by Australian Research Council Fellowship
\#F69700172 and Discovery Project Grant \#DP0208430}
\author{A. V. Kitaev}
\address{
Steklov Mathematical Institute, Fontanka 27, St. Petersburg,
191023, Russia} \email{kitaev@pdmi.ras.ru} 
%\commby{XXX}
\date{13 February 2003; revised 06 July 2004}
\begin{abstract}
We develop a qualitative theory for real solutions of the equation
$y''=6y^2 -x$. In this work a restriction $x\leq0$ is assumed.
An important ingredient of our theory is the introduction of several new
transcendental functions of one, two, and three variables that describe
different properties of the solutions. In particular, the results obtained
allow us to completely analyse the Dirichlet boundary value problem
$y(a)=y^0$, $y(b)=y_0$ for $a<b\leq0$.
\end{abstract}
\maketitle

\section{Introduction}
\subsection{}%{\bf\large1.1}
In 1898, in the framework of the theory of
functions of one complex variable and the analytic theory of ODEs,
P. Painlev\'e \cite{pp:first} introduced the equation
\begin{equation}
 \label{p1}
        y'' = 6y^2 -x,
\end{equation}
now called the first Painlev\'e equation. The way it was introduced suggested
the following two main directions of research:\\
(1)
The proof of the transcendency of the solutions of Equation~(\ref{p1}),
i.e., that they are actually ``new'' functions that cannot be expressed
(within a certain class of operations) in terms of already known functions;\\
(2) The study of solutions as the functions of complex variable
$x$. This concerned mainly the singularity structure of its
solutions, that is:
(i) the proof of the Painlev\'e property,
i.e., all solutions are meromorphic functions in $\mathbb C$
possessing only second order poles of strength unity (at any pole
$x_0$ the leading term of the Laurent expansion is $1/(x-x_0)^2$);
(ii) the distribution of the poles; and (iii) the structure of the
essential singularity at the infinity, i.e., asymptotic behaviour
of the solutions in a neighbourhood of $x=\infty$.

Both directions were substantially worked out by Painlev\'e. The singularity
structure and asymptotic behaviour in a neighbourhood of infinity is due to
Boutroux \cite{boutroux:I}. Contemporary researchers actually have to make
considerable efforts in clarifying the many statements made by Painlev\'e
and Boutroux and putting them on a solid mathematical foundation. We do not
review here the subsequent modern studies, this requires a special focus
and considerable space, whilst they do not concern directly this work. To
complete the general picture of the ``complex theory'', the connection
formulae, which allow to connect asymptotics as $x\to\infty$ for different
$\arg x$ of the solutions should be mentioned (see \cite{K} for a review),
as they were not discussed in the classical works.

\subsection{}%{\bf\large1.2}
Real solutions of Equation~(\ref{p1}) appear in applications.
One of the sources is similarity reductions of integrable PDEs.
In particular, there are such reductions for the famous Korteweg de Vries,
Boussinesq, and Kadomtsev-Petviashvili equations (see \cite{AC}). To study
related similarity regimes a knowledge of the behaviour of real
solutions for finite values of $x$ is important.

There is another, probably even more important source of applications of
equation~(\ref{p1}) viz. in asymptotic analysis of nonlinear PDEs, ODEs, and
difference equations. This concerns both integrable and non-integrable
situations. In this analysis Equation~(\ref{p1})
serves as a model equation in the description of asymptotics in
transition layers and caustic-type domains \cite{H,K2}, so that particular
solutions of Equation~(\ref{p1}) at some {\it finite} point $x$ appears in
the corresponding asymptotic formulae. The point $x$ in such formulae serves
as a parameter describing the domain of validity of these asymptotics.
To further use these asymptotics, say, in the techniques related with the
matching of asymptotic expansions, the knowledge of behaviour of real
solutions of Equation~(\ref{p1}) for finite $x$ is important.

\subsection{}%{\bf\large1.3}
Despite the growing importance of the real solutions for applications,
there are only a few papers about them in the literature. Some information
concerning qualitative behaviour of real solutions of Equation~(\ref{p1})
in the finite domain can be found in the following works
\cite{bart:73,HT,HS,jablonski1,JK,N,jablonski2}. Among them the main emphasis
in \cite{HS,jablonski2} is made on asymptotic analysis of solutions for
$x\geq0$ and the only paper \cite{bart:73} concerns the solutions for
$x\leq0$.

The major problem involved here is that the standard methods of
qualitative analysis for ODEs applicable to Equation~(\ref{p1}), based on
the maximum principle \cite{PW} or equivalent statements, work
only for solutions whose graphs in $(x,y)$-plane are located in the
upper half plane $y\geq0$, i.e., positive solutions. At the same time integral
curves corresponding to the most number of interesting solutions cross
the the $x$-axis and than return back to the upper half-plane.
Therefore, whilst the analysis in the upper half plane is almost trivial,
still little is known about the global behaviour of the integral curves of
Equation~(\ref{p1}). Another interesting problem, how to connect asymptotic
description of the solutions, which is now well known for all solutions
complex and real, with their Cauchy initial values, even on the qualitative
level, still remains unsolved.

We address these problems in our series of works on the qualitative
studies of Equation~(\ref{p1}). The first work \cite{JK} is devoted to the
most famous Boutroux tritronqu\'ee solution of Equation~(\ref{p1}), which
very often appears in different applications. This is the second work,
where we develop the theory of the integral curves for $x\leq0$ and our
third work \cite{JK3} provides a qualitative description of the integral
curves for $x\geq0$ and, in particular, of the Cauchy initial
value problem for this semi-axis. The works
\cite{HT,HS,jablonski1,JK,N,jablonski2} concern real solutions on the
non-negative semi-axis. We give their overview in our \cite{JK3}. The only
work devoted to the qualitative studies of solutions of Equation~(\ref{p1})
on the non-positive semi-axis, we were able to find, is \cite{bart:73}.
The results of the last work are reviewed in Subsection~\ref{subsec:1.5}
below, along with a description of our work.

\subsection{}%{\bf\large1.4}
Before explaining our results it is important to stress that the
word {\bf solution} is used throughout this paper to denote a
{\it real-valued solution of Equation}~(\ref{p1}) which is defined in
the {\bf interval of existence}, i.e., {\it the largest connected open
domain in $\mathbb R$ where the solution has no poles}. Moreover,
{\it solutions that have different intervals of existence are referred
to as different solutions}, despite the fact that some of them
could be analytic continuations of each other in the complex $x$-plane,
i.e., represent the same meromorphic solution from the
standpoint of analytic theory of ODEs. In a few places where we mean %In few->In a few
solutions in the latter sense we call them {\bf complex
solutions}. Note that this convention is different from the one in
our previous work \cite{JK}, where by real solutions we meant
those complex solutions whose imaginary part vanishes on the real
axis. Real solutions in the latter sense have an infinite number of  %last->latter
poles on the real axis.

\subsection{}%{\bf\large1.5}
 \label{subsec:1.5}
In our previous work we remarked (cf. \cite{JK}, Remark 2) that the
intervals of existence of real solutions, whose closure contains $x=0$,
are uniformly bounded on the left. In \S~\ref{sect:2} we prove the
uniform lower boundedness of intervals of existence of solutions that
are regular at any point $x_0\leq0$, i.e., for any point $x_0$ there
exists a minimal finite interval $(X(x_0),x_0)$, such that all solutions
regular at $x_0$ have a pole in this interval. Further, for brevity, we
call this property {\it the property of uniform boundedness}. This is a
specific property of the real solutions; it is not valid for the
restriction of general complex solutions to the real axis. To see it,
we recall the special solutions, which Boutroux called
``int\'egrales tronqu\'ees'' (see \cite{boutroux:I,JK}). Their
restrictions to the real axis have intervals of existence which
contain the whole negative semi-axis. The latter property is not
valid for complex solutions even if we exclude the int\'egrales
tronqu\'ees as there are complex solutions that approximate them
on intervals of arbitrary length.

The property of uniform boundedness has not been observed before and does
not follow from the results known about the real solutions in the
literature.
It is known that asymptotics as $x\to-\infty$ of all real solutions, if
we understand them in the sense of analytic continuation, are
given by a Boutroux type asymptotic formula in terms of the Weierstrass
$\wp$-function with fixed period parallelogram. One can deduce from this
that all such solutions have an infinite number of poles on the real axis
and find their asymptotic distribution \cite{boutroux:I,njmdk:direct,K}.
However, the property of uniform boundedness cannot be derived from the
asymptotic results since, there is no universal value $X$ such that for all
$x<X$ the correction term to the asymptotic formula would have a uniformly
(with respect to all real solutions) small estimate. If we suppose that such
$X$ exists, we immediately arrive at a contradiction by taking for any
$x_0<X$ initial data, such that the corresponding solution does not have the
universal behaviour prescribed by the Boutroux asymptotic formulae in the
neighbourhood of $x_0$, though, of course, for sufficiently large $x$ this
solution fits the Boutroux asymptotic regime. The existence of such data
is easy to observe: if we assume the uniform estimate of the asymptotic
correction term, we can estimate the initial data of the Weierstrass
$\wp$-function in terms of the initial data for the Equation~(\ref{p1}).
The latter initial data can be taken arbitrary, whilst the data for the
$\wp$-function are related via the Weierstrass differential equation of the
first order with fixed invariants.

It is worth mentioning that Bartashevich \cite{bart:73} was very close to
formulation of the property of uniform boundedness. He observed a partial
uniformity with respect to the initial slope of the solutions.
More precisely, the property he noticed can be formulated as follows.
Let $(a,b)$ with $a<0$ be the interval of existence of solution $y(x)$ of
Equation~(\ref{p1}) with initial data: $y(x_0)=y_0$ and $y'(x_0)=y_1$ at
$x_0<0$. Then $X(x_0,y_0)\equiv\underset{y_1}{\inf}\{a(x_0,y_0,y_1)\}$ is
finite. However, the proof given in \cite{bart:73} is incomplete, as it
addresses only the simplest case, in our notation, $y(x)>y_0>0$ for $x<x_0$
or, equivalently, $y_1\leq0$. Moreover, his estimate is divergent as
$y_0\to+0$.

To prove and further study the property of uniform boundedness we introduce
in \S~\ref{sect:2} three special functions: $X(x_0)$, $X_-(x_0)$, and
$X_{min}(x_0)$. We not only prove the lower boundedness of these
functions but also find the first nontrivial term of their asymptotic
expansion as $x_0\to-\infty$. This term for the function $X(x_0)$ looks
very similar to the leading term of asymptotics of the difference of two
neighbouring poles (the pole spacing) for real solutions that is governed by
the Boutroux asymptotic regime. This leading term of the pole spacing
is universal for all real solutions.
However, a more precise comparison (to appear elsewhere) shows that
asymptotics of the pole spacing has a smaller coefficient than that for
$X(x_0)-x_0$, i.e., $2C$ in Equation~(\ref{eq:X_asymptotics}). So that as
$x_0\to-\infty$ the asymptotic spacing of the poles is not the maximal
possible on the real axis.

\subsection{}%{\bf\large1.6}
The central technical idea of the paper is to reparametrize solutions,
$y(x)$, which are originally considered as the functions of initial data,
$y_0=y(x_0)$ and $y_1=y'(x_0)$, i.e., the functions $y(x)=y(x;x_0,y_0,y_1)$
(or $y(x)=y(x;x_0,c)$, where $c$ is the so called pole parameter, if $x_0$
is a pole), in terms, as we call it, the {\it level parametrization}:
$y(x)=y(x;x_0,y_0,y_l)$ (or, respectively, $y(x)=y(x;x_0,y_l)$), where $y_l$
is the minimum value of the corresponding solution. Most of the results in
Section~\ref{sect:3} aim to establish a possibility of this
reparametrization and different properties of the functions
$y_1=y_1(x_0,y_0,y_l)$ and $c=c(x_0,y_l)$ defining it. In particular, we
found that $c=c(x_0,y_l)$ is a strictly monotonically increasing function
of the second argument. This property can be viewed as a ``visualization''
of the parameter $c$ which is, in analytical definition, hidden being the
fourth coefficient of the Laurent expansion at $x_0$.

Our theory is essentially real-valued. Even where we prove the smoothness
of our reparameterization, we substantially exploit specific properties of
real solutions. Unexpectedly, one pure analytic fact plays an
important role in our study of the geometry of the integral curves.
This is the existence of an analytic mapping of the parameters $(x_0,c_0)$
defining the Laurent expansion of the general complex solution at $x_0$ to the
parameters $(x_1,c_1)$ of the Laurent expansion of the same solution at
$x_1$. This mapping was actually discussed by Boutroux in his landmark
work \cite{boutroux:I}, where he made deep comments about this mapping
lying in the heart of the Painlev\'e method\footnote{In Boutroux's own words
(p.261 of \cite{boutroux:I}): ``C'est l\`a un fait gros de cons\'equences:
en approfondissant l'\'etude de la fonction $\bar X_1(\bar X_0,C)$ et des
fonctions connexes, il semble que nous touchions au c{\oe}ur des nouveaux
\^etres analytiques introduits dans la Science par M. Painlev\'e.''
Here the notation $\bar X_1(\bar X_0,C)$ coincides with our
$x_1(x_0,c_0)$.}.
In particular, he indicated
that an attempt to understand the structure of its singularities can be
viewed as the motivation for his studies. Boutroux, actually, did not give
a definition of this mapping, possibly, considering it as self-evident.
We give an accurate definition of the mapping in
Proposition~\ref{Prop:analyticity} of Section~\ref{sect:3}. It shows that
the existence of such mapping relies upon a specific dependence of the Laurent
expansion on parameters $x_0$ and $c_0$, so that for equations of a
non-Painlev\'e type there could be some situations, where this mapping does
not exist, or has some unusual properties.

\subsection{}%{\bf\large1.7}
One of the basic technical tools that make it possible to study the
functions $X(x_0)$ and  $X_-(x_0)$ as well as to establish many properties %function->functions
of the integral curves is Lemma~\ref{Lem:4.3} and its two limiting cases
Lemmas~\ref{Lem:4.1} and \ref{Lem:4.2}. These Lemmas concern the behaviour
of the integral curves to the left of their minima. We provide slightly
different proofs for the last two lemmas in Section~\ref{sect:4}, both due
to their importance and having in mind possible generalizations. In
Sections~\ref{sect:4} and \ref{sect:5} we explore the opportunities suggested
by these Lemmas to study the functions $X(x_0)$ and,  correspondingly,
$X_-(x_0)$. In Section~\ref{sect:5} we actually study the function
inverse to $X_-(x_0)$ since we are considering there solutions as being
initially defined at the left bound of their interval of existence.

\subsection{}%{\bf\large1.8}
Another very important technical tool, the Moore--Nehari
Lemma~\ref{Lem:MN}, is invoked in Section~\ref{sect:main}. This Lemma
concerns the number of intersections of the integral curves. It is convenient
for us to call two solutions having a common pole,
{\it solutions intersecting at the point at infinity}. In particular,
we prove that the Moore--Nehari lemma can be generalized to include the
points of intersection at infinity. We call this generalization
of Moore--Nehari lemma Projective Lemma~\ref{Lem:MN}.

The latter Lemma allows us to prove the important statement that for any
$x_0$ the solution which has the maximal interval of existence
$(X(x_0),x_0)$ is unique. Another interesting property that we establish
in Section~\ref{sect:main} is that the integral curves can be viewed as the
lines of a model of geometry where: any two lines can intersect only at
two, one, or zero points. Moreover, the usual duality principle of projective
geometry, with the interchange of lines and points, is valid: any two
points can belong only to two, one, or zero lines.

The latter result can be reformulated in terms of the Dirichlet boundary
value problem; we pursue this in Section~\ref{sect:boundary}. To make
such a reformulation we introduce and study an auxiliary function,
$Z=Z(x_0,y_0,y^0)$. This function defines the largest segment $[Z,x_0]$
where the boundary value problem for Equation~(\ref{p1}), $y(x_0)=y_0$ and
$y(Z)=y^0$ has a solution. We prove that the problem with the same boundary
values on the segments $[z,x_0]$ have: two solutions if $z\in(Z,x_0)$,
one solution if $z=Z$, and no solutions if $z<Z$.

The function  $Z(x_0,y_0,y^0)$ for finite values $y_0$ and $y^0$ can be
calculated numerically with \textsc{Maple 8} code by using the
\texttt{dsolve} procedure for boundary value problems. However, it is
important to know an actual number of the solutions as this procedure, in
case of existence, gives always only one solution.

\subsection{}%{\bf\large1.9}
 \label{subsec:1.9}
We provide a numerical illustration of some results obtained in
\S\S~\ref{sect:2}--\ref{sect:5} in the last Section~\ref{sect:numerics}.
This illustration requires rather precise calculations as some data
for different solutions are very close.

In this Section we also formulate four conjectures. The first one
concerns an interesting property of {\it approximate symmetry} for a
function $f(x_0,y_0,y_l)$ describing initial slope of solutions at $x_0$
in terms of their initial, $y_0$, and the minimum, $y_l$, values.
The other three are about uniqueness and behaviour of special solutions
corresponding to the functions $X_{min}(x_0)$ and $\Xi_{min}(x_0)$.

\subsection{}%{\bf\large1.10}
In our next work~\cite{JK3} we, in particular, consider the continuation of
the functions $X_{min}(x_0)$, $X_-(x_0)$, and $X(x_0)$ to the whole real
axis. Some interesting further questions, apart of the conjectures mentioned
above, concern further terms in the asymptotic expansions of
these functions derived in Section~\ref{sect:2}. For example, the leading
terms of asymptotics for the functions $X_{min}(x_0)$ and  $X_-(x_0)$
coincide, however, the next terms should be
different. Another interesting further question might be about
analytic continuation of these functions to the complex plane. The natural
question about the Neumann and mixed Neumann--Dirichlet boundary value
problems on the non-positive semi-axis is closely related with the uniqueness
mentioned in the last sentence of Subsection~\ref{subsec:1.9}.

An interesting development would be a ``continuation'' of the
theory developed here to a wider class of second order ODEs.
\section{Uniform Bound on Intervals of Existence}
 \label{sect:2}
\begin{proposition}
 \label{Prop:real pole}
The interval of existence of any solution of Equation~{\rm(\ref{p1})} is
bounded below.
\end{proposition}
\begin{proof}
Consider solutions whose intervals of existence have non-empty intersection
with the negative semi-axis.
Change variables to $t=-x>0$, $u(t)=y(x)$, then Equation (\ref{p1}) becomes
\begin{equation}
 \label{p1t}
u''(t)=6u(t)^2+t.
\end{equation}
Assume there is no pole for $t>0$. Since $u''(t)>t$, there exists
a point $t_0$ such that for
$t>t_0$, $u(t)$ is monotonically increasing to $+\infty$.
We assume that $t_0$ is chosen such that $u_0=u(t_0)>0$, $u_1=u'(t_0)>0$.

We have
$
u''(t)>6u(t)^2.
$
Integration gives
$
u'(t)^2/2>2u(t)^3+N(u_0,u_1),
$
where
$
N(u_0,u_1)=u_1^2/2-2u_0^3.
$
Choose now $t_1\geq t_0$ such that for all for $t\geq t_1$,
$
3u(t)^3/2+N(u_0,u_1)>0.
$
Hence
$
u'(t)^2>u(t)^3
$
for $t>t_1$.

Integration gives
\[ - \frac1{\sqrt{u(t)}} + \frac1{\sqrt{u(t_1)}} > \frac{t-t_1}2 \]
which is a contradiction.
\end{proof}
\begin{remark}{\rm
 \label{Rem:1}
The proof of Proposition~\ref{Prop:real pole} shows that
every solution with a positive derivative at some point $x_0\leq0$
has a unique minimum, $x_{min}$, achieved to the left of $x_0$.
}\end{remark}
\begin{definition}
 \label{Def:2.1}
Let $x_0\leq0$, define
$$
{X}_{min}(x_0)\equiv
\underset{y_0\in\mathbb{R},y_1\geq0}{\inf}\,x_{min},
$$
where $x_{min}=x_{min}(x_0)<x_0$ is the minimum of the solution $y(x)$ with
the initial values: $y(x_0)=y_0$ and $y'(x_0)=y_1\geq0$.
\end{definition}
\begin{lemma}
 \label{Lem:2.1}
The function ${X}_{min}(x_0)$ is finite for all $x_0\leq0$.
\end{lemma}
\begin{proof}
Without loss of generality we can assume that $x_0<0$: the boundedness of
${X}_{min}(0)$ follows from the boundedness of ${X}_{min}(x_0)$ for
$x_0\neq0$.
Actually, for any $\varepsilon>0$, either a solution with initial data given
at $x_0=0$ blows up in the semi-segment $[-\varepsilon,0)$, or its minimum
is lower-bounded by
${X}_{min}(-\varepsilon)\leq {X}_{min}(0)<0$.
In fact, as $\varepsilon\to0$, by continuity,
$X_{min}(-\varepsilon)\to X_{min}(0)$, therefore
$X_{min}(0)=\underset{0<\varepsilon}{\sup}\{X_{min}(-\varepsilon)\}$.

We again use notation (\ref{p1t}) for Equation~(\ref{p1}).
Denote $t_0=-x_0>0$, $t_{min}=-x_{min}$, $u(t_0)=y_0\equiv u_0$, and
$u'(t_0)=-y_1\equiv u_1<0$.
Consider now Equation~(\ref{p1t}) in the segment $[t_0,t_{min}]$.
Multiplying it by $u'(t)$, integrating then from $t_0$ to $t$,
and solving for $u'(t)$, we obtain,
\begin{equation}
 \label{eq:u'}
-u'(t)=\sqrt{4u^3+2N(u_0,u_1)+2\int_{t_0}^tsu'(s)ds},
\end{equation}
where $N(u_0,u_1)=u_1^2/2-2u_0^3$ and
the positive branch of the square root is assumed.
Putting $t=t_{min}$ in Equation~(\ref{eq:u'}) we get,
\begin{equation}
 \label{eq:umin}
0=4u_{min}^3+2N(u_0,u_1)+2\int_{t_0}^{t_{min}}tu'(t)dt.
\end{equation}
Subtracting now the right-hand side of Equation~(\ref{eq:umin})
from the expression under the square root in Equation~(\ref{eq:u'}),
dividing both sides of the latter equation by this square root,
and integrating from $t_0$ to $t$ we arrive at
\begin{equation}
 \label{eq:u-int}
t-t_0=\int_{u(t)}^{u_0}
\frac{d u}{\sqrt{4(u^3-u_{min}^3)+2\int_{u_{min}}^{u}t(\hat u)d\hat u}}.
\end{equation}
Note that $u(t)$ is monotonic on $[t_0,t_{min}]$, therefore its
inverse, $t(u)$, is properly defined.
Now, substituting $t=t_{min}$ into Equation (\ref{eq:u-int}), noting that
$t(\hat u)\geq t_0$, and changing variables to  $u=\sqrt{t_0}v$,
we get the following estimate
\begin{eqnarray}
 \label{eq:t-v-est}
&&t_{min}-t_0<\frac{I(v_0,v_{min})}{t_0^{1/4}},\\
\label{eq:t-v-est-second}
&&I(v_0,v_{min})\equiv
\int_{v_{min}}^{v_0}\frac{d v}{\sqrt{4(v^3-v_{min}^3)+2(v-v_{min})}},
\end{eqnarray}
where, $v_{min}=u_{min}/\sqrt{t_0}$ and  $v_0=u_0/\sqrt{t_0}$.

The function $I(v_0,v_{min})$ is bounded for all values of $v_0\in\mathbb R$
and $v_{min}<v_0$. Actually,
$
I(v_0,v_{min})<I(+\infty,v_{min}).
$
The function $I(+\infty,v_{min})$ is a continuous function of
$v_{min}\in\mathbb{R}$ which vanishes as $v_{min}\to\pm\infty$:
$$
{\rm for}\quad
|v_{min}|\geq\epsilon>0,\quad{\rm we\;have}\quad
I(v_0,v_{min})<\frac{I_\nu}{2\sqrt{|v_{min}|}},\quad
\phantom{}\nu={\rm sign}\,\{v_{min}\}1,
$$
where
\begin{equation}
 \label{eq:intdelta}
I_\nu:=\int_\nu^{+\infty}\frac{dw}{\sqrt{w^3-\nu^3}}=
2\int_0^{+\infty}\frac{dv}{\sqrt{v^4+3\nu v^2+3\nu^2}}.
\end{equation}
Thus,
$
\underset{v_{min}<v_0\in\mathbb{R}}{\sup}\,I(v_0,v_{min})
$
exists and
$$
\underset{v_{min}<v_0\in\mathbb{R}}{\sup}\,I(v_0,v_{min})=
\underset{v_{min}\in\mathbb{R}}{\max}\,I(+\infty,v_{min}),
$$
where the maximum is achieved at some finite value of $v_{min}$.
Under the change of variables, $\sqrt{v-v_{min}}=w$ one finds that
$$
I(+\infty,v_{min})=\int^{+\infty}_0
\frac{dw}{\sqrt{w^4+3v_{min}w^2+3v^2_{min}+1/2}},
$$
so that it is clear that the maximum is achieved at a negative value of
$v_{min}$. Now, introducing the constant
\begin{equation}
 \label{eq:c}
C=\underset{v_{min}<0}{\max}\int^{+\infty}_0
\frac{dw}{\sqrt{w^4+3v_{min}w^2+3v^2_{min}+1/2}}
\end{equation}
we get (cf. (\ref{eq:t-v-est})) the estimate
\begin{equation}
 \label{eq:est-xmin}
x_0-\frac C{|x_0|^{1/4}}<X_{min}(x_0)<x_0<0.
\end{equation}
\end{proof}
\begin{remark}
 \label{Rem:2.2}
We can further specify a location of the maximum, $v_{min}^{max}$, in
Equation~(\ref{eq:c}).
Indeed, making the change of variables: $w=\sqrt{-v_{min}}u$,
$2x=1/v_{min}^2$, one rewrites Equation~(\ref{eq:c}):
$
C=\underset{x>0}{\max}\,(2x)^{\frac14}\int^{+\infty}_0
\frac{du}{\sqrt{u^4-3u^2+3+x}}.
$
Now, differentiating the function under the $\max$-sign with respect to $x$
and combining all the resulting terms under one integral, one observes that
the integrand is non-negative for $x\leq3/4$. Thus, all extrema of the
integral in Equation~(\ref{eq:c}) are located in the interval,
$v_{min}\in(-\sqrt{2/3},0)\subset(-0.8165,0)$.
Numerical calculations (with \textsc{Maple 8}) show that there is only
one extremum, the global maximum, which is achieved at
$$
v_{min}^{max}=-0.22600387635302095\ldots,\;\;
x_{max}=\frac1{2(v_{min}^{max})^2}=9.78899773742578347\ldots\,.
$$
This many ($18$) decimal digits allows us to calculate $36$ decimal digits
of $C$:
\begin{equation}
 \label{eq:c_num}
C=2.32470720434237566413065947435242998\ldots\;.
\end{equation}
Our attention to the precision of calculation of $C$ is clarified by
Corollary~\ref{Cor:X_min_asymptotics}.
\end{remark}
\begin{remark}
 \label{Rem:X_min_0}
For small $-1\leq x_0<0$ the estimate for $X_{min}(x_0)$ in
Equation (\ref{eq:est-xmin}) becomes too rough. The reason lies in the
rough approximation we made of the integral inside the square root in
Equation (\ref{eq:u-int}). For these values of $x_0$, we can improve
the estimate: $-1-C<X_{min}(x_0)<x_0$. The improved estimate is also
valid for $x_0=0$. This follows by similar arguments as in the first
paragraph of the proof of Lemma \ref{Lem:2.1}. The numerical
value of $C$ (\ref{eq:c_num}) allows us to make a slight improvement of
the above small-$x_0$-estimate. Actually, the function $x_0-C/|x_0|^{1/4}$
has the minimum $\hat x_0=-\left(C/4\right)^{4/5}=
-0.64780846\ldots\in[-1,0]$.
Thus for $x_0\leq\hat x_0$ Inequality~(\ref{eq:est-xmin}) still delivers
the best available estimate, while for $x\in[\hat x_0,0]$ the better
estimate is $-5\left(C/4\right)^{4/5}<X_{min}(x_0)<x_0$. However, this is
really a minor improvement as the numerical values
$5\left(C/4\right)^{4/5}=3.23904230\ldots$ and $1+C=3.32470720\ldots$ are
very close.
\end{remark}
\begin{corollary}
 \label{Cor:X_min_asymptotics}
For $x_0\leq0$, denote by $\eta(x_0)$ the unique positive root of the
equation $\eta^5-|x_0|\eta=C$. Then,
\begin{equation}
 \label{ineq:eta_est}
X_{min}(x_0)<-\eta(x_0)^4=x_0-\frac C{\eta(x_0)}.
\end{equation}
For $x_0<0$,
\begin{equation}
 \label{eq:x_min_double_est}
x_0-\frac C{|x_0|^{1/4}}<X_{min}(x_0)<x_0-
\frac C{\big||x_0|+C/|x_0|^{1/4}\big|^{1/4}}.
\end{equation}
In particular, as $x_0\to-\infty$,
$$
X_{min}(x_0)=x_0-\frac C{|x_0|^{1/4}}+
\mathcal{O}\left(\frac1{|x_0|^{3/2}}\right),
$$
where
$0<{\mathcal O}\big(\frac1{|x_0|^{3/2}}\big)<\frac{C^2}{4|x_0|^{3/2}}$
for $|x_0|>C^{4/5}$.
\end{corollary}
\begin{proof}
We turn back to the notation of Lemma~\ref{Lem:2.1} and consider
Equation~(\ref{eq:u-int}). However, now we estimate the function
$t(\hat u)$ from above $t(\hat u)<t_{min}$ and arrive at the following
estimate,
\begin{equation}
 \label{ineq:tmin_below}
t_{min}-t_0>\frac{I(v_0,v_{min})}{t_{min}^{1/4}},
\end{equation}
where now, $v_{min}=u_{min}/\sqrt{t_{min}}$ and  $v_0=u_0/\sqrt{t_{min}}$.
The key point now is that there exists a solution of Equation~(\ref{p1})
which corresponds to the value $v_{min}=v_{min}^{max}$ (see the last
paragraph of this proof). For this solution,
Inequality~(\ref{ineq:tmin_below}) reads
$(t_{min}-t_0)t_{min}^{1/4}>C$. Now, noting
that $\sup t_{min}=|X_{min}(x_0)|$ is finite according to
Lemma~\ref{Lem:2.1}, we obtain
$
(|X_{min}(x_0)|-|x_0|)|X_{min}(x_0)|^{1/4}>C.
$
So that $|X_{min}(x_0)|^{1/4}>\eta(x_0)$ and  we arrive at
estimate~(\ref{ineq:eta_est}). The next simpler inequality is obtained
by substitution of the first Inequality~(\ref{eq:est-xmin}) instead of
$|X_{min}(x_0)|^{1/4}$ in the derivation above.

We now prove the existence of the solution of Equation~(\ref{p1})
corresponding to the value $v_{min}=v_{min}^{max}$ for all $x_0$. This
proof is based on the results established in \S~\ref{sect:3} and can be
omitted at
the first reading as it is not used anywhere in the subsequent sections.
This solution must have a pole at $x_0$, and such minimum, $x_{min}$, that
$y_l=v_{min}^{max}\sqrt{|x_{min}|}$, where $y_l$ is the corresponding
minimum value. From Theorem~\ref{Thm:3.1} it follows that for any $x_0\leq0$
and $y_l\in\mathbb{R}$, there is a unique solution of Equation~(\ref{p1}).
The corresponding minimum $x_{min}=x_{min}(x_0,y_l)<x_0\leq0$ for all values
of its arguments is a continuous function (see the proof of
Corollary~\ref{Cor:3.3}). Moreover, it is bounded below according to
Lemma~\ref{Lem:2.1}. Thus, the function $y_l/\sqrt{|x_{min}(x_0,y_l)|}$ is a
continuous function of $y_l\in\mathbb{R}$, which maps $\mathbb{R}$ onto
$\mathbb{R}$ and therefore achieves at some point the value $v_{min}^{max}$.
\end{proof}
\begin{remark}
The lower bound estimate of $X_{min}(x_0)$ in terms of $\eta(x_0)$ is
much sharper for small values of $x_0$ than the upper bounds, see
\S~\ref{sect:numerics} for a comparison with the numerical results.
\end{remark}
\begin{definition}
 \label{Def:2.2}
Let $(a,b)$ with $a<0$ be the interval of existence of the solution,
$y(x)$, of Equation~(\ref{p1}) with initial data:
$y(x_0)=y_0$ and $y'(x_0)=y_1$ at $x_0\le0$.
Define functions:
$$
X(x_0)\equiv\underset{y_0,y_1\in\mathbb R}{\inf}\{a(x_0,y_0,y_1)\},\quad
X_-(x_0)\equiv\underset{y_0\in\mathbb R,y_1\leq0}{\inf}\{a(x_0,y_0,y_1)\}
$$
which map $(-\infty,0]$ into $[-\infty,0)$.
\end{definition}
\begin{lemma}
 \label{Lem:2.2}
Let $x_0\leq0$ and $C$ and $\eta(x_0)$ be the same as in
Corollary~{\rm\ref{Cor:X_min_asymptotics}}. Then
\begin{equation}
 \label{ineq:X_-eta_est}
X_-(x_0)<-\eta(x_0)^4=x_0-\frac C{\eta(x_0)}.
\end{equation}
For any $x_0\leq0$ the function $X_-(x_0)$ is finite and for $x_0<0$,
\begin{equation}
 \label{eq:X_-_double_est}
x_0-\frac C{|x_0|^{1/4}}<X_-(x_0)<x_0-
\frac C{\big||x_0|+C/|x_0|^{1/4}\big|^{1/4}}.
\end{equation}
In particular, as $x_0\to-\infty$,
$$
X_-(x_0)=x_0-\frac C{|x_0|^{1/4}}+
\mathcal{O}\left(\frac1{|x_0|^{3/2}}\right),
$$
where
$0<{\mathcal O}\big(\frac1{|x_0|^{3/2}}\big)<\frac{C^2}{4|x_0|^{3/2}}$ for
$|x_0|>C^{4/5}$.
\end{lemma}
\begin{proof}
The proof is analogous to that of Lemma~\ref{Lem:2.1} and
Corollary~\ref{Cor:X_min_asymptotics}: we exclude
from consideration the point $x_0=0$, as done in the first paragraph
of the proof of Lemma~\ref{Lem:2.1}, and switch to notation (\ref{p1t}).
As $u'(t)>0$ now, instead of Equation (\ref{eq:u'}), we get:
\[
u'(t)=\sqrt{4u^3+2N(u_0,u_1)+2\int_{t_0}^tsu'(s)ds}>
\sqrt{4(u^3-u_0^3)+2t_0(u-u_0)}.
\]
This results in the following upper bound for the position of the nearest
(to the right of $t_0$) pole, $t_p$, of the solution $u(t)$,
\begin{equation}
 \label{ineq:tp-v-est}
t_p-t_0<\frac{I(+\infty,v_0)}{t_0^{1/4}},
\end{equation}
where $I(+\infty,v_0)$ is defined by Equation~(\ref{eq:t-v-est-second})
with the change $v_0\mapsto+\infty$ and $v_{min}\mapsto v_0$. The
left Inequality~(\ref{eq:X_-_double_est}) follows from
Inequality~(\ref{ineq:tp-v-est}).
The right Inequality~(\ref{eq:X_-_double_est}) is a consequence of
Inequality~(\ref{ineq:X_-eta_est}). The latter follows from the estimates
analogous to the ones in the first paragraph of the proof of
Corollary~\ref{Cor:X_min_asymptotics} (with $t_{min}\to t_p$). Thus we
arrive at the estimate analogous to (\ref{ineq:tmin_below}),
$t_p^{1/4}(t_p-t_0)>J$, where $J$ is an integral like $I(+\infty,v_0)$,
but with an extra term $+u_1^2/t^{3/2}_p$ under the square root. This term
disappears under the passage to the supremum.

The proof of existence of the solution with $y_0=v_{min}^{max}\sqrt{t_p}$
is even easier here than in Corollary~\ref{Cor:X_min_asymptotics} as it
does not require references to any specific results: the pole
$t_p=|x_p|>|x_0|\geq0$ is a continuous function of the
initial data, $y_0$ and $y_1=-u_1$ at $x_0$, $t_p=t_p(x_0,y_0,y_1)$.
Moreover, as established above it is bounded above. Thus,
$y_0/|x_p(x_0,y_0,0)|$ is a continuous function of $y_0$ which maps
$\mathbb{R}$ onto $\mathbb{R}$.
\end{proof}
\begin{remark}
 \label{Rem:3}
Concerning the behaviour of $X_-(x_0)$ for small $x_0$, one can make the
same comments and estimates as in Remark \ref{Rem:X_min_0}. It is
interesting that for two different functions, $X_{min}(x_0)$ and $X_-(x_0)$,
we have exactly the same bounds and asymptotics; see \S~\ref{sect:numerics}
for a numerical comparison.
\end{remark}
\begin{theorem}
 \label{Thm:1}
Let $x_0\leq0$ and the constant $C$ and function $\eta(\cdot)$ be the same
as in Corollary~{\rm\ref{Cor:X_min_asymptotics}}. Then
\begin{equation}
 \label{ineq:1thm1}
X_-(X_{min}(x_0))\leq X(x_0)<x_0-\frac{2^{4/5}C}{\eta(x_0/2^{4/5})},
\end{equation}
in particular, $X(x_0)$ is finite.
Moreover, for $x_0<0$
\begin{eqnarray}
 \label{ineq:2thm1}
\phantom{AAA}x_0-\frac C{|x_0|^{1/4}}-\frac C{\eta(x_0)}\!\!&\!\!<\!\!&\!\!X_{min}(x_0)-
\frac C{|X_{min}(x_0)|^{1/4}}<X(x_0)\\
&\!\!<\!\!&\!\!x_0-\frac{2^{4/5}C}{|X_{min}(x_0/2^{4/5})|^{1/4}}\!<\!
x_0-\frac{2C}{\big||x_0|+2C/|x_0|^{1/4}\big|^{1/4}},
\label{ineq:3thm1}
\end{eqnarray}
in particular, as $x_0\to-\infty$,
\begin{equation}
 \label{eq:X_asymptotics}
X(x_0)=x_0-\frac{2C}{|x_0|^{1/4}}+
\mathcal{O}\left(\frac1{|x_0|^{3/2}}\right),
\end{equation}
where
$0<{\mathcal O}\big(\frac1{|x_0|^{3/2}}\big)<\frac{C^2}{|x_0|^{3/2}}$
for $|x_0|>(2C)^{4/5}$.
\end{theorem}
\begin{proof}
The first Inequality~(\ref{ineq:1thm1}) is an evident consequence of
Definitions~\ref{Def:2.1} and \ref{Def:2.2}. To prove the second
Inequality~(\ref{ineq:1thm1}), we have to combine corresponding
estimates for the upper bounds from Corollary~\ref{Cor:X_min_asymptotics},
$$
t_{min}-t_0>I(+\infty,v_{min1})t_{min}^{-1/4},
$$
with $v_{min1}=u_{min}/\sqrt{t_{min}}$, and  Lemma~\ref{Lem:2.2},
$$
t_p-t_{min}>I(+\infty,v_{min2})t_p^{-1/4},
$$ where $v_{min2}=u_{min}/\sqrt{t_p}$ and we have assumed that $u_1=0$.
Since $t_p>t_{min}$, we can never get for the same solution
$v_{min1}=v_{min2}=v_{min}^{max}$. Therefore, we make the first estimate
above a bit more rough, i.e., in the original integral (\ref{eq:u-int}) we put
$t_p$ instead of $t(\hat u)$. This formally results in changing
$t_{min}^{-1/4}\to t_p^{-1/4}$ and $v_{min1}\to v_{min2}$ in the first
estimate. Then summing the above estimates, after the same procedure of
passing to the supremum as described in Lemma~\ref{Lem:2.2}, we arrive at
the inequality, $|X(x_0)|^{1/4}(|X(x_0)|-|x_0|)>2C$.

The first Inequality~(\ref{ineq:2thm1}) follows from
Corollary~\ref{Cor:X_min_asymptotics}. The second
Inequality~(\ref{ineq:2thm1}) is the consequence of the first
Inequality~(\ref{ineq:1thm1}) and Lemma~\ref{Lem:2.2}. The first
Inequality~(\ref{ineq:3thm1}) follows from the second
Inequality~(\ref{ineq:1thm1}) and Corollary~\ref{Cor:X_min_asymptotics},
finally, the last Inequality~(\ref{ineq:3thm1}) also results from
Corollary~\ref{Cor:X_min_asymptotics}.
\end{proof}
\begin{remark}
The reader can find numerical values of $X(0)$ and $X(-1)$ in
\S~\ref{sect:numerics}.
\end{remark}

\section{The Functions ${\mathcal X}_{min}(x_0, y_l)$ and $X_{min}(x_0)$}
 \label{sect:3}
\begin{proposition}
 \label{Prop:3.1}
Given $y_0$ and $y_l$ with $y_0> y_l$, there exists only one solution $y(x)$
such that $y(x_0)=y_0$ and $y_{min}=y_l$, where $y_{min}$ is the minimum of
$y(x)$ achieved to the left of $x_0\le 0$.
\end{proposition}
\begin{proof}
It is convenient to use the notation $t=-x$, $u(t)=y(x)$, with
$t_0=-x_0>0$, $u(t_0)=y_0\equiv u_0$, $u'(t_0)=-y'(x_0)\equiv u_1$ and
$y_l=u_l$. See Equation (\ref{p1t}).

Consider a straight line $L:=\{(t,u): u=u_l\}$. The graph of the solution
with initial data
$u(t_0)=u_0$, $u'(t_0)=0$ does not intersect $L$. On the other hand, the
solutions with initial data $u(t_0)=u_0$, $u'(t_0)=-\mu<0$ all intersect
$L$ for large enough $\mu$. This follows from the inequality
\begin{equation}
 \label{eq:Prop:3.1}
-4{u_{min}}^3+2t_{min}(u_0-u_{min})>2N(u_0,\mu),
\end{equation}
which  can be obtained from Equation (\ref{eq:umin}), and the fact that
$2N=\mu^2-4u_0^3\to+\infty$ as $\mu\to+\infty$. By continuity, it follows
that there exists a solution tangent to $L$.

This solution is unique. Actually, suppose
there are two such solutions $u_1(t)$, $u_2(t)$ with
$u_1(t_0)=u_2(t_0)=u_0$, ${u_1}'(t_0)<{u_2}'(t_0)<0$ and the same minimum
$u_{1,min}=u_{2,min}=u_l$. Since their derivatives are non-zero before
their respective minima, $u_1(t)$, $u_2(t)$ are monotonic and we can
define the inverse functions $t_1(u)$, $t_2(u)$ on $[u_l,u_0]$. In some
neighbourhood of $(t_0,u_0)$, $t_1(u)<t_2(u)$. Therefore,
${u_1}''\bigl(t_1(u)\bigr)<{u_2}''\bigl(t_2(u)\bigr)$. Integrating
the last inequality, we get
\begin{equation}
 \label{eq:u_ineq}
0<{{u'}_1}^2(t_0)-{{u'}_2}^2(t_0)<{{u'}_1}^2\bigl(t_1(u)\bigr)
-{{u'}_2}^2\bigl(t_2(u)\bigr).
\end{equation}
{}From this inequality follows that the graphs of $t_1(u)$, $t_2(u)$
cannot cross. Therefore, Inequality (\ref{eq:u_ineq}) holds on the whole
interval $[u_l,u_0]$.
Substituting $u=u_l$ into (\ref{eq:u_ineq}), we get a contradiction, since
the right side is zero while the left side is positive.
\end{proof}
\begin{lemma}
 \label{Lem:3.1}
Let $y_1(x)$ and $y_2(x)$ be solutions of Equation {\rm(\ref{p1})} whose
graphs intersect at $x_0$:
$$
y_1(x_0)=y_2(x_0)=y_0,\quad y'_2(x_0)>y'_1(x_0)\geq0.
$$
Denote by $b_1$ and $b_2$ the right bounds of the interval of existence of
the solutions $y_1(x)$ and $y_2(x)$, respectively. Then $b_1>b_2$ and
$y_2(x)>y_1(x)$ for $x\in(x_0,b_2)$.
\end{lemma}
\begin{proof}
In some right neighbourhood of $x_0$, $y'_2(x)>y'_1(x)>0$. Therefore,
the functions $x_1(y)$ and $x_2(y)$
inverse to the solutions $y_1(x)$ and $y_2(x)$, respectively,
are defined in
some right neighbourhood of $y_0$. Moreover, there exists
$\hat y>y_0$ such that for $y\in[y_0,\hat y)$
\begin{equation}
 \label{eq:cond1}
x'_1(y)>x'_2(y).
\end{equation}
Thus, $0>x_1(y)>x_2(y)$ for $y\in(y_0,\hat y)$. Denote
\[\hat y_0=\sup\{y:\;{\rm Condition (\ref{eq:cond1})\ holds}\}.\]
Suppose, $\hat y_0$ is finite. Then $x'_1(\hat y_0)=x'_2(\hat y_0)$ and
$y_2''(x_2(y))\geq y_1''(x_1(y))$ for $y\in[y_0,\hat y_0]$. Integrating
the last inequality from $y_0$ to $y\leq\hat y_0$ we get
$$
{y'_2}^2(x_2(y))-{y'_1}^2(x_1(y))\geq
{y'_2}^2(x_2(y_0))-{y'_1}^2(x_1(y_0))>0,
$$
implying ${y'_2}(x_2(y))>{y'_1}(x_1(y))>0$. Therefore we deduce that
Inequality (\ref{eq:cond1}) actually holds for any $y\leq\hat y_0$.
But for $y=\hat y_0$ we get a contradiction. Consequently, Inequality
(\ref{eq:cond1}) holds for all $y\in[y_0,+\infty)$. Integrating it
now from $y_0+\epsilon$ to $y>y_0+\epsilon$, where $\epsilon>0$ is some
fixed number, we have
$$
x_1(y)-x_2(y)>x_1(y_0+\epsilon)-x_2(y_0+\epsilon)\equiv\varkappa>0.
$$
This inequality means that there are no intersections of the graphs for
any $y>y_0$. Finally, taking the limit $y\to+\infty$ in this inequality
we get $b_1-b_2\geq\varkappa>0$.
\end{proof}
\begin{proposition}
 \label{Prop:3.2}
Given $y_0$, $\tilde y_0$, and $y_l$ with $\tilde y_0>y_0>y_l$,
consider solutions $y(x)$, $\tilde y(x)$ such that $y(x_0)=y_0$,
$\tilde y(x_0)=\tilde y_0$ and $y_{min}=\tilde y_{min}= y_l$.
Then $x_{min}>\tilde x_{min}$ and $\tilde y'(x_0)>y'(x_0)>0$. Moreover,
$\tilde y(x)>y(x)$ for $x\in[x_{min},x_0]$.
\end{proposition}
\begin{proof}
For $\mu\in{\mathbb R}$, consider the solutions
$\tilde y_\mu(x)$ with the same initial value $\tilde y_0$
and $\mu=\tilde y_\mu'(x_0)$. From the proof of the previous proposition,
it is clear that for sufficiently large $\mu>M>0$, these solutions cross
the graph of $y(x)$, while for all $\mu\le 0$ they do not.

Continuously increasing $\mu$ from $0$ to $M$, the first intersection of
the graphs $\tilde y_\mu(x)$ with $y(x)$ cannot occur for
$x\in[x_{min},x_0]$. Otherwise, there would exist $\mu_1$, such that
the graphs of $\tilde y_{\mu_1}(x)$ and $y(x)$ are tangent at some point
in this interval. This contradicts the uniqueness of solution of the
Cauchy problem at that point.

The first ``intersection'' occurs for some $\mu=\mu_2$ at $x=x_p$, the
pole of both solutions $y(x)$ and $\tilde y_{\mu_2}(x)$. Actually, the
solutions are asymptotic to the vertical line $x=x_p$ without intersection.%asymptote->asymptotic
By increasing $\mu$ from $\mu_2$, we observe that a finite intersection
point $A$ with coordinates $(A_x,A_y)$ appears. This point $A=A(\mu)$
moves continuously on the graphs of both $y(x)$ and $\tilde y_{\mu}(x)$
as $\mu$ increases. For $\mu>0$, the graph of $\tilde y_{\mu}(x)$ has a
minimum $\bigl\{x_{min}(\mu), y_{min}(\mu)\bigr\}$. For large $\mu$,
$A_x\in(x_{min}, x_0)$. Thus, for some value of $\mu=\mu_{min}$, the point
$A$ coincides with the minimum point of the graph of $y(x)$. For this
value of $\mu=\mu_{min}$, the minimum point of the graph of
$\tilde y_{\mu}(x)$ is lower and to the left of the minimum point of the
graph of $y(x)$. Consequently, for some $\mu<\mu_{min}$, $A(\mu)$ crossed
the minimum point of the graph of $\tilde y_{\mu}(x)$. Such a crossing
may occur more than once. We choose the largest such value of $\mu$ and
denote it by $\tilde\mu_{min}$. Therefore, for all
$\mu\in\bigl(\tilde\mu_{min},\mu_{min}\bigr)$,
$\tilde x_{min}(\mu)<A_x<x_{min}$
and, for some $\mu$ in this interval, the minimum value of
$\tilde y_\mu(x)$ equals $y_l$. This proves the existence of the desired
solution $\tilde y(x)$ with $\tilde x_{min}<x_{min}$. It is unique by
Proposition \ref{Prop:3.1}.

Denote by $\tilde x$ the point where the straight line $y=y_0$ first
crosses the graph of $\tilde y(x)$. Arguments similar to that of the
previous proposition show that ${\tilde y}'(\tilde x)>y'(x_0)$;
otherwise, $\tilde y_{min}>y_l$.
It is clear that ${\tilde y}'(x_0)>{\tilde y}'(\tilde x)$.
\end{proof}
\begin{remark}
There is another proof of Proposition \ref{Prop:3.2} in the spirit
of Lemma \ref{Lem:3.1}.
\end{remark}
\begin{corollary}
 \label{Cor:3.1}
Consider the solution with initial value $y(x_0)=y_0$ and minimum value
$y_l$ achieved to the left of $x_0$.
The initial slope $y_1=f(x_0;y_0,y_l)$ at $x_0$ is a continuous function
defined for $x_0\leq0$ and $y_0\geq y_l\in\mathbb R$ with $f(x_0;y_0,y_0)=0$.
Moreover, $f(x_0;y_0,y_l)$ is a smooth function of all its variables
for $x\leq0$ and $y_l<y_0$ with $\partial_{y_l}f(x_0;y_0,y_l)<0$.

For fixed $x_0$ and $y_l$ it is a monotonically increasing
function with the asymptotic behaviour
\[
f(x_0;y_0,y_l)\underset{y_0\to+\infty}{=}2y_0\sqrt{y_0}+
{\mathcal O}\left(\frac1{\sqrt{y_0}}\right).
\]
For fixed $x_0$ and $y_0$, it is monotonically decreasing with the
asymptotics
\[
f(x_0;y_0,y_l)\underset{y_l\to-\infty}{=}-2y_l\sqrt{-y_l}+
{\mathcal O}\left(\frac1{\sqrt{-y_l}}\right).
\]
For fixed $y_0$ and $y_l$, $f(x_0;y_0,y_l)$ is a
monotonically decreasing function of $x_0$ with the
asymptotics
$$
f(x_0;y_0,y_l)\underset{x_0\to-\infty}{=}\sqrt{2|x_0|(y_0-y_l)}
\left(1+\frac{y_0^2+y_l^2+y_0y_l}{|x_0|}+
{\mathcal O}\left(\frac1{|x_0|^{5/4}}\right)\right).
$$
\end{corollary}
\begin{proof}
The continuity, actually  even the smoothness of the inverse function
$y_l=g(x_0;y_0,y_1)$, is evident as $y_l=y(x_{min};x_0,y_0,y_1)$,
where $y(x;x_0,y_0,y_1)$ is a solution of Equation~\ref{p1} with
initial data $y_0$ and $y_1$ given at $x_0$, and
$x_{min}=x_{min}(x_0,y_0,y_1)$ is a smooth function by the implicit function
theorem as the solution of the equation $y'(x_{min};x_0,y_0,y_1)=0$, with
$y''(x_{min};x_0,y_0,y_1)=6y_l^2-x_{min}>0$. Continuity of $f(x_0;y_0,y_l)$
follows from continuity of $g(x_0;y_0,y_1)$ by standard arguments,
we recall them below in the particular case of continuity on $y_l$.
Consider a convergent sequence $\{y_{ln}\}$ of values of $y_l$ and
let $y_{l0}=\lim y_{ln}$. Define $y_{1n}=f(x_0;y_0,y_{ln})$, then the
sequence $y_{1n}$ is bounded, because unbounded slopes, $y_1\to-\infty$,
correspond to unbounded ``levels'', $y_l\to-\infty$. If we suppose that
$\lim y_{1n}$ does not exist, then there exist two subsequences, $y_{1n_k}$
and $y_{1n_m}$, $k,m=1,2,\ldots$, convergent to $y_{10}$ and $y_{11}$,
respectively, where $y_{10}\neq y_{11}$.
However, using the continuity of $g$ we obtain:
$$
g(y_{10})=g(\lim y_{1n_k})=\lim y_{ln_k}=y_{l0},\quad
g(y_{11})=g(\lim y_{1n_m})=\lim y_{ln_m}=y_{l0}.
$$
This gives a contradiction, since it implies that $y_{10}=y_{11}$.
Hereafter we suppress the dependence of $f$ and $g$ on $x_0$ and $y_0$.
Now, again by continuity of $g$, we find:
$$
\lim f(y_{ln})=\lim y_{1n}=f(g(\lim y_{1n}))=
f(\lim g(y_{1n}))=f(\lim y_{ln})=f(y_{l0}).
$$
To prove smoothness of the function $f$, we rewrite
Equation~(\ref{eq:umin}) in original variables, $y$ and $x$, and
differentiate the latter with respect to
$y_1$ recalling that $y_l=g(x_0;y_0,y_1)$:
\begin{equation}
 \label{eq:diff_yl}
0=y_l'(6y_l^2-x_m)+y_1-\int\limits_{y_0}^{y_l}x'(u;x_0,y_0,y_1)du,
\end{equation}
where the prime is the derivative with respect to $y_1$ and the minimum
$x_m=x(y_l;x_0,y_0,y_1)$. A key observation now is that $x(u;x_0,y_0,y_1)$
is a monotonically increasing function of $y_1$, for fixed values of the
other variables (cf. the last paragraph of the proof of
Proposition~\ref{Prop:3.1}), so that the integrand, $x'$, is nonnegative.
Therefore, the integral in Equation~(\ref{eq:diff_yl}) is convergent, as
$y_l'$ exists, and (strictly) negative since we assume that $y_l<y_0$.
At the same time $y_1>0$ under the same condition, $y_l<y_0$, and the
factor multiplying $y_l'$ is positive as $-x_m>-x_0\geq0$. Thus
Equation~(\ref{eq:diff_yl}), together with continuity of $f$, implies that
$y_l'$ is continuous and strictly negative for $x_0\leq0$ and $y_l<y_0$.
So, the derivative of the function inverse to $y_l$, i.e., the function
$f$, has the same properties. The proof of smoothness with respect to
$x_0$ and $y_0$ is analogous with the help of the observation that
$\partial_{x_0}x\geq0$ and $\partial_{y_0}x\leq0$, so that the
corresponding integrals exist and are strictly positive or, respectively,
negative.

Rewriting Equation (\ref{eq:umin}) in terms of $y(x)$, we get the estimate
\begin{equation}
 \label{eq:cor31}
2|x_{min}|(y_0-y_l)-4y_l^3>y_1^2-4y_0^3>2|x_{0}|(y_0-y_l)-4y_l^3
\end{equation}
The first two asymptotics follow immediately, whilst the proof of the
last one requires a reference to estimate (\ref{eq:est-xmin}).
\end{proof}
\begin{remark}
 \label{Rem:notation}
It follows from Corollary~\ref{Cor:3.1} that instead of a function, say,
$h(x_0,\!y_0,\!y_1)$ smoothly depending on initial data $x_0$, $y_0$, and
$y_1$, we can always consider a function
$\hat h(x_0,y_0,y_l)=h(x_0,y_0,f(x_0;y_0,y_l))$, which depends smoothly
on $x_0$, $y_0$, and $y_l$. Where it does not cause confusion, we admit
an abuse of notation and denote both functions, $h$ and $\hat h$, by the
same symbol, $h$.
\end{remark}
\begin{remark}
 \label{Rem:symmetry}
Corollary \ref{Cor:3.1} shows that the functions
$f(x_0;y_0,y_l)$ and $f(x_0;-y_l,-y_0)$ have similar monotonicity
properties. Moreover, the difference
$\Delta(x_0;y_0,y_l)=f(x_0;y_0,y_l)-f(x_0;-y_l,-y_0)=o(1)$ when
one of the variables $|x_0|$, $|y_0|$ or $|y_l|$ is large.
This estimate also holds when two or all three of the variables
are large, as can be deduced from the corresponding asymptotics of
$f(x_0;y_0,y_l)$. It is straightforward to obtain these asymptotics
from Inequality~(\ref{eq:cor31}) and estimate (\ref{eq:est-xmin}).
They, clearly, inherit from (\ref{eq:cor31}) an
invariance under the change of variables $y_0,y_l\to-y_l,-y_0$.
Therefore the function $\Delta(x_0;y_0,y_l)$ is bounded in its
domain of definition,
${\mathcal D}:=\{x_0\leq0,y_l\leq y_0\}\subset\mathbb R^3$.
Note also that $\Delta(x_0;y_0,\pm y_0)=0$.
Some other numerical values of $\Delta$ are given in \S~\ref{sect:numerics}.

One can define also the function $f^+(x_0;y_0,y_l)$ as initial slope
of the solution with the initial value $y(x_0)=y_0$ and
minimum value $y_l\leq y_0$ achieved to the right of $x_0$. It is clear
that $f^+<0$ and $|f^+(x_0;y_0,y_l)|<f(x_0;y_0,y_l)$.
Actually, the function $|f^+|$ has properties similar to those of $f$.
Moreover, using $f^+$ one can make a continuous prolongation of $f$ from
${\mathcal D}$ to $\mathbb R_-\times\mathbb R^2$, by putting
$f(x_0;y_0,y_l):=f^+(x_0;y_l,y_0)$ for $y_l\geq y_0$. It would be interesting
to know whether this continuation is smooth as $x'$ in
Equation~(\ref{eq:diff_yl}) is unbounded as $y_l\to y_0$, and thus a
finite limit of $y_l'$ is possible.
\end{remark}
\begin{definition}
 \label{Def:3.1}
Given $x_0\le0$, $y_0>y_l$, consider the solution $y(x)=y(x;x_0,y_0,y_l)$
with initial value $y(x_0)=y_0$ and minimum value $y_l$ at
$x_{min}=x_{min}(x_0,y_0,y_l)<x_0$. Define
\[{\mathcal X}_{min}(x_0,y_l)=\underset{y_0>y_l}{\inf}\,
\left\{x_{min}(x_0,y_0,y_l)\right\}\]
\end{definition}
\begin{proposition}
 \label{Prop:3.3}
${\mathcal X}_{min}(x_0, y_l)$ is finite. The solution $y_m(x)$ with
initial data
\[y_m\bigl({\mathcal X}_{min}(x_0, y_l)\bigr)=y_l,\quad
y_m'\bigl({\mathcal X}_{min}(x_0, y_l)\bigr)=0\]
has a pole at $x_0$.
\end{proposition}
\begin{proof}
Definitions \ref{Def:2.1} and \ref{Def:3.1} imply
\[X_{min}(x_0)\leq{\mathcal X}_{min}(x_0, y_l)<x_0.\]
The solution $y_m(x)$ cannot be regular at $x_0$ as it is a limit of
solutions with initial value $y_0\to+\infty$, $y_0'\to+\infty$ at $x_0$.
(See Proposition \ref{Prop:3.2} and Corollary \ref{Cor:3.1}.) Therefore,
it has a pole $x_p\leq x_0$. Suppose $x_p<x_0$.

For any given $\epsilon>0$ and a sequence $x_{n}\nearrow x_p$,
one finds $y_{0n}\nearrow+\infty$ such that for all values of
$y_0\geq y_{0n}$, we have
$$
\left|y_m(x)-y(x;x_0,y_0,y_l)\right|\le\epsilon,\quad
\left|y'_m(x)-y'(x;x_0,y_0,y_l)\right|\le\epsilon,
$$
for $x\in[{\mathcal X}_{min}(x_0,y_l),x_n]$. In particular, this means that
\[
y(x_{n};x_0,y_{0n},y_l)\underset{y_0\to+\infty}=L_n\to+\infty, \quad
y'(x_{n};x_0,y_{0n},y_l)\underset{y_0\to+\infty}=M_n\to+\infty.
\]
Arguments analogous to that of Lemma \ref{Lem:2.2} show that
any solution with initial values $y(x_{n})=L_n$, $y'(x_{n})\ge0$ has
a pole $x_{pn}>x_{n}$, whose location can be estimated by
$$
\frac1{2\sqrt{L_n}}\int_1^{+\infty}\frac{dv}{\sqrt{v^3-1}}>x_{pn}-x_n.
$$
Therefore, for all sufficiently large $L_n$, $x_{min}<x_{pn}<x_0$. This
contradicts the fact that $y(x;x_0,y_0,y_l)$ is finite on $[x_{min},x_0]$.
\end{proof}
\begin{lemma}
 \label{Lem:3.2}
Given $y_l\in{\mathbb R}$ there is only one solution that has minimum
value $y_l$ with a pole at $x_0$ and an interval of existence to the
left of $x_0$.
\end{lemma}
\begin{proof}
Existence of the solution $y(x;x_0,y_l)$ follows from Proposition
\ref{Prop:3.3}. Any solution with a pole at $x_0$ has the following
convergent Laurent expansion:
\begin{equation}
 \label{laurent}
y(x)=\frac{1}{(x-x_0)^2}+\frac{x_0}{10}(x-x_0)^2+\frac{1}{6}(x-x_0)^3
+c(x-x_0)^4+\ldots
\end{equation}
where $c\in{\mathbb R}$ is a parameter characterising the solution.
Suppose $y_1(x)$, $y_2(x)$ are solutions corresponding to $c_1$ and
$c_2$ with $c_1>c_2$. Consider the difference
\begin{eqnarray}
 \label{Eqn:1}
y_1-y_2&\underset{x\to x_0}{=}&(c_1-c_2)(x-x_0)^4+
{\mathcal O}\bigl((x-x_0)^5\bigr),\\
y_1'-y_2'&\underset{x\to x_0}{=}&4(c_1-c_2)(x-x_0)^3+
{\mathcal O}\bigl((x-x_0)^4\bigr).
\end{eqnarray}
{}From this it follows that
\begin{equation}
 \label{ineq:lau}
y_1(x)>y_2(x),\quad 0<y_1'(x)<y_2'(x)
\end{equation}
in some (small) left neighbourhood of $x_0$. Consider a horizontal line
which crosses the graphs of $y_1(x)$, $y_2(x)$ at points $x_1$, $x_2$ in
the neighbourhood of $x_0$ where Inequalities (\ref{ineq:lau}) hold.
We have
\[
y_1'(x_1)<y_2'(x_2).
\]
Repeating almost exactly the proof of Proposition \ref{Prop:3.1}, we see
that the minimum values of the solutions satisfy
${y_1}_{min}>{y_2}_{min}$.
\end{proof}
\begin{theorem}
 \label{Thm:3.1}
Any solution with an interval of existence that is contained in the
non-positive semi-axis can be uniquely characterised by the position of
its right pole $x_0$ and minimum value $y_l$. Moreover, it achieves its
minimum value at ${\mathcal X}_{min}(x_0,y_l)$.
\end{theorem}
\begin{proof}
The proof follows from the above propositions and Lemma \ref{Lem:3.2}.
\end{proof}
\begin{remark}
 \label{Rem:correspondence}
{}From Theorem \ref{Thm:3.1} we see that there is a natural one-to-one
correspondence between solutions of Equation (\ref{p1}) and the functions
${\mathcal X}_{min}(x_0,y_l)$. The solution $y(x; x_0, y_l)$ that
corresponds to the function ${\mathcal X}_{min}(x_0,y_l)$ has the following
extremal property: its minimum ${\mathcal X}_{min}(x_0,y_l)<x_{min}$
where $x_{min}$ is the minimum of any solution regular at $x_0$ with
minimum value $y_l$. Moreover, the properties of the solutions can be
reformulated as monotonicity properties of the functions
${\mathcal X}_{min}(x_0,y_l)$.
\end{remark}
\begin{corollary}
 \label{Cor:3.2}
For any $y_l$, if $x_1<x_2<0$ then
\begin{eqnarray*}
{\mathcal X}_{min}(x_1,y_l)&<&{\mathcal X}_{min}(x_2,y_l)
\end{eqnarray*}
\end{corollary}
\begin{proof}
If ${\mathcal X}_{min}(x_1,y_l)={\mathcal X}_{min}(x_2,y_l)$
then the two corresponding solutions (see Remark \ref{Rem:correspondence})
would coincide.
Suppose we have ${\mathcal X}_{min}(x_1,y_l)>{\mathcal X}_{min}(x_2,y_l)$.
Then there is a solution regular at $x_1$, with
$x_{min}<{\mathcal X}_{min}(x_1,y_l)$. This contradicts the definition
of ${\mathcal X}_{min}(x_0,y_l)$.
\end{proof}
\begin{proposition}
 \label{Prop:analyticity}
For any pair of numbers $(x_0,c_0)\in{\mathbb C^2}$ consider a Laurent
expansion~{\rm(\ref{laurent})} (with $c\mapsto c_0$) as a germ defining
a complex solution $y(x)$ by analytic continuation on $x$. Denote
such a complex solution by $y(x;x_0,c_0)$. It is a meromorphic function of
$x$ which, for any choice of $(x_0,c_0)$, has an infinite discrete set of
second order poles, ${\mathcal P}={\mathcal P}(x_0,c_0)$. For any
$x\in\mathbb{C}\setminus\mathcal P$ the function $y(x;x_0,c_0)$ is
an analytic (locally holomorphic) function of the parameters $x_0$
and $c_0$.

Fix any pair $(\hat x_0,\hat c_0)\in\mathbb{C}^2$. If $\hat x_1$ is any
other pole of $y(x;\hat x_0,\hat c_0)$ and $\hat c_1$ is the corresponding
pole parameter of expansion~{\rm(\ref{laurent})}, then there exist a unique
pair of functions $x_1=x_1(x_0,c_0)$, $c_1=c_1(x_0,c_0)$ which is
{\rm(1)} holomorphic at $(\hat x_0,\hat c_0)$;
{\rm(2)} $\hat x_1=x_1(\hat x_0,\hat c_0)$,
$\hat c_1=c_1(\hat x_0,\hat c_0)$; and
{\rm(3)} $x_0$, $x_1$ and $c_0$, $c_1$ are respectively the poles and
corresponding pole parameters of the complex solution $y(x;x_0,c_0)$.
\end{proposition}
\begin{proof}
Let $r_k$, $k=0,1$ be the radii of convergence of the Laurent expansions
at $x_k$ and $\epsilon>0$ be chosen sufficiently small, $\epsilon<r_k$.
Then in the annuli $R_k:=\{x:\epsilon<|x-x_k|<r_k\}$, the complex
solutions $y(x;x_k,c_k)$ are (locally) holomorphic functions of variables
$x_k$ and $c_k$ since the coefficients of the Laurent
expansions~{\rm(\ref{laurent})} are polynomials in $x_k$ and $c_k$.
Take any points $\tilde x_k\in R_k$ and consider the analytic functions of
$x_k$ and $c_k$: $y_k=y(\tilde x_k;x_k,c_k)$ and
$y_{1k}=y'(\tilde x_k;x_k,c_k)$, as the initial data for the complex
solution $y(x;x_k,c_k)$, where now $x$ is an arbitrary point in
${\mathbb C}\setminus{\mathcal P}$. For each $k$, $y(x;x_k,c_k)$ is
an analytic function of $x_k$ and $c_k$ as a composition of analytic
functions $x_k,c_k\mapsto y_k,y_{1k}$ and $y_k,y_{1k}\mapsto y(x;x_k,c_k)$.

Consider now the system of equations which guarantees that the complex
solutions $y(x;x_k,c_k)$ have the same ``initial data'' at $x=\tilde x_1$:
\begin{equation}
 \label{eq:system}
y(\tilde x_1;x_0,c_0)=y(\tilde x_1;x_1,c_1),\qquad
y'(\tilde x_1;x_0,c_0)=y'(\tilde x_1;x_1,c_1),
\end{equation}
where the left (respectively, right) sides of the equations are considered
as holomorphic functions of $x_0$ and $c_0$ (respectively, $x_1$ and $c_1)$.
We claim that the unique solvability of this system in neighbourhoods of
$(\hat x_k,\hat c_k)$ and the existence of the corresponding derivatives
follow from the implicit function theorem. In fact, consider the Jacobian,
$$
J(\tilde x_1)=\left(
\begin{matrix}
\partial_{x_1}y(\tilde x_1;x_1,c_1)&\partial_{c_1}y(\tilde x_1;x_1,c_1)\\
\partial_{x_1}y'(\tilde x_1;x_1,c_1)&\partial_{c_1}y'(\tilde x_1;x_1,c_1)
\end{matrix}\right).
$$
This Jacobian is independent of $\tilde x_1$ as it coincides with the
Wronskian of two complex solutions $\partial_{x_1}y(\tilde x_1;x_1,c_1)$ and
$\partial_{c_1}y(\tilde x_1;x_1,c_1)$ of the linearization of
Equation~(\ref{p1}), $Y''(\tilde x_1)=12y(\tilde x_1;x_1,c_1)Y(\tilde x_1)$.
Direct calculation with the help of Laurent expansion~(\ref{laurent}) gives
$J(\tilde x_1)=14$.
Finally, we remark that the functions $x_1(x_0,c_0)$ and $c_1(x_0,c_0)$
do not actually depend on the ``connection'' point $\tilde x_1$. In fact,
if these functions solve System~(\ref{eq:system}) at some point
$\tilde x_1$, then they solve it at any other point in the annulus $R_1$,
as this means that $y(x;x_0,c_0)$ and $y(x;x_1(x_0,c_0),c_1(x_0,c_0))$ have
the same initial data at $\tilde x_1$ and thus they coincide for all $x$ due
to the uniqueness of solution of the Cauchy initial value problem. Suppose
now that there is another pair of functions $x_1={\bf x}_1(x_0,c_0)$ and
$c_1={\bf c}_1(x_0,c_0)$ with the properties (1)--(3).
Then clearly, these functions solve System~(\ref{eq:system}) and therefore
coincide with $x_1(x_0,c_0)$ and $c_1(x_0,c_0)$, respectively.
\end{proof}
\begin{remark}
To establish most of the qualitative properties of the solutions
we need continuity rather than analyticity of the corresponding functions.
However, in \S\S~\ref{sect:main} and \ref{sect:boundary} we actually use
complex analyticity of the function $x_1=x_1(x_0,c_0)$ with respect to
$c_0$. Below we use the adjective smooth to indicate that corresponding
real function has one continuous derivative; in fact, most of our functions
are real analytic.
\end{remark}
\begin{remark}
A solution with a pole at $x_0$ and interval of existence to the left
is uniquely characterised on the one hand by the parameter $c$ in
Laurent expansion (\ref{laurent}) and on the other hand, by its minimum
$y_l$ achieved at $x_{min}=x_{min}(x_0, c)$.
Therefore, for any given $x_0\leq0$ there exists a bijection
$y_l(x_0,c):{\mathbb R}\to{\mathbb R}$.
\end{remark}
\begin{corollary}
 \label{Cor:3.3}
The function $y_l(x_0,c)$ is a smooth function of both variables.
For any $x_0\leq0$, $y_l(x_0,c)$ is a monotonically growing function of
$c\in\mathbb R$. Moreover, $y_l(x_0,c)\to\pm\infty$ as $c\to\pm\infty$,
respectively. The inverse function $c(x_0,y_l)$ has the same properties.
The function $x_{min}=x_{min}(x_0,c)$ is a smooth function of both variables.
\end{corollary}
\begin{proof}
To denote more explicitly the dependence of the solution on the parameters
$x_0$ and $c$, we use the notation $y(x)=y(x;x_0,c)$. As follows from
Proposition~\ref{Prop:analyticity} both $y(x;x_0,c)$ and $y'_1(x;x_0,c)$,
where the subscript $1$ means that the derivative is taken with respect to
the first variable, are smooth functions of $x_0$ and $c$.

Now $x_{min}(x_0,c)$ is the unique solution of the equation
$y'_1(x;x_0,c)=0$.
By the implicit function theorem, both derivatives $x_{min,k}'(x_0,c)$,
$k=1,2$, exist and are given by
\[
x_{min,k}'(x_0,c) =-\,\frac{y''_{1(k+1)}(x_{min};x_0,c)}
{y''_{11}(x_{min};x_0,c)},
\]
where the denominator on the right is always positive because
$x_{min}<x_0\le 0$.
Clearly, $y_l(x_0,c)=y(x_{min}(x_0,c);x_0,c)$ is also a smooth function of
$x_0$ and $c$. Monotonicity of $y_l(x_0,c)$ with respect to $c$ is proved
in fact in Lemma~\ref{Lem:3.2}.

Let $c\searrow-\infty$. Then $y_l(x_0,c)$ monotonically decreases to a limit
$\hat y(x_0)$. Suppose $\hat y(x_0)>-\infty$.
In this case, $c(x_0,y_l)$ would not be defined in the interval
$(-\infty, \hat y(x_0))$. An analogous argument shows that
$y_l(x_0,c)\nearrow+\infty$ as $c\nearrow+\infty$.

One now proves continuity of $c(x_0,y_l)$ by almost a repetition of
the analogous proof for continuity of the function $y_1=f(x_0,y_0,y_l)$
in Corollary~\ref{Cor:3.1}. Since the function $x_{min}=x_{min}(x_0,c)$
is smooth, as proved above, this continuity of $c(x_0,y_l)$ implies
continuity of $x_{min}(x_0,y_l)$.

To prove smoothness of $c(x_0,y_l)$, one starts with the
initial value problem at the minimum of the solution $y(x;x_0,c)$,
$y(x_{min};x_0,c)=y_l$, $y'_1(x_{min},x_0,c)=0$.
The solution of this initial value problem for Equation~(\ref{p1})
is a smooth function of the parameters $x_{min}$ and $y_l$. So we have
another parameterization of the same solution,
\begin{equation}
 \label{eq:reparametrization}
y(x;x_0,c)=y(x;x_{min},y_l)
\end{equation}
with a slight abuse of notation. From this equation one deduces that
the functions $x_0=x_0(x_{min},y_l)$ and $c=c(x_{min},y_l)$ are
smooth functions of their arguments. To prove this, define the functions
$F_k(x_0,x_{min},y_l)=
\frac1{2\pi i}\oint\frac{y(x;x_{min},y_l)}{(x-x_0)^{1+2k}}dx
$,
$k=1,2$, where integration is taken anti-clockwise along a circle centred
at $x_0$. Then, from Laurent expansion~(\ref{laurent}) one finds the
following representation for the function $c$, $c=F_2(x_0,x_{min},y_l)$,
where $x_0=x_0(x_{min},y_l)$ solves equation $F_1(x_0,x_{min},y_l)=x_0/10$.
Using again Laurent expansion~(\ref{laurent}) we prove that
$\partial_{x_0}F_1=3/6=1/2$, thus the implicit function theorem implies
the smoothness of $x_0(x_m,y_l)$.

Now we differentiate Equation~(\ref{eq:reparametrization}) with respect to
$x$ to get analogous equation for the $x$-derivatives and consider both
equations as a system to determine the functions $x_{min}(x_0,y_l)$ and
$c(x_0,y_l)$. To employ the implicit function theorem we have to check that
the corresponding Jacobian,
$
J=y'_3(x;x_0,c)y''_{12}(x;x_{min},y_l)-
y'_2(x;x_{min},y_l)y''_{13}(x;x_0,c)
$,
is nonzero for all $x_0\leq0$ and real $y_l$. It is straightforward to see
that $J$ is independent of $x$ (see Proposition~\ref{Prop:analyticity}).
Considering it for $x$ in a proper neighbourhood of $x_0$, we find that
$J=-14\partial_{x_{min}}x_0(x_{min},y_l)$. Our goal now is to prove that
$\partial_{x_{min}}x_0(x_{min},y_l)>0$. The fact that
$\partial_{x_{min}}x_0(x_{min},y_l)\geq0$ is easy to deduce from
Proposition~\ref{Prop:3.1}: it is a monotonically growing function. Let
us rewrite Equation~(\ref{eq:u-int}) by turning back to the the original
$x$, $y$ variables, putting now that $t_0=-x_0$ is the pole, so that
$u_0=+\infty$,
\begin{equation}
 \label{eq:int_x0-xmin}
x_0-x_{min}=\int_{y_l}^{+\infty}\frac{dy}\kappa,\qquad
\kappa=\sqrt{4(y^3-y_l^3)-
2\int_{y_l}^yx(\tilde y;x_{min},y_l)d\tilde y}.
\end{equation}
Differentiating this equation with respect to $x_{min}$ we obtain
\begin{equation}
 \label{eq:derivative}
\partial_{x_{min}}x_0-1=\int_{y_l}^{+\infty}\frac{dy}{\kappa^3}
\int_{y_l}^{y}\partial_{x_{min}}x(\tilde y;x_{min},y_l)d\tilde y.
\end{equation}
We note that $\partial_{x_{min}}x(\tilde y;x_{min},y_l)$ exists for
$\tilde y>y_l$ and from Proposition~\ref{Prop:3.1} follows that
$\partial_{x_{min}}x(\tilde y;x_{min},y_l)\geq0$.
The integral~(\ref{eq:derivative}) is improper, however it is
straightforward to see that it converges, obviously at infinity and
at $y=y_l$ in virtue of the following limit:
$
\underset{\tilde y\to y_l}\lim\partial_{x_{min}}x(\tilde y;x_{min},y_l)=1.
$
Thus, $\partial_{x_{min}}x_0(x_{min},y_l)>1$ as the right-hand side of
Equation~(\ref{eq:derivative}) is positive.
\end{proof}
\begin{remark}
 \label{Rem:smoothxmin}
{}From the above corollary it follows that the function
${\mathcal X}_{min}(x_0,y_l)=x_{min}(x_0,c(x_0,y_l))$
is a smooth function of both variables.
\end{remark}
\begin{proposition}
 \label{Prop:3.5}
For any $y_l$ and $x_1<x_2<0$,
\begin{equation}
 \label{ineq:calxmin}
0<{\mathcal X}_{min}(x_2,y_l)-{\mathcal X}_{min}(x_1, y_l)< x_2-x_1.
\end{equation}
\end{proposition}
\begin{proof}
Let $y_1(x)$, $y_2(x)$ be the solutions with poles at $x_1$, $x_2$ and the
same minimum value $y_l$.
The inverse functions $x_1(y)$ with range $[x_{1,min},x_1)$,
and $x_2(y)$ with range $[x_{2,min},x_2)$ exist and are differentiable on
the interval $(y_l,+\infty)$. Of course, by this definition
$x_{i,min}=x_i(y_l)$ for $i=1,2$. Moreover, for $y\in(y_l,+\infty)$, we have
\[
y_1''(x_1(y))>y_2''(x_2(y))\ \Rightarrow\ y_1'(x_1(y))>y_2'(x_2(y))>0
\]
where we have used the fact that $y_i'(x_i(y_l))=0$. This implies
\[x_1'(y)<x_2'(y).\]
Integration from $y_l$ to $+\infty$ yields the right
Inequality~(\ref{ineq:calxmin}).
The left inequality follows from Corollary \ref{Cor:3.2}.
\end{proof}
\begin{proposition}
 \label{Prop:3.6}
For $|y_l|\geq\epsilon>0$
$$
0<x_0-{\mathcal X}_{min}(x_0,y_l)<\frac{I_\nu}{2\sqrt{|y_l|}},
$$
where $\nu=\mathrm{sign}\,\{y_l\}1$ and $I_\nu$ is defined by
Equation~{\rm(\ref{eq:intdelta})}.
\end{proposition}
\begin{proof}
Neglecting in Equation~(\ref{eq:u-int}) the integral term under the
square root, sending the upper limit of the remaining integral
to $+\infty$ and $t\to t_{min}$, making the change of variables,
$u=|y_l|w$ ($u_{min}=y_l$), and using Proposition~\ref{Prop:3.2} and
Definition~\ref{Def:3.1} we arrive at the stated result.
\end{proof}
\begin{proposition}
 \label{Prop:3.7}
For any $x_1<x_2<0$,
\begin{equation}
\label{ineq:xmin}
0<X_{min}(x_2)-X_{min}(x_1)<x_2-x_1.
\end{equation}
\end{proposition}
\begin{proof}
{}From Remark \ref{Rem:smoothxmin}, we have
for any $y_l$ and $x_0<0$, that ${\mathcal X}_{min}(x_0,y_l)$ is a
continuous function of $y_l$. Moreover, as follows from
Proposition~\ref{Prop:3.6},
\[
0<x_0-{\mathcal X}_{min}(x_0,y_l)\to0,\quad{\rm as}\ y_l\to\pm\infty.
\]
Therefore, $\exists$ $\hat y_l(x_0)$ ($\le 0$) such that
\[X_{min}(x_0)={\mathcal X}_{min}(x_0, \hat y_l(x_0)).\]
The left Inequality~(\ref{ineq:xmin}) follows from
\[
X_{min}(x_2)={\mathcal X}_{min}(x_2,\hat y_l(x_2))>
{\mathcal X}_{min}(x_1,
\hat y_l(x_2))\ge X_{min}(x_1)
\]
and the right one can be proved analogously:
\begin{eqnarray*}
x_1-X_{min}(x_1)&=&x_1-{\mathcal X}_{min}(x_1, \hat y_l(x_1))\\
                &<& x_2-{\mathcal X}_{min}(x_2, \hat y_l(x_1))
                \le x_2-X_{min}(x_2)
\end{eqnarray*}
where we have again used Proposition \ref{Prop:3.5}.
\end{proof}
\begin{remark}
 \label{Rem:conjecture}
It is easy to prove that ${\mathcal X}_{min}(x_0,y_l)$ is monotonically
decreasing as $y_l\searrow +0$.
We expect that this monotonic decrease continues until some negative value of
$y_l=\hat y(x_l)$.
After that value, ${\mathcal X}_{min}(x_0,y_l)$ monotonically grows
with $y_l{\searrow} -\infty$.
Therefore, we conjecture that the solution with minimum at $X_{min}(x_0)$
is unique. This uniqueness was not assumed in the proof of the above
proposition.
\end{remark}
\begin{proposition}
 \label{Prop:der_mathcalXmonotonic}
For any $y_0>y_l$ and four points $x_k$, $k=1,2,3,4$, such that
$x_4<x_3<x_1\leq0$, $x_4<x_2<x_1$, and $0<x_1-x_2\leq x_3-x_4$,
consider four solutions $y_k(x;x_k,y_0,y_l)$. Denote their minima
$x_{min}(x_k,y_0,y_l)$. Then,
\begin{equation}
 \label{eq:der_xminmonotonic}
x_{min}(x_1;y_0,y_l)-x_{min}(x_2;y_0,y_l)\leq
x_{min}(x_3;y_0,y_l)-x_{min}(x_4;y_0,y_l)
\end{equation}
and
\begin{equation}
 \label{eq:der_mathcalXminmonotonic}
\mathcal X_{min}(x_1,y_l)-\mathcal X_{min}(x_2,y_l)\leq
\mathcal X_{min}(x_3,y_l)-\mathcal X_{min}(x_4,y_l).
\end{equation}
\end{proposition}
\begin{proof}
Restrict our solutions $y_k(x;x_k,y_0,y_l)$ on the segments
$[x_{min}(x_k;y_0,y_l),x_k]$, then their inverse functions,
$x_k(y;x_k,y_0,y_l)$, are properly defined on the segment $y\in[y_l,y_0]$.
For brevity, we denote them $x_k(y)$. In particular we have $x_k=x_k(y_0)$.
It is convenient to introduce the following notation:
$$
\Delta_{km}(y):=x_k(y)-x_m(y),\quad
w_k(y):=4(y^3-y_l^3)-2\int^{y}_{y_l}x_k(u)du.
$$
Using them we can rewrite Equation~(\ref{eq:u-int}) as follows,
$$
x_k(y)-x_k(y_0)=-\int^{y_0}_y\frac{du}{\sqrt{w_k(u)}},
$$
and obtain equation for the differences,
\begin{equation}
 \label{eq:delta_km}
\Delta_{km}(y_0)-\Delta_{km}(y)=2\int^{y_0}_y
\frac{\int^y_{y_l}\Delta_{km}(\tilde u)d\tilde u}
{\sqrt{w_k(u)w_m(u)}(\sqrt{w_k(u)}+\sqrt{w_m(u)})}.
\end{equation}
Assume that
\begin{equation}
 \label{ineq:wrong}
\int^{y}_{y_l}\Delta_{12}(\tilde u)d\tilde u>
\int^{y}_{y_l}\Delta_{34}(\tilde u)d\tilde u
\qquad{\rm for}\quad y=y_0.
\end{equation}
Then by continuity this inequality holds for
all $y:\;y^*< y\leq y_0$, for some $y^*\geq y_l$. We further assume that
$y^*$ denotes the infimum of such numbers $y^*$.
Then Equation~(\ref{eq:delta_km}) implies for this values of $y$,
$
\Delta_{12}(y_0)-\Delta_{12}(y)>\Delta_{34}(y_0)-\Delta_{34}(y),
$
since by our conditions $w_1(u)<w_3(u)$ and $w_2(u)<w_4(u)$ for
all $u\in[y_0,y_l)$. Our assumption for points $x_k$ reads
$\Delta_{12}(y_0)\leq\Delta_{34}(y_0)$, therefore we obtain
$\Delta_{12}(y)<\Delta_{34}(y)$ for $y:\;y^*< y\leq y_0$.
If $y^*=y_l$, then we integrate the last inequality for the differences
from $y_l$ to $y_0$ and arrive at a contradiction with
Inequality~(\ref{ineq:wrong}). Thus $y^*>y_l$. In this case we have
$
\int^{y^*}_{y_l}\Delta_{12}(\tilde u)d\tilde u=
\int^{y^*}_{y_l}\Delta_{34}(\tilde u)d\tilde u
$
and
$
\int^{y^0}_{y^*}\Delta_{12}(\tilde u)d\tilde u<
\int^{y^0}_{y^*}\Delta_{34}(\tilde u)d\tilde u
$.
Now summing up the last inequality and equation we again arrive at a
contradiction with Inequality~(\ref{ineq:wrong}).

Thus conditions of our proposition imply that
\begin{equation}
 \label{ineq:right}
\int^{y_0}_{y_l}\Delta_{12}(\tilde u)d\tilde u\leq
\int^{y_0}_{y_l}\Delta_{34}(\tilde u)d\tilde u.
\end{equation}
Now we have two logical possibilities:\\
(1) $\Delta_{12}(y)\leq\Delta_{34}(y)$ for $y:\;y^*\leq y\leq y_0$; and\\
(2) $\Delta_{12}(y)\geq\Delta_{34}(y)$ for $y:\;y^*\leq y\leq y_0$;\\
where again $y^*<y_0$ denotes the infimum of such values of $y$ that the
corresponding condition holds. Consider the second case. If we suppose
that $y^*=y_l$, then we integrate the second condition from $y_l$ to $y_0$
to arrive at
$
\int^{y_0}_{y_l}\Delta_{12}(\tilde u)d\tilde u\geq
\int^{y_0}_{y_l}\Delta_{34}(\tilde u)d\tilde u
$.
The latter inequality does not contradict Inequality~(\ref{ineq:right}) only
if $\Delta_{12}(y)=\Delta_{34}(y)$ for all $y\in[y_0,y_l]$. Otherwise
$y^*>y_l$. In both situations of case (2) we proved that there exists
$y^*<y_0$ and $\Delta_{12}(y^*)=\Delta_{34}(y^*)$. Now note that if we
combine this result with the statement of case (1) we arrive at the
existence of a point $y^*:\;y_l\leq y^*<y_0$ where
\begin{equation}
 \label{ineq:y*}
\Delta_{12}(y^*)\leq\Delta_{34}(y^*).
\end{equation}
Moreover, as the curves $x_k(y)$ for $y\in[y_l,y_0]$ do not intersect.
We observe, that the assumptions of our proposition holds at $y^*$.

Now we consider the set $Y^*=\{y^*:\;y^*\in[y_l,y_0]$,
Inequality~(\ref{ineq:y*}) holds$\}$. Clearly, the set $Y^*$ is non-empty
and closed, moreover by previous construction if a point $y\in Y^*$ and
$y<y_l$, then there is a point $y^*<y$ that belongs to $Y^*$. This means that
$\inf Y^*=y_l$.

Finally, Inequality~(\ref{eq:der_mathcalXminmonotonic}) follows from
Inequality~(\ref{eq:der_xminmonotonic}) by taking a limit $y_0\to+\infty$.
\end{proof}
\begin{corollary}
In the notation of Proposition~{\rm\ref{Prop:der_mathcalXmonotonic}} and
Corollary~{\rm\ref{Cor:3.1}} the following inequality is valid
\begin{equation}
 \label{eq:f_der_xmin}
f(x_2;y_0,y_l)\frac\partial{\partial y_0}x_{min}(x_2;y_0,y_l)\leq
f(x_1;y_0,y_l)\frac\partial{\partial y_0}x_{min}(x_1;y_0,y_l),
\end{equation}
where $x_2<x_1\leq0$.
\end{corollary}
\begin{proof}
Consider four solutions:
$y_1(x)=y(x;x_1,y_0,y_l)$, $y_2(x)=y(x;x_1,y_1,y_l)$,
$y_3(x)=y(x;x_2,y_0,y_l)$, and $y_4(x)=y(x;x_2,\tilde y_1,y_l)$,
where $\tilde y_1>y_1>y_0>y_l$. Consider the inverse functions,
$x_k(y)$, as in the proof of Proposition~\ref{Prop:der_mathcalXmonotonic}.
We have $x_1=x_1(y_0)$, $x_2=x_3(y_0)$. Now we impose an additional
condition on the parameter $\tilde y_1$,
$x_2(y_0)-x_1(y_0)=x_3(y_0)-x_4(y_0)$ and write
Equation~(\ref{eq:der_xminmonotonic}) for the case under consideration
$$
x_{min}(x_1;y_0,y_l)-x_{min}(x_1;y_1,y_l)\leq
x_{min}(x_2;y_0,y_l)-x_{min}(x_2;\tilde y_1,y_l).
$$
Divide now both parts of the last inequality by $\tilde y_1-y_0$, and
make a passage to a limit $\tilde y_1\to y_0$, taking into account that
$\lim(\tilde y_1-y_0)/(x_3(y_0)-x_4(y_0))=f(x_2;y_0,y_l)$ and
$\lim(y_1-y_0)/(x_1(y_0)-x_2(y_0))=f(x_1;y_0,y_l)$.
\end{proof}
\begin{remark}
Actually, we expect that there is a strict monotonicity in
Equations~(\ref{eq:der_xminmonotonic}), (\ref{eq:der_mathcalXminmonotonic}),
and (\ref{eq:f_der_xmin}).
\end{remark}
\section{The Functions ${\mathcal X}(x_0, y_l)$ and $X(x_0)$}
 \label{sect:4}
\begin{lemma}
 \label{Lem:4.1}
Let $y_1(x)$ and $y_2(x)$ be solutions that intersect at $x_0$, i.e.,
$y_1(x_0)=y_2(x_0)$, such that $y'_1(x_0)<y'_2(x_0)\leq0$. Denote by $a_k$
the left end of the interval of existence of $y_k(x)$
{\rm (}$k=1, 2${\rm )}. Then $a_2<a_1$ and
\[
y_1(x)>y_2(x), y_1'(x)<y_2'(x)<0,\quad {\rm for}\ x\in(a_1, x_0).
\]
\end{lemma}
\begin{proof}
The conditions on the solutions can be reformulated as $u_1(t_0)=u_2(t_0)$,
$u_1'(t_0)>u_2'(t_0)\leq0$ in our usual notation $t=-x$, $u(t)=y(x)$.
Denote $\hat a_k=-a_k$. In some open right neighbourhood of $t_0$ we have
\begin{equation}
 \label{eq:K}
u_1(t)>u_2(t)\quad {\rm and}\quad u_1'(t)>u_2'(t).
\end{equation}
Denote
\[
\hat a=\sup\left\{t\Bigm| {\rm Condition}\ (\ref{eq:K})\ {\rm holds}\right\}.
\]
Suppose $\hat a<\min\{a_1, a_2\}$.
We have $u_1'(\hat a)=u_2'(\hat a)$ and $u_1(\hat a)>u_2(\hat a)$.
Moreover,
\begin{equation}
 \label{eq:ham}
{u_k'}^2(\hat a)-4{u_k}^3(\hat a)={u_k'}^2(t_0)-4{u_k}^3(t_0)
+2\int_{t_0}^{\hat a}\,tu_k'(t)\,dt
\end{equation}
for $k=1,2$. Subtracting Equation (\ref{eq:ham}) for $k=2$ from that for
$k=1$, we get
\[
4\left(u_2^3(\hat a)-u_1^3(\hat a)\right)={u_1'}^2(t_0)-{u_2'}^2(t_0)
+2\int_{t_0}^{\hat a}\,t\bigl(u_1'(t)-u_2'(t)\bigr)\,dt
\]
This is a contradiction as the right-hand side of the equation
is positive whilst the left-hand side is negative.
Therefore, $\hat a=\min\{a_1, a_2\}$. Since
$u_1(t)>u_2(t)$, we have $\hat a_1\le \hat a_2$. Suppose that
$\hat a_1=\hat a_2=\hat a$. Since
$u_1'(t)>u_2'(t)$, for $t\le \hat a$, we have $u_1(t)-u_2(t)$
grows as $t\to \hat a$. However, according to Equation (\ref{Eqn:1}),
$u_1(t)-u_2(t)\to0$ as $t\to\hat a$.
\end{proof}
\begin{lemma}
 \label{Lem:4.2}
Let $y_1(x)$ and $y_2(x)$ be solutions with the same minimum value $y_l$.
More precisely,
$y_1(x_1)=y_2(x_2)=y_l$, $y_1'(x_1)=y_2'(x_2)=0$, for $x_1<x_2$.
Then $a_2> a_1$ and
\[
y_1(x)<y_2(x),\quad y_2'(x)<y_1'(x)<0,\quad {\rm for}\quad x\in(a_2, x_1).
\]
\end{lemma}
\begin{proof}
In terms of the variables $t=-x$, $u(t)=y(x)$, the conditions read:
$t_1>t_2>0$, $u_1(t_1)=u_2(t_2)=y_l=:u_l$, $u_1'(t_1)=u_2'(t_2)=0$.
In some open right neighbourhood of $t_1$ we have
\begin{equation}
 \label{eq:L}
u_2(t)>u_1(t)\quad {\rm and}\quad u_2'(t)>u_1'(t).
\end{equation}
Denote
\[
\hat a=\sup\left\{t\Bigm| {\rm Condition}\ (\ref{eq:L})\ {\rm holds}\right\}.
\]

Suppose $\hat a<\min\{a_1, a_2\}$.
Then $u_1'(\hat a)=u_2'(\hat a)$ and $u_2(\hat a)>u_1(\hat a)$.
Integrating Equation (\ref{p1t}) from $t_k$, we have
\[
{u_k'}^2(\hat a)-4{u_k}^3(\hat a)=-4{u_l}^3
+2\int_{t_k}^{\hat a}\,tu_k'(t)\,dt
\]
for $k=1,2$. Subtracting and rearranging the integral, we get
\begin{equation}
 \label{eq:int_tk}
{u_1'}^2(\hat a)-{u_2'}^2(\hat a)=4\left(u_1^3(\hat a)-u_2^3(\hat a)\right)
+2\int_{t_1}^{\hat a}\,t\bigl(u_1'(t)-u_2'(t)\bigr)\,dt
-2\int_{t_2}^{t_1}\,t\,u_2'(t)\,dt
\end{equation}
The right-hand side of Equation (\ref{eq:int_tk}) is strictly negative.
Hence we have ${u_1'}^2(\hat a)<{u_2'}^2(\hat a)$
which is a contradiction. The remainder of the proof is the same as
the proof of the previous lemma with interchange of the indices $1$ and
$2$.
\end{proof}
\begin{lemma}
 \label{Lem:4.3}
Let $y_1(x)$ and $y_2(x)$ be solutions satisfying the following conditions:
$y_1(x_1)=y_2(x_2)=y_l$, $0\geq y_1'(x_1)\geq y_2'(x_2)$, for $x_1<x_2$.
Then $a_2>a_1$ and
\[
y_1(x)<y_2(x),\quad y_2'(x)<y_1'(x)<0,\quad {\rm for}\quad x\in(a_2, x_1).
\]
\end{lemma}
\begin{proof}
This is a generalization of the previous Lemma. The proof is literally
the same. The only difference is an additional term,
${u'_1}^2(t_1)-{u'_2}^2(t_2)={y'_1}^2(x_1)-{y'_2}^2(x_2)\leq0$,
which does not spoil the proof.
\end{proof}
\begin{definition}
 \label{Def:4.1}
Let $x_0\le 0$ and $x_{min}$ be the minimum of the solution
$y(x)$ corresponding to initial
data $(y_0, y_0')$ at $x_0$. Define
\[
{\mathcal X}(x_0, y_l)=
\underset{y_0,y_0'}{\inf}\left\{a(x_0, y_0, y_0')\Bigm|
y(x_{min})=y_l\right\}.
\]
\end{definition}
\begin{proposition}
 \label{Prop:4.1}
For any $x_0\le 0$, $y_l\in{\mathbb R}$, ${\mathcal X}(x_0, y_l)$ is
finite. The unique solution with a pole at $x_0$ and
$x_{min}={\mathcal X}_{min}(x_0,y_l)$ has a pole at
${\mathcal X}(x_0, y_l)$. Moreover, ${\mathcal X}(x_0, y_l)$
is a smooth function of $(x_0, y_l)$.
\end{proposition}
\begin{proof}
Since $X(x_0)\le{\mathcal X}(x_0,y_l)$, it follows from
Theorem~\ref{Thm:1}
that ${\mathcal X}(x_0, y_l)$ is bounded for any $x_0\le 0$,
$y_l\in{\mathbb R}$. The statement that the solution with minimum at
${\mathcal X}_{min}(x_0,y_l)$ has a pole at ${\mathcal X}(x_0, y_l)$
follows from Theorem \ref{Thm:3.1} and Lemma \ref{Lem:4.2}. The smoothness
of ${\mathcal X}(x_0, y_l)$ is a consequence of
Proposition \ref{Prop:analyticity} and Corollary \ref{Cor:3.3}.
\end{proof}
\begin{corollary}
 \label{Cor:4.1}
For $y_l\in{\mathbb R}$, if $x_1<x_2<0$ then
\begin{eqnarray*}
{\mathcal X}(x_1,y_l)&<&{\mathcal X}(x_2,y_l).
\end{eqnarray*}
\end{corollary}
\begin{proof}
The proof follows from the monotonicity result for
${\mathcal X}_{min}(x_0, y_l)$ (Corollary~\ref{Cor:3.2})
and Lemma~\ref{Lem:4.2}.
\end{proof}
\begin{proposition}
 \label{Prop:4.2}
The function ${\mathcal X}(x_0, y_l)$ satisfies the following inequalities,
\begin{eqnarray}
 \label{ineq:4.5}
x_1<x_2\leq0:&&
0<\mathcal X(x_2, y_l)-\mathcal X(x_1, y_l)
<\mathcal X_{min}(x_2, y_l)-\mathcal X_{min}(x_1, y_l),\\
\label{ineq:4.6}
x_0\leq0:&&
\mathcal X_{min}(x_0, y_l)<\frac12(x_0+\mathcal X(x_0, y_l)).
\end{eqnarray}
\end{proposition}
\begin{proof}
The left Inequality~(\ref{ineq:4.5}) is proved in Corollary \ref{Cor:4.1}.
The right Inequality~(\ref{ineq:4.5}) follows from the fact that the graph
of the solution starting at $\mathcal X_{min}(x_1, y_l)$ goes to its
pole more steeply than the one starting at $\mathcal X_{min}(x_2, y_l)$.
Analogous arguments also justify Inequality~(\ref{ineq:4.6}); if we denote
$y_1(x)=y(x)$ for $x\in(x_0,\mathcal X_{min}(x_0, y_l)]$ and
$y_2(x)=y(x)$ for $x\in[\mathcal X_{min}(x_0, y_l),\mathcal X(x_0, y_l))$,
then $y_2(x)$ is steeper than $y_1(x)$.
The formal proof is analogous to that of Proposition~\ref{Prop:3.5} with
corresponding changes of notation.
\end{proof}
\begin{proposition}
 \label{Prop:4.3}
For $|y_l|\geq\epsilon>0$
$$
0<x_0-{\mathcal X}(x_0,y_l)<\frac{I_\nu}{\sqrt{|y_l|}},
$$
where $\nu=\mathrm{sign}\,\{y_l\}1$ and $I_\nu$ is defined by
Equation~{\rm(\ref{eq:intdelta})}.
\end{proposition}
\begin{proof}
Note that
$
0<{\mathcal X}_{min}(x_0,y_l)-{\mathcal X}(x_0,y_l)<
x_0-{\mathcal X}_{min}(x_0,y_l),
$
as the graph of the solution with poles at $x_0$ and $\mathcal{X}(x_0,y_l)$
is steeper to the left of its minimum $x_{min}=\mathcal{X}_{min}(x_0,y_l)$
(see Proposition~\ref{Prop:4.1}) than that to the right. The formal proof
is very similar to the one of Proposition~\ref{Prop:3.5}.
Now, the result follows from the identity,
$
x_0-{\mathcal X}(x_0,y_l)<
x_0-{\mathcal X}_{min}(x_0,y_l)+
{\mathcal X}_{min}(x_0,y_l)-{\mathcal X}(x_0,y_l),
$
and Proposition~\ref{Prop:3.6}.
\end{proof}
\begin{proposition}
 \label{Prop:lim_der_x}
For any $y_l\in\mathbb R$ and $x_0\leq0$,
\begin{equation}
 \label{ineq:der01}
0\leq\frac\partial{\partial x_0}
\mathcal X_{min}(x_0,y_l)<1,\qquad
0\leq\frac\partial{\partial x_0}
\mathcal X(x_0,y_l)<1.
\end{equation}
Moreover,
\begin{equation}
 \label{eq:der_limit}
\underset{\kappa_-+\kappa_+\to+\infty}{\lim}\;\;\frac\partial{\partial x_0}
\mathcal X_{min}(x_0,y_l)=1,\qquad
\underset{\kappa_-+\kappa_+\to+\infty}{\lim}\;\;\frac\partial{\partial x_0}
\mathcal X(x_0,y_l)=1,
\end{equation}
where $\kappa_-=-x_0$, $\kappa_+=|y_l|$.
\end{proposition}
\begin{proof}
The fact that the values of the derivatives belong to $[0,1]$
follows from Propositions~\ref{Prop:3.5} and \ref{Prop:4.2}.
If we consider $x_{min}$ in Equation~(\ref{eq:int_x0-xmin}) as a
function of $x_0$ and $y_l$, then it is nothing but
$\mathcal X_{min}(x_0,y_l)$. Differentiating now both parts of the
equation with respect to $x_0$, we get on the left
$1-\frac\partial{\partial x_0}\mathcal X_{min}(x_0,y_l)$ and on the
right an integral similar to the one in the right part of
Equation~(\ref{eq:derivative}) but with $x_{min}$ changed to $x_0$. Taking
now into account that $0\leq\frac\partial{\partial x_0}x(y;x_0,y_l)\leq1$,
one finds that the integral vanishes as $x_0\to-\infty$ and/or
$|y_l|\to+\infty$. The partial derivative cannot be identically equal to
$0$ for all $y$, as it implies that the partial derivative of $y(x;x_0,y_l)$
(with respect to $x_0$) identically vanishes for all $x$. This contradicts
Laurent expansion~(\ref{laurent}). Thus, the right
Inequalities~(\ref{ineq:der01}) are strict.

The proof for $\frac\partial{\partial x_0}\mathcal X(x_0,y_l)$ is
similar; instead of one integral in the right-hand side
of the equation for $1-\frac\partial{\partial x_0}\mathcal X(x_0,y_l)$
we get two similar integrals, and then, the same argument as above shows
that both vanish as $x_0\to-\infty$.
\end{proof}
\begin{remark}
In Proposition~\ref{Prop:der_strict} of \S~\ref{sect:5} we prove that
both left Inequalities~(\ref{ineq:der01}) are also strict for all
$x_0$ and $y_l\in\mathbb R$.
\end{remark}
\begin{proposition}
 \label{Prop:existence}
For any $x_0\le0$, there exists a solution with interval of existence
$(X(x_0),x_0)$. If some solution has an interval of existence
$I\supset(X(x_0),x_0)$, then $I=(X(x_0),x_0)$.
\end{proposition}
\begin{proof}
By definition, we have
$
\underset{y_l}{\inf}\,{\mathcal X}(x_0,y_l)\ge X(x_0).
$
On the other hand, for any $\epsilon>0$, there exists a solution
$y(x)$ with initial data $y(x_0)=y_0$ and $y'(x_0)=y_1$ such that its left
pole $a(x_0,y_0,y_1)<X(x_0)+\epsilon$. Let $y_l$ be its minimum.
Then ${\mathcal X}(x_0,y_l)< a(x_0,y_0,y_1)< X(x_0)+\epsilon$. Therefore
\begin{equation}
 \label{eq:X=infcalX}
X(x_0)=\underset{y_l}{\inf}\,{\mathcal X}(x_0,y_l).
\end{equation}
Since ${\mathcal X}(x_0,y_l)$ is differentiable in $y_l$ and by virtue of
Proposition~\ref{Prop:4.3} ${\mathcal X}(x_0,y_l)\to x_0$ as
$y_l\to\pm\infty$, it achieves its supremum at some finite value $\hat y_l$.
Actually $\hat y_l\leq0$. Thus, $X(x_0)={\mathcal X}(x_0,\hat y_l)$, and
Proposition~\ref{Prop:4.1} implies that there is a solution with the
interval of existence $(X(x_0),x_0)$.

Denote by $y_1(x)$ the solution with interval of existence $I$. If $y_1(x)$
has a pole at $x_0$, then it has another pole $p\geq X(x_0)$, as follows from
Equation~(\ref{eq:X=infcalX}) and thus its interval of existence
$I\subset(X(x_0),x_0)$.

Suppose now that $y_1(x)$ is finite at $x_0$, then by
Proposition~\ref{Prop:3.2} and Definition~\ref{Def:3.1},
${\mathcal X}_{min}(x_0,y_l)<x_{min}$, where $y_l$ and $x_{min}$ are the
minimum value and minimum of $y_1(x)$. Denoting by $a$ the left bound of $I$,
and applying Lemma~\ref{Lem:4.2} and Definition~\ref{Def:4.1}, we get that
${\mathcal X}(x_0,y_l)<a$. This inequality contradicts
Equation~(\ref{eq:X=infcalX}) as on the other hand $a=X(x_0)$: by definition
of $X(x_0)$ we have $a\geq X(x_0)$, whilst by definition of $I$,
$a\leq X(x_0)$.
\end{proof}
\begin{remark}
If $I\subset(-\infty,0]$, then the last paragraph of the proof of
Proposition~\ref{Prop:existence} can be simplified: Denote by $b$ the right
bound of $I$, then Corollary~\ref{Cor:4.1} implies
${\mathcal X}(x_0,y_l)<{\mathcal X}(b,y_l)\equiv a\leq X(x_0)$
and we arrive at a contradiction with Equation~(\ref{eq:X=infcalX}).
\end{remark}
\begin{corollary}
$
x_0\leq0:\quad
X(x_0)<X_{min}(x_0)<\frac12(x_0+X(x_0)).
$
\end{corollary}
\begin{proof}
Immediate consequence of Propositions~\ref{Prop:existence} and
\ref{Prop:4.2}.
\end{proof}
\begin{proposition}
 \label{Prop:XLipschitz}
The function $X(x_0)$ is monotonically increasing with increase of
$x_0\leq0$ and satisfies the Lipschitz condition,
$$
0<X(x_2)-X(x_1)<x_2-x_1.
$$
\end{proposition}
\begin{proof}
The proof follows from Proposition~\ref{Prop:4.2} and
Proposition~\ref{Prop:3.5} by similar arguments as in the proof of
Proposition~\ref{Prop:3.7}.
\end{proof}
\begin{remark}
The solution introduced in Proposition~\ref{Prop:existence} is unique
and smooth. The proof is given in \S~\ref{sect:main}, see
Theorem~\ref{Thm:uniqueness} and Corollary~\ref{Cor:smoothness}, respectively.
Most probably instead of the Lipschitz condition proved in
Proposition~\ref{Prop:XLipschitz}, a stronger inequality holds, analogous
to the one proved for the level-functions in Proposition~\ref{Prop:4.2}
(see Conjecture~\ref{Con:5} in \S~\ref{sect:numerics}).
\end{remark}
\section{The functions $\varXi_{min}(x_0,y_l)$, $\Xi_{min}(x_0)$,
$\varXi(x_0,y_l)$, $\Xi(x_0)$, and $X_-(x_0)$}
 \label{sect:5}
\begin{remark}
 \label{Rem:5.1}
In those statements below where no condition is given on $x_0$,
we assume that $x_0\leq{\mathcal X}(0,y_l)$.
This condition guarantees that each solution under consideration has an
interval of existence which belongs to the negative semi-axis.
\end{remark}
\begin{proposition}
 \label{Prop:5.1}
Given $y_0$ and $y_l$ with $y_0>y_l$, there exists only one solution $y(x)$
such that $y(x_0)=y_0$ and $y_{min}=y_l$ where $y_{min}$
is the minimum of $y(x)$ achieved to the right of $x_0$.
\end{proposition}
\begin{proof}
The existence of the solution follows by arguments similar to those
in the first two paragraphs of the proof of Proposition \ref{Prop:3.1},
with $t_0$ substituted for $t_{min}$ in Equation (\ref{eq:Prop:3.1}).
The uniqueness follows from Lemma \ref{Lem:4.2}.
\end{proof}
\begin{proposition}
 \label{Prop:5.2}
Given $y_0$, $\tilde y_0$, and $y_l$ with $\tilde y_0>y_0> y_l$,
consider solutions $y(x)$, $\tilde y(x)$ such that $y(x_0)=y_0$,
$\tilde y(x_0)=\tilde y_0$ and with the same minimum value
$y_l=y_{min}=\tilde y_{min}$ achieved to the right of $x_0$.
Then $x_0<x_{min}<\tilde x_{min}$ and $\tilde y'(x_0)<y'(x_0)<0$.
\end{proposition}
\begin{proof}
Follows from Lemma \ref{Lem:4.2}.
\end{proof}
\begin{corollary}
 \label{Cor:5.1}
Consider the solution with initial value $y(x_0)=y_0$ and minimum value
$y_l$ achieved to the right of $x_0$. The initial slope
$y_0'\equiv f^+(x_0;y_0,y_l)$ {\rm(}see Remark {\rm\ref{Rem:symmetry}}{\rm)}
at $x_0$ is a function of $y_0$ and, for fixed $y_l$, has asymptotic
behaviour
\[
\frac{{y_0'}^2}{y_0}\underset{y_0\to+\infty}{=}4{y_0}^2+O(1).
\]
\end{corollary}
\begin{proof}
Literally the same as the proof of Corollary \ref{Cor:3.1} with the
corresponding change of $x_0\leftrightarrow x_{min}$.
\end{proof}
\begin{definition}
 \label{Def:5.1}
Given $x_0$, $y_0>y_l$, consider the solution
$y(x)=y(x;x_0,y_0,y_l)$ with initial value $y(x_0)=y_0$, minimum value
$y_l$ at $x_{min}=x_{min}(x_0,y_0,y_l)>x_0$, and the right bound of its
interval of existence $b=b(x_0,y_0,y_l)$. Define
$$
\varXi_{min}(x_0,y_l)=\underset{y_0}{\sup}\,x_{min}(x_0,y_0,y_l),\qquad
\varXi(x_0,y_l)=\underset{y_0}{\sup}\,b(x_0,y_0,y_l).
$$
\end{definition}
\begin{remark}
The functions $\varXi_{min}(x_0,y_l)$ and $\varXi(x_0,y_l)$ are finitely
defined for $y_l\in\mathbb{R}$, $-\infty<x_0\leq\mathcal X(0,y_l)$,
 and satisfy evident inequalities:
$$
x_0<\varXi_{min}(x_0,y_l)<\varXi(x_0,y_l)<0.
$$
\end{remark}
\begin{proposition}
 \label{Prop:5.3}
The solution $y_\mu(x)$ with initial data
\[y_\mu\bigl(\varXi_{min}(x_0,y_l)\bigr)=y_l,\quad
y_\mu'\bigl(\varXi_{min}(x_0,y_l)\bigr)=0\]
has an interval of existence $(x_0,\varXi(x_0,y_l))$.
\end{proposition}
\begin{proof}
{}From Proposition~\ref{Prop:5.2} and Definition~\ref{Cor:5.1} it follows
that the solution $y_\mu(x)$ has a pole $x_p$,
$x_0\leq x_p<\varXi_{min}(x_0,y_l)$. The equality $x_p=x_0$ can be
confirmed by making an analogous construction to that in the last
paragraph of the proof of Proposition~\ref{Prop:3.3}.

According to Proposition~\ref{Prop:5.2}, any solution with finite initial
value at $x_0$ has a minimum strictly less than $\varXi_{min}(x_0,y_l)$.
The graph of such a solution crosses the graph of $y_\mu(x)$ at a point to
the left of its minimum $(\varXi_{min}(x_0,y_l),y_l)$ and does not cross
to the right of it by virtue of Proposition~\ref{Prop:3.1}. Thus, the right
bound of the interval of existence, $b$, of any solution finite at $x_0$
is less than $b_\mu$, the right bound of the interval of existence of
$y_\mu(x)$. On the other hand, by similar arguments as above we establish
that any sequence of solutions $y_n(x)$ with finite initial data
$y_n(x_0)\nearrow+\infty$, has a sequence of minima
$x_{min,n}\nearrow\varXi_{min}(x_0,y_l)$ and the sequence of the
right bounds of their intervals of existence, $b_n\nearrow b_\mu$.
\end{proof}
\begin{corollary}
 \label{Cor:5.2}
$x_0={\mathcal X}(\varXi(x_0,y_l),y_l)$. In other words,
$\varXi(x_0,y_l)={\mathcal X}^{-1}(x_0,y_l)$, where
${\mathcal X}^{-1}(x_0,y_l)$ denotes the section of the fibre
${\mathcal X}^{-1}(x_0)$ by the straight line $y=y_l$.
\end{corollary}
\begin{proof}
By virtue of Theorem~\ref{Thm:3.1} the solution with a pole at
$\varXi(x_0,y_l)$ and minimum value $y_l$ is unique. Therefore, it
coincides with the solution defined in Proposition~\ref{Prop:5.3}.
From Proposition~\ref{Prop:4.1} follows that this solution has a
pole at ${\mathcal X}(\varXi(x_0,y_l),y_l)$ while
Proposition~\ref{Prop:5.3} implies that this pole is $x_0$.
\end{proof}
\begin{proposition}
 \label{Prop:5.4}
Given $x_0\leq0$ and $y_l\in{\mathbb R}$ and  $x_0\leq\mathcal X(0,y_l)$,
there is only one solution with a pole at $x_0$, an interval of existence to
the right of $x_0$, and minimum value is $y_l$.
\end{proposition}
\begin{proof}
Existence of the solution follows from Proposition \ref{Prop:5.3},
uniqueness from Lemma \ref{Lem:4.2}.
\end{proof}
\begin{theorem}
 \label{Thm:5.1}
Any solution with an interval of existence that is contained in the
non-positive semi-axis can be uniquely characterised by the position of
the left pole $x_0$ and its minimum value $y_l$. Moreover, it has the
following extremal property. Its minimum $\varXi_{min}(x_0,y_l)>x_{min}$
where $x_{min}$ is the minimum of any solution regular at $x_0$ with
minimum value $y_l$.
\end{theorem}
\begin{proof}
The proof follows from Propositions \ref{Prop:5.1}--\ref{Prop:5.4}.
\end{proof}
The above results can be formulated as monotonicity properties of
$\varXi_{min}(x_0,y_l)$.
\begin{corollary}
 \label{Cor:5.3}
If $x_1<x_2<{\mathcal X}(0,y_l)$, where $y_l\in\mathbb R$, then
\begin{eqnarray*}
\varXi_{min}(x_1,y_l)&<&\varXi_{min}(x_2,y_l)
\end{eqnarray*}
\end{corollary}
\begin{proof}
If $\varXi_{min}(x_1,y_l)=\varXi_{min}(x_2,y_l)$
then the two corresponding solutions (see Proposition \ref{Prop:3.3})
would coincide.
Suppose we have $\varXi_{min}(x_1,y_l)>\varXi_{min}(x_2,y_l)$.
Then there is a solution regular at $x_2$, with
$x_{min}>\varXi_{min}(x_2,y_l)$. This contradicts the definition
of $\varXi_{min}(\cdot,y_l)$.
\end{proof}
\begin{proposition}
 \label{Prop:5.5}
For any $x_0\leq0$, $x_1<x_2\leq{\mathcal X}(0,y_l)$, where
$y_l\in\mathbb R$:
\begin{eqnarray}
 \label{eq:5.1}
{\mathcal X}_{min}(x_0,y_l)\!\!&=&
\!\!\varXi_{min}({\mathcal X}(x_0,y_l),y_l)<x_0,\\
 \label{ineq:5.2}
0<x_2-x_1\!\!&<&\!\!\varXi_{min}(x_2,y_l)-\varXi_{min}(x_1,y_l)<
\varXi(x_2,y_l)-\varXi(x_1,y_l),\\
\label{ineq:5.3}
\varXi_{min}(x_0,y_l)&<&\frac12\left(x_0+\varXi(x_0,y_l)\right).
\end{eqnarray}
\end{proposition}
\begin{proof}
Equation (\ref{eq:5.1}) holds because the solutions corresponding to its
left- and right-hand sides coincide.
To prove Inequality (\ref{ineq:5.2}) it is enough to notice that
the graph of the solution starting at $\varXi_{min}(x_1,y_l)$ goes to its
pole more steeply than the one starting at $\varXi_{min}(x_2, y_l)$.
A similar argument works to prove Inequality~(\ref{ineq:5.3}); a
graph of the solution to the left of $\varXi_{min}(x_0,y_l)$ is steeper
than that to the right.
The formal proof of Inequalities~(\ref{ineq:5.2}) and (\ref{ineq:5.3}) is
analogous to that of Proposition \ref{Prop:3.5} with corresponding changes
of notation.
\end{proof}
\begin{remark}
For $x_0\leq\mathcal X(0,y_l)$ we also have
$\varXi_{min}(x_0,y_l)=\mathcal X_{min}(\varXi(x_0,y_l),y_l)<0$.
\end{remark}
\begin{remark}
A solution with a pole at $x_0$ and interval of existence to the right of
it, is uniquely characterised on the one hand by the parameter $c$ in the
Laurent expansion (\ref{laurent}) and on the other hand, by its minimum
$y_l$ achieved at $x_{min}=x_{min}(x_0,c)$.
Therefore, there exists a bijection $y_l(c):{\mathbb R}\to{\mathbb R}$.
\end{remark}
\begin{corollary}
 \label{Cor:5.4}
For any $c$, $y_l(c)$ is a smooth monotonically growing function of $c$.
The inverse function $c(y_l)$ has the same properties.
Moreover, $y_l(c)\to\pm\infty$ as $c\to\pm\infty$.
For any $x_0\leq{\mathcal X}(0,y_l)$, $x_{min}(c)=x_{min}(x_0,c)$ is a
smooth function of $c$.
\end{corollary}
\begin{proof}
The proof is analogous to the one for Corollary \ref{Cor:3.3}.
\end{proof}
\begin{corollary}
 \label{Cor:5.5}
The functions $\varXi_{min}(x_0,y_l)$ and $\varXi(x_0,y_l)$
are smooth functions of both variables.
\end{corollary}
\begin{proof}
Follows from Propositions~\ref{Prop:5.3}, \ref{Prop:analyticity} and
Corollary \ref{Cor:5.4}.
\end{proof}
\begin{definition} For $x_0\leq X(0)$ define
$$
\Xi_{min}(x_0)=\underset{y_0,y_l}{\sup}\,x_{min}(x_0;y_0,y_l),\qquad
\Xi(x_0)=\underset{y_0,y_l}{\sup}\,b(x_0;y_0,y_l)
$$
where $x_{min}(x_0;y_0,y_l)$ and $b(x_0;y_0,y_l)$ are the same as in
Definition \ref{Def:5.1}.
\end{definition}
\begin{proposition}
 \label{Prop:5.6}
For any $x_0\leq X(0)$: {\rm(1)} There exists a solution of
Equation {\rm(\ref{p1})} that has a minimum at $\Xi_{min}(x_0)$.
Any such solution has a pole at $x_0$; {\rm(2)} There exists a solution
with the interval of existence $(x_0,\Xi(x_0))$.
\end{proposition}
\begin{proof}
Taking into account that $\varXi_{min}(x_0,y_l)-x_0>0$ and
$\varXi(x_0,y_l)-x_0>0$ vanish as $|y_l|\to\infty$
(Proposition~\ref{Prop:4.3} and Corollary~\ref{Cor:5.3}) and
continuity of $\varXi_{min}(x_0,y_l)$ and $\varXi(x_0,y_l)$
with respect to $y_l$, one proves the following representations for
$\Xi_{min}(x_0)$ and $\Xi(x_0)$:
\begin{equation}
 \label{eq:5.4}
\Xi_{min}(x_0)=\underset{y_l}{\max}\,\varXi_{min}(x_0,y_l),\quad
\Xi(x_0)=\underset{y_l}{\max}\,\varXi(x_0,y_l)
\end{equation}
The maxima are achieved at some finite values of $y_l$.
\end{proof}
\begin{remark}
 \label{Rem:conjecture2}
We conjecture that the solution from Proposition \ref{Prop:5.6}
corresponding to $\Xi_{min}(x_0)$ is unique; the argument is analogous
to the one in Remark~\ref{Rem:conjecture}. The solution corresponding
to $\Xi(x_0)$ is unique. This follows from Proposition~\ref{Prop:5.7}
and Theorem~\ref{Thm:uniqueness}.
\end{remark}
\begin{proposition}
 \label{Prop:5.7}
$\Xi_{min}(X_-(x_0))=x_0$, $X_-(\Xi_{min}(x_0))=x_0$,
$\Xi(X(x_0))=x_0$, $X(\Xi(x_0))=x_0$,
i.e., the functions $\Xi_{min}(x_0)$ and $\Xi(x_0)$ are inverse to
$X_-(x_0)$ and $X(x_0)$ respectively.
\end{proposition}
\begin{proof}
Lemma~\ref{Lem:4.1} implies that in the definition
of $X_-(x_0)$ (see Definition~\ref{Def:2.2}) the infimum is achieved
for $y_1=0$. Moreover, for any $x_0\leq0$ there exists a solution $y(x)$
with $y'(x_0)=0$ and a pole at $X_-(x_0)$. Actually, denote by $y_n$ a
sequence of initial data at $x_0$ such that the corresponding solutions
$y_n(x)$ ($y_n(x_0)=y_n$ and $y'_n(x_0)=0$) have the left bounds
$a_n\to X_-(x_0)$. The sequence $y_n$ is bounded as
$0<x_0-a_n<x_0-X(x_0)<I_\nu/|y_n|$ (see Proposition~\ref{Prop:4.3}). So
there is a convergent subsequence $y_{n_k}\to\hat y$. The solution $y(x)$
with initial data $y(x)=\hat y$ and $y'(x_0)=0$ has a pole at $X_-(x_0)$.

It follows from the previous paragraph that $\Xi_{min}(X_-(x_0))\geq x_0$.
Suppose that $\Xi_{min}(X_-(x_0))>x_0$. Then, as it follows
from Proposition~\ref{Prop:5.6}, there is a solution $y_1(x)$ with the pole
at $X_-(x_0)$ whose graph crosses the vertical line $x=x_0$ with a slope
$y_1'(x_0)<0$. According to Lemma~\ref{Lem:4.1} a solution $y_2(x)$ with the
initial data $y_2(x_0)=y_1(x_0)$ and  $y'_2(x_0)=0$ has a left bound of the
interval of existence less than $X_-(x_0)$; this contradicts the
definition of $X_-(x_0)$.
Now, we prove that $X_-(\Xi_{min}(x_0))=x_0$. Proposition~\ref{Prop:5.6}
and the definition of the function $X_-(\cdot)$ imply that
$X_-(\Xi_{min}(x_0))\leq x_0$. If we suppose that $X_-(\Xi_{min}(x_0))<x_0$,
then this means that there exists a solution finite at $x_0$ with the minimum
at $\Xi_{min}(x_0)$. This contradicts the definition of $\Xi_{min}(x_0)$ by
virtue of Proposition~\ref{Prop:5.2}.
The proof that $X(x_0)=\Xi^{-1}(x_0)$ is very similar.
\end{proof}
\begin{corollary}
The functions $\Xi_{min}(x_0)$ and $\Xi(x_0)$ are continuous and satisfy
the following inequalities:
\begin{eqnarray}
 \label{eq:5.6}
\Xi_{min}\left(\mathcal X(x_0,y_l)\right)\geq\mathcal X_{min}(x_0,y_l),&&
\Xi_{min}\left(X(x_0)\right)\geq X_{min}(x_0),\\
 \label{eq:5.7}
0<x_2-x_1<\Xi_{min}(x_2)-\Xi_{min}(x_1),&&0<x_2-x_1<\Xi(x_2)-\Xi(x_1)\\
\Xi_{min}(x_0)<\frac12(x_0+\Xi(x_0)).
\end{eqnarray}
\end{corollary}
\begin{proof}
The continuity follows from the fact that both functions, $\Xi(x_0)$ and
$\Xi_{min}(x_0)$ are the inverses of the continuous monotonic functions
$X(x_0)$ and $X_-(x_0)$, respectively (concerning $X_-(x_0)$ see
Remark~\ref{Rem:XX_-}). The inequalities follow from
Equation~(\ref{eq:5.1}), Inequalities~(\ref{ineq:5.2}) and (\ref{ineq:5.3}),
and Equations~(\ref{eq:5.4}).
\end{proof}
\begin{remark}
 \label{Rem:XX_-}
The function $X_-(x_0)$ has very similar properties to those of
the function $X_{min}(x_0)$. To prove this one can develop the technique
based on the introduction of the level function, $X_-(x_0,y_l)$, using
Lemmas~\ref{Lem:4.1}-\ref{Lem:4.3} as the main instrument. Analogously,
one can define a function
$\Xi_+(x_0):=\underset{y_0\in\mathbb R,y_1\geq0}{\sup}\{b(x_0,y_0,y_1)\}$,
where $b(x_0,y_0,y_1)>x_0$ is the right bound of the interval of existence of
solution $y_+(x;x_0,y_0,y_1)$. The function $\Xi_+(x_0)$ is the inverse to
$X_{min}(x_0)$ and has properties similar to those of $\Xi_{min}(x_0)$.
The proof again can be based on  the corresponding level function,
$\varXi_+(x_0,y_l)$ with Lemma~\ref{Lem:3.1} as the main instrument.
None of these constructions require any new ideas and we omit them.

Introduction of these functions has not only some notational
convenience, as we see below in Proposition~\ref{Prop:der_strict}, it
allows us to prove an important property of the derivatives of
$\mathcal X$- and $\varXi$-functions, that otherwise would be difficult to
establish directly by working with the corresponding integral equation.
\end{remark}
\begin{proposition}
 \label{Prop:der_strict}
For any $y_l\in\mathbb R$ and $x_0\leq0$,
\begin{equation}
 \label{ineq:Xlevel_der_strict}
0<\frac\partial{\partial x_0}\mathcal X_{min}(x_0,y_0)<1,\qquad
0<\frac\partial{\partial x_0}\mathcal X(x_0,y_0)<1.
\end{equation}
For any $y_l\in\mathbb R$ and $x_0\leq\mathcal X(x_0,y_l)$,
\begin{equation}
 \label{ineq:Xilevel_der_strict}
1<\frac\partial{\partial x_0}\varXi_{min}(x_0,y_0)<+\infty,\qquad
1<\frac\partial{\partial x_0}\varXi(x_0,y_0)<+\infty.
\end{equation}
\end{proposition}
\begin{proof}
We have only to prove that the derivatives are strictly positive
as the rest is proved in Proposition~\ref{Prop:lim_der_x}. Consider
the functions $\mathcal X(x_0,y_0)$ and $\varXi(x_0,y_0)$. According
to Corollary~\ref{Cor:5.2} these functions are mutually inverse.
Proposition~\ref{Prop:4.1} and Corollary~\ref{Cor:5.5} imply
that they are smooth for all values of $x_0$ and $y_l$. Thus, their
derivatives cannot vanish. The proof for $\mathcal X_{min}(x_0,y_0)$ and
$\varXi_{min}(x_0,y_0)$ is similar with the help of the functions
$\varXi_+(x_0,y_0)$ and $X_-(x_0,y_0)$, respectively
(see Remark~\ref{Rem:XX_-}).
\end{proof}
\begin{definition}
 \label{Def:Lipschitz}
$$
L_{\Xi_{min}}:=\underset{y_l\in\mathbb R,\,x_0\leq\mathcal X(0,y_l)}\sup\;
\frac\partial{\partial x_0}\varXi_{min}(x_0,y_0),\qquad
L_{\Xi}:=\underset{y_l\in\mathbb R,\,x_0\leq\mathcal X(0,y_l)}\sup\;
\frac\partial{\partial x_0}\varXi(x_0,y_0).
$$
\end{definition}
\begin{proposition}
 \label{Prop:XiLipschitz}
The numbers $L_{\Xi_{min}}$ and $L_{\Xi}$ are finite and satisfy
the following inequalities:
$$
1<L_{\Xi_{min}}\leq L_{\Xi}<+\infty.
$$
Moreover, for $x_0<\hat x_0\leq X(0)$,
\begin{equation}
 \label{ineq:Xi_min_Lipschitz}
0<\hat x_0-x_0<\varXi_{min}(\hat x_0,y_0)-\varXi_{min}(x_0,y_0)<
L_{\Xi_{min}}(\hat x_0-x_0),
\end{equation}
\begin{equation}
 \label{ineq:XiLipschitz}
0<\hat x_0-x_0<\varXi(\hat x_0,y_0)-\varXi(x_0,y_0)<L_{\Xi}(\hat x_0-x_0).
\end{equation}
\end{proposition}
\begin{proof}
Consider $\frac\partial{\partial x_0}\mathcal X(x_0,y_0)$. It is a
bounded continuous function in the whole left half-plane $x_0\leq0$.
It follows from Proposition~\ref{Prop:lim_der_x}
(see Equation~(\ref{eq:der_limit})), that the infimum of this function, which
is strictly less than $1$, is achieved at some finite point of the half-plane.
This infimum is a positive number, as
$\frac\partial{\partial x_0}\mathcal X(x_0,y_0)>0$ at any finite point
according to Proposition~\ref{Prop:der_strict}. Exactly the same conclusion
is valid for the domain
$(x_0,y_l): y_l\in\mathbb R,\;x_0\leq\mathcal X(0,y_l)$. In this domain
$1>\inf\frac\partial{\partial x_0}\mathcal X(x_0,y_0)=1/L_{\Xi}>0$.
Now from Inequalities~(\ref{ineq:5.2}) follows that
$1\leq L_{\Xi_{min}}\leq L_{\Xi}<+\infty$. The fact that $L_{\Xi_{min}}\neq1$
is clear from the first Inequality~(\ref{ineq:Xilevel_der_strict}).
The uniform bounds for the derivatives imply the Lipschitz estimates
in Inequalities~(\ref{ineq:Xi_min_Lipschitz}) and (\ref{ineq:XiLipschitz}).
\end{proof}
\begin{remark}
In \S~\ref{sect:main}, Corollary~\ref{Cor:smoothness} and
Remark~\ref{Rem:smoothness}, we prove that the functions $X(x_0)$ and
$\Xi(x_0)$ are smooth and in Corollary~\ref{Cor:der_X_Xi_bounded} state
the properties of their derivatives similar to those for the corresponding
level functions proved above.
\end{remark}
\section{Geometry of Integral Curves}
 \label{sect:main}
\begin{definition}
 \label{Def:6.1}
The graphs of two solutions are said to intersect at the point at infinity
iff these solutions have a common pole, i.e., if their
intervals of existence have a common end point.
\end{definition}
\begin{theorem}
 \label{Thm:intersections}
The graphs of any two (different) solutions of Equation {\rm(\ref{p1})} have
at most two points of intersection, counting the points at infinity,
with non-positive abscissas.
\end{theorem}
\begin{proof}
Given two intersecting solutions, we count intersection points starting
from the one with abscissa closest to the origin and counting to the left.
The first intersection point exists. To prove this, denote by $x_0\leq0$
the supremum of abscissas of the set of intersection points. If $x_0$
belongs to the interval of existence of both solutions, then it is also
an intersection point by continuity.
Suppose that $x_0$ is a pole of one of the solutions, then it should be
also a pole of the other. From the Laurent expansion~(\ref{laurent})
it follows that any two solutions with the pole at the same point do not
intersect in some neighbourhood of this point. Therefore in the latter case
$x_0$ is an isolated intersection point at infinity.

Now, we have to consider several different cases:
\begin{enumerate}
\item
The solutions have a common pole: one solution has the interval of existence
to the left of the pole and the other to the right. The theorem is obviously
true.
\item
The solutions have a common pole at the right end of their intervals
of existence. Then they do not intersect until both of them achieve their
minimum values (see Lemma~\ref{Lem:3.1}). One of these minimum values is
larger than the other and we call the corresponding solution the
upper-solution. The finite point of intersection may occur
at the minimum value of the upper solution or after it.
In the first case we use Lemma \ref{Lem:4.3} and in the second
Lemma \ref{Lem:4.1} to prove that there are no other intersection points.
\item
The first point of intersection is finite and solutions at the intersection
point have positive derivatives. The respective location of graphs of the
solutions after the first intersection is the same as in the second item.
\item
At the first intersection point one solution has a positive derivative
and the other non-positive. Then the second intersection point may occur
and at the latter point, both solutions have negative derivative and we use
Lemma \ref{Lem:4.1} to prove that there are no other intersection points.
\item
At the first intersection point one solution has non-positive derivative
and the other a negative one. These solutions have no other intersection
points according to Lemma~\ref{Lem:4.1}.
\end{enumerate}
\end{proof}
\begin{lemma}
 \label{Lem:MN}
Let $y_k(x)$ for $k=1,2,3$, be solutions of Equation {\rm(\ref{p1})} whose
intervals of existence contain the segment $[\alpha,\beta]\in(-\infty,0]$.
Suppose that: {\rm(i)} $y_1(x)>y_3(x)$ and $y_2(x)>y_3(x)$ for
$x\in(\alpha,\beta)$; and {\rm(ii)} $y_1(x)\geq y_3(x)$ and
$y_2(x)\geq y_3(x)$ at $x=\alpha,\beta$. Then the graphs of solutions
$y_1(x)$ and $y_2(x)$ have no more than one point of intersection with
abscissa $x\in[\alpha,\beta]$.
\end{lemma}
\begin{proof}
If the conditions $y_1(x)>y_3(x)$ and $y_2(x)>y_3(x)$ hold for all
$x\in[\alpha,\beta]$, then this statement is a particular case of a general
lemma proved by Moore and Nehari (Lemma 2 in \cite{MN}).
Therefore, we consider only the case when all three solutions intersect
at $\alpha$ and/or $\beta$, since for all other possibilities our Lemma
follows from the one by Moore and Nehari \cite{MN} by applying the latter
to possibly a smaller segment $[\alpha',\beta']\subset[\alpha,\beta]$.

Suppose $y_k(\alpha)=y_0$ for $k=1,2,3$ and solutions $y_1(x)$ and
$y_2(x)$ have one more point of intersection $\gamma\in(\alpha,\beta)$.
For definiteness, let $y_1(x)>y_2(x)$ for $x\in(\alpha,\gamma)$.
In that case for $\epsilon>0$ consider a family of solutions
$y_\epsilon(x)$ with the initial data
$y_\epsilon(\alpha)=y_3(\alpha)-\epsilon$,
$y'_\epsilon(\alpha)=y'_3(\alpha)$. We claim that for all rather small
$\epsilon$, $y_2(x)>y_\epsilon(x)$ for $x\in(\alpha,\gamma)$. Otherwise,
there is a sequence $\epsilon_n\to0$ such that solutions
$y_{\epsilon_n}(x)$ cross $y_2(x)$ at least twice, at the points:
$x_2>x_1\in(\alpha,\gamma)$, say, as $y_{\epsilon_n}(x)$ have to be close
to $y_3(x)$ and they cannot be tangent to $y_2(x)$.
Before crossing $y_2(x)$ they have to cross $y_3(x)$ as $y_3(x)$ separated
them from $y_2(x)$. Therefore, on the segment $[x_1,x_2]$ we have
$y_2(x)>y_3(x)$ and $y_{\epsilon_n}(x)>y_3(x)$ for all large enough $n$.
This contradicts the standard Moore--Nehari Lemma (see the first
paragraph of this proof) as solutions $y_2(x)$ and $y_{\epsilon_n}(x)$
have two points of intersection in $[x_1,x_2]$. Thus, the triple of
solutions $y_1(x)$, $y_2(x)$ and $y_\epsilon(x)$ for any sufficiently small
$\epsilon$ satisfies the standard Moore--Nehari Lemma on $[\alpha,\gamma]$.
This proves that $\gamma\notin(\alpha,\beta)$, since
$y_1(\alpha)=y_2(\alpha)$.

We are left with the case when
$y_1(x)>y_2(x)>y_3(x)$ for $x\in(\alpha,\beta)$ and
$y_1(x)=y_2(x)=y_3(x)$ at $x=\alpha, \beta$. In this case consider a
solution $y_4(x)$ with $y_4(\alpha)=y_k(\alpha)$ and
$y'_2(\alpha)>y'_4(\alpha)> y'_3(\alpha)$. This solution cannot cross
$y_2(x)$ for $x\in(\alpha,\beta)$ as follows from the previous paragraph
(for a triple $y_2(x)\geq y_4(x)\geq y_3(x)$).
Suppose $y_4(x)$ crosses $y_3(x)$ at $x=\gamma<\beta$. In this case
$y_4(x)<y_3(x)$ for $x\in(\gamma,\beta]$ as it follows from
Theorem~\ref{Thm:intersections}. Therefore, we arrive at the
situation considered in the second paragraph of this proof but for the
triple $y_1(x)\geq y_2(x)\geq y_4(x)$ for $x\in[\alpha,\beta]$.
The only situation that is left is that for all slopes $y'_4(\alpha)$:
$y'_2(\alpha)\geq y'_4(\alpha)\geq y'_3(\alpha)$ solution $y_4(x)$ is
such that $y_4(\beta)=y_1(\beta)$, i.e., $y_4(\beta)$ is independent of
the initial slope $y'_4(\alpha)\in[y'_3(\alpha),y'_2(\alpha)]$. Taking
into account that $y_4(\beta)$ is actually a complex analytic function of
$y'_4(\alpha)$ we get that $y_4(\beta)=y_1(\beta)$ for all
$y'_4(\alpha)\in\mathbb C$. That cannot be the case as for
$y'_4(\alpha)\to+\infty$ $y_4(x)$ has a pole $x_p\to\alpha$.

Consider now the case $y_k(\beta)=y_1$ for $k=1,2,3$. Without
loss of generality we suppose now that $y_1(\alpha)>y_2(\alpha)>y_3(\alpha)$.
According to the standard Moore--Nehari lemma $y_1(x)$ and $y_2(x)$ has
at most one intersection point $\gamma\in(\alpha,\beta)$.
Then we have
$y_2(x)>y_1(x)>y_3(x)$ for $x\in(\gamma,\beta)$ and
$y_1(\gamma)=y_2(\gamma)>y_3(\gamma)$. Now consider an auxiliary solution
$y_4(x)$ with initial data $y_4(\gamma)=y_1(\gamma)=y_2(\gamma)$ and
$y'_2(\gamma)>y_4'(\gamma)>y_1'(\gamma)$. The graph of this solution cannot
cross either $y_1(x)$ or $y_2(x)$ on the interval $(\gamma,\beta)$
(by the standard Moore--Nehari lemma). Thus, $y_4(\beta)=y_1$ and we arrive
to the situation considered in the previous paragraph.
\end{proof}
\begin{remark}
 \label{Rem:projective_MN}
In Lemma~\ref{Lem:MN} ``the point of intersection'' is supposed to be
a finite point of intersection. However this Lemma can be generalized
to include the points of intersection at infinity in the sense of
Definition~\ref{Def:6.1}. In the latter case we have to allow $\alpha$
($\beta$) to be the left (right) bound of the interval of existence of
one of the functions $y_k(x)$, $k=1,2,3$. Otherwise all points of
intersection are finite. If $\alpha$ ($\beta$) is the left (right) bound
of the function $y_3(x)$, then inequality $y_k(x)\geq y_3(x)$ is violated at
$\alpha$ ($\beta$): we formally have
$y_1(\alpha)=y_2(\alpha)=y_3(\alpha)=+\infty$ (the same condition at
$\beta$). As proved below in this situation we have to
count the point at infinity with abscissa $\alpha$ ($\beta$) as the
intersection point. We call Lemma~\ref{Lem:MN} equipped with the above
extended conditions as {\bf Projective Lemma~\ref{Lem:MN}}. The proof is
given below after Theorem~\ref{Thm:uniqueness} and
Lemma~\ref{Lem:approximation}.
\end{remark}
\begin{theorem}
 \label{Thm:uniqueness}
For any $x_0\leq0$ there is a unique solution with the interval of
existence $\left(X(x_0),x_0\right)$.
\end{theorem}
\begin{proof}
Suppose there are at least two such solutions $Y_1(x)$ and $Y_2(x)$.
Denote by $Y_{lk}$, $k=1,2$ their respective minimum values. Suppose
$Y_{l1}>Y_{l2}$. From Lemmas \ref{Lem:3.1}, \ref{Lem:4.1}, and
\ref{Lem:4.3} it follows that solutions do not intersect. Therefore,
$Y_1(x)>Y_2(x)$ for $x\in\left(X(x_0),x_0\right)$.
Suppose there is a solution, $Y_0(x)$ with the pole
at $x_0$ with the minimum value $Y_{l0}\in(Y_{l1},Y_{l2})$, such that the
left bound of its interval of existence is larger than $X(x_0)$. Therefore,
the graph of $Y_0(x)$ crosses at some finite point, $z$, the graph of
$Y_1(x)$, however, it does not cross the graph of $Y_2(x)$. This crossing
point, $z$, may occur before, after, or exactly at the minimum of
$Y_1(x)$.

Consider now a solution $y(x)$ with the finite initial value
$y(x_0)=y_0>Y_{l1}$ and minimum value $Y_{l1}$ at $x=x_{min}$. Denote by
$x_{1min}$ the minimum of $Y_1(x)$. Recall that
$x_{min}\searrow x_{1min}$ as $y_0\nearrow+\infty$. If $x_{1min}<z$, then
increasing, if necessary $y_0$, we can assume that $x_{1min}<x_{min}<z$.
In this case, the graph of $y(x)$ first crosses the graph of $Y_2(x)$ at
some point with abscissa $\beta$. Then it crosses the graph of $Y_0(x)$ at
some point with abscissa $x_1<\beta$. At $\beta$ and $x_1$ both solutions
$y(x)$ and $Y_2(x)$ or, respectively, $y(x)$ and $Y_0(x)$ have positive
derivatives. If one of these crossings does not occur, say $x_1$, then the
graph of $y(x)$ for $x\in[x_{min},x_0]$ will be located below the graph of
$Y_0(x)$. Therefore, the minimum point of $y(x)$ will be on finite distance
from the minimum point of $Y_1(x)$, despite $y_0\nearrow+\infty$.
At point $z$, however, $Y'_0(z)<0$, as $Y_0(x)$ and
$Y_1(x)$ have the same pole $x_0$: they cannot, according to
Lemma \ref{Lem:3.1}, intersect both having positive derivatives. Solution
$Y_0(x)$ has the minimum value $Y_{l0}<Y_{l1}$, therefore, to have a
negative derivative at $z$ it should cross the level $Y_{l1}$ at some
point with abscissa $z_1>z$ with $Y'_0(z_1)<0$. Thus, the graph of $y(x)$
crosses the graph of $Y_0(x)$ at a point with abscissa $x_2\in(z,z_1)$.
Moreover, $x_2\neq x_1$ since $Y'_0(x_1)Y'_0(x_2)<0$.

In case $z\leq x_{1min}$ the second crossing point, with abscissa $x_2$,
of $y(x)$ and $Y_0(x)$ occurs in the domain above the graph of $Y_1(x)$,
where both derivatives $Y'_0(x)$ and $y'(x)$ are negative. Note, that the
first intersection point always belongs to the domain below the graph of
$Y_1(x)$. More precisely, in this case, the graph of $y(x)$ crosses the
graph of $Y_1(x)$ at some point with abscissa located in the interval
$(x_{1min},x_{min})$. It cannot cross $Y_1(x)$ anymore according to
Lemma~\ref{Lem:4.2}. However, for any $\epsilon_1,\epsilon_2>0$ for all
large enough $y_0$ we have $0<y(x)-Y_1(x)<\epsilon_1$ for
$x\in[X(x_0)+\epsilon_2,x_{1min}]$. For all rather small $\epsilon_2$,
$z\in[X(x_0)+\epsilon_2,x_{1min}]$, thus $0<y(z)-Y_0(z)<\epsilon_1$.
On the other hand, $Y_0(x)$ has a pole at
$\alpha\in[X(x_0)+\epsilon_2,x_{1min}]$;
therefore the difference $y(x)-Y_0(x)<0$ as $x\nearrow\alpha$.
This proves existence of the second intersection point with abscissa $x_2$.

Thus, in all cases when $Y_0(x)$ has its interval of existence less
than $(X(x_0),x_0)$, we constructed three solutions: $y(x)$,
$Y_0(x)$, and $Y_2(x)$, such that $y(x)>Y_2(x)$ and $Y_0(x)>Y_2(x)$ for
$x\in[\alpha,\beta-\epsilon]$, with small enough $\epsilon$,
and the graphs of solutions $y(x)$ and $Y_0(x)$ have two (finite)
intersection points with abscissas
$x_1$ and $x_2\in[\alpha,\beta-\epsilon]$.
This contradicts Lemma~\ref{Lem:MN}.

Therefore, all solutions with the pole at $x_0$ and minimum value
$Y_{l0}\in[Y_{l2},Y_{l1}]$ have the interval of existence $(X(x_0),x_0)$
and their graphs do not intersect. Recall that $\mathcal X(x_0,y_l)$
is a (complex) analytic function of the pole parameter $c$; it coincides
with the function $x_1=x_1(x_0,c_0)$ for $c=c_0$ introduced in
Proposition~\ref{Prop:analyticity}. The fact that
$\mathcal X(x_0,y_l(c))=X(x_0)$ on some segment $c\in[c_1,c_2]$ where
$c_1<c_2\in\mathbb R$ implies that $\mathcal X(x_0,y_l(x_0,c))=X(x_0)$
for all $c\in\mathbb C$ and, in particular, for $c\in\mathbb R$. This
contradicts the fact that
$\mathcal X(x_0,y_l(x_0,c))\underset{c\to\pm\infty}{\to}x_0$ as follows
from Corollary~\ref{Cor:3.3} and Proposition~\ref{Prop:4.3}.
\end{proof}
\begin{corollary}
 \label{Cor:smoothness}
The function $X(x_0)$ is smooth. There exists a continuous function
$\hat y_l(x_0)$, such that
$
X'(x_0)=\frac\partial{\partial x_0}\mathcal X(x_0,\hat y_l(x_0)).
$
\end{corollary}
\begin{proof}
In Proposition~\ref{Prop:existence} it was proved that
for any $x_0\leq0$ there exists $\hat y_l$ such that
$X(x_0)=\mathcal X(x_0,\hat y_l)$. Now by Theorem~\ref{Thm:uniqueness} we
can define a single valued function $\hat y_l(x_0)$ with the property
$X(x_0)=\mathcal X(x_0,\hat y_l(x_0))$. The function $\hat y_l(x_0)$ is
continuous. Actually, consider a sequence $x_n\to x_0$. The sequence
$\hat y_l(x_n)$ is bounded, since for an unbounded subsequence
$\hat y_l(x_{n_k})$ we would have
$\mathcal X(x_{n_k},\hat y_l(x_{n_k}))-x_{n_k}\to0$
(see Proposition~\ref{Prop:4.3}) rather than to a finite value about
$2C/|x_0|^{1/4}$ (cf. Theorem~\ref{Thm:1}). If we suppose that there is
a subsequence $\hat y_l(x_{n_m})$ that does not converge to $\hat y_l(x_0)$,
then we get that a subsequence $\hat y_l(x_{n_{m_k}})$ converges to some
number $\hat y\neq\hat y_l(x_0)$. The functions $\mathcal X(x_0,y_l)$ are
smooth with respect to both variables (see Proposition~\ref{Prop:4.1}),
therefore $X(x_0)=\mathcal X(x_0,\hat y_l(x_0))=\mathcal X(x_0,\hat y)$
in contradiction with Theorem~\ref{Thm:uniqueness}.

In the following part of the proof we consider, as at the end of the
proof of Theorem~\ref{Thm:uniqueness}, the analytic functions
$\mathcal X(x_0,y_l(x_0,c))$, which we denote for simplicity as
$\mathcal X(x_0,c)$. Note that from the previous paragraph via the smooth
bijection $c(x_0,y_l)$ (see Corollary~\ref{Cor:3.3}), we have a continuous
function $\hat c(x_0)=c(x_0,\hat y_l(x_0))$ with the property
$X(x_0)=\mathcal X(x_0,\hat c(x_0))$.
Consider variation of the function $X(x_0)$ under the shift $x_0+\Delta x_0$,
\begin{equation}
 \label{eq:variation}
\Delta X(x_0)=\frac\partial{\partial x_0}
\mathcal X(\xi,\hat c(x_0+\Delta x_0))\Delta x_0+
\frac\partial{\partial c}\mathcal X(x_0,\eta)\Delta\hat c,
\end{equation}
where the numbers $\xi\in(x_0,x_0+\Delta x_0)$,
$\eta\in(\hat c(x_0+\Delta x_0),\hat c(x_0))$, and
$\Delta\hat c=\hat c(x_0+\Delta x_0)-\hat c(x_0)$, where for definiteness
we assume $\Delta x_0>0$ (the case $\Delta x_0<0$ differs only
by an obvious modification of the notation).
Since $\hat c(x_0)$ is a continuous function there is a point
$x_0+\Delta x_1\in(x_0,x_0+\Delta x_0)$ such that
$\hat c(x_0+\Delta x_1)=\eta$. From the definition of the
function $\hat c(x_0)$ it follows that:
$
\frac\partial{\partial c}\mathcal X(x_0,\hat c(x_0))=0
$
and
$
\frac\partial{\partial c}\mathcal X(x_0+\Delta x_1,\eta)=0.
$
The last equation can be rewritten as
$$
\frac{\partial^2}{\partial x_0\partial c}\mathcal X(\xi_1,\eta)\Delta x_1+
\frac\partial{\partial c}\mathcal X(x_0,\eta)=0,
$$
where $\xi_1\in(x_0,x_0+\Delta x_1)$. Now substituting
$
\frac\partial{\partial c}\mathcal X(x_0,\eta)
$
from the last equation into Equation~\ref{eq:variation} and taking into
account that $0<\Delta x_1<\Delta x_0$, $\hat c(x_0)$ is a continuous
function, and $\mathcal X(x_0,c)$ is analytic in both variables, and
taking a limit $\Delta x_0\to0$, we arrive at the announced result.
\end{proof}
\begin{remark}
 \label{Rem:smoothness}
Only notational modifications to the above proof are required to prove
the smoothness of $\Xi(x_0)$. If the uniqueness Conjectures~\ref{Con:2}
and \ref{Con:3} of \S~\ref{sect:numerics} are justified literally the same
proof works to establish smoothness of the functions $X_{min}(x_0)$ and
$\Xi_{min}(x_0)$.
\end{remark}
\begin{corollary}
 \label{Cor:der_X_Xi_bounded}
For all $x_0\leq0$:
$$
0<X'(x_0)<1,\qquad
\underset{x_0\to-\infty}\lim X'(x_0)=1,
$$
For all $x_0\leq X(0)$:
$$
0<\Xi'(x_0)\leq L_{\Xi},\qquad
\underset{x_0\to-\infty}\lim \Xi'(x_0)=1.
$$
\end{corollary}
\begin{proof}
The inequalities and the limit for $X'(x_0)$ follows from the formula for this
derivative in terms of the partial $x_0$-derivative of the corresponding level
function (see Corollary~\ref{Cor:smoothness}), and
Inequalities~(\ref{ineq:Xlevel_der_strict}) and double
limit~(\ref{eq:der_limit}), respectively.

The limit of $\Xi'(x_0)$ as $x_0\to-\infty$ is the same as the limit of
$X'(x_0)$ by virtue of the formula $X'(\Xi(x_0))\Xi'(x_0)=1$ and the fact
that $\Xi(x_0)\to-\infty$ when $x_0\to-\infty$. To prove the inequalities we
recall that, as mentioned in Remark~\ref{Rem:smoothness}, there is an
analogous representation for $\Xi'(x_0)$ as the one for $X'(x_0)$:
$\Xi'(x_0)=\frac\partial{\partial x_0}\Xi(x_0,\tilde y_l(x_0))$, where
$\tilde y_l(x_0)$ is a continuous function. To finish the proof we refer
to Propositions~\ref{Prop:der_strict}, \ref{Prop:XiLipschitz}.
\end{proof}
\begin{remark}
We expect that the functions $X'(x_0)$ and  $\Xi'(x_0)$ are monotonic.
\end{remark}
\begin{lemma}
 \label{Lem:approximation}
Any solution $y(x)$ of Equation~{\rm(\ref{p1})} with the interval of
existence $(a,b)$ can be presented $($non-uniquely$)$ as a limit of
a sequence of solutions $y_n(x)$ that are regular at point $a$ $($or $b)$
with initial data $y_n(a)\to+\infty$ and $y'_n(a)\to+\infty$
$(y_n(b)\to+\infty$ and $y'_n(b)\to-\infty)$. More precisely, this means
that for any $\epsilon,\delta>0$ and $x\in[a+\delta,b-\delta]$, there is
$N$ such that for all $n\geq N$: $|y(x)-y_n(x)|<\epsilon$.
\end{lemma}
\begin{proof}
{}From Theorem \ref{Thm:3.1} one deduces that the minimum of any solution
with the pole at $b$ and minimum value $y_l$ can be presented as
$\mathcal X_{min}(b,y_l)$. Using Definition~\ref{Def:3.1} and
Proposition~\ref{Prop:3.1} one constructs, for any sequence
$y_n(b)\to+\infty$ as $n\to+\infty$, a sequence of solutions with
the properties stated in this Lemma. These properties follow
from Proposition~\ref{Prop:3.2}, Corollary~\ref{Cor:3.1}, and standard
facts known about dependence of the solution on the initial data.
Existence of the sequence under discussion at point $a$ follows by
analogous arguments from the results collected in \S~\ref{sect:5}.
\end{proof}
\begin{remark}
In fact one can suggest many other constructions of the sequences
satisfying the requirements of Lemma~\ref{Lem:approximation}. Say,
consider any point $(x_0,y(x_0))$ of the graph of $y(x)$ with
$y'(x_0)>0$ (or $y'(x_0)<0$).
Take any sequence of numbers $a_n\to y'(x_0)$. Consider a sequence of
solutions $y_n(x)$ with initial data $y_n(x_0)=y(x_0)$ and $y'_n(x_0)=a_n$.
All these solutions are regular at $b$ (or, respectively, $a$), this
follows from Lemma~\ref{Lem:3.1} (Lemma~\ref{Lem:4.1}), and have the
properties stated in Lemma~\ref{Lem:approximation}.
\end{remark}
\begin{proof}[Proof of Projective Lemma~\ref{Lem:MN}]
We have to consider several cases:

(1) We have three solutions $y_k(x)$, with the graphs
located as in Lemma \ref{Lem:MN}. Solutions $y_1(x)$ and $y_2(x)$ have
a pole at $x=\beta$ and one more finite point of intersection at
$x_1<\beta$.
Suppose, first, that $y_3(\beta)$ is finite. Then, consider two sequences
of solutions $\hat y_k^n(x)$, $k=1,2$ and $n=1,2,\ldots$, regular at
$x=\beta$ and approximating respectively the given solutions $y_k(x)$,
in the sense of Lemma~\ref{Lem:approximation}. To shorten the proof, we
suppose that $y_1^n(\beta)=y_2^n(\beta)\to+\infty$: it follows from
Proposition~\ref{Prop:3.1} that it is always possible. (However, this
assumption is not important and one can modify this proof by considering
more general sequences of solutions and using only
Lemma~\ref{Lem:approximation}.) For any $n$, solutions $y_1^n(x)$ and
$y_2^n(x)$ have one point of intersection $\beta$. Then increasing $n$,
we see that there exist $N\in\mathbb N$ such that every solution
$y_1^n(x)$ will be intersecting with all solutions $y_2^m(x)$,
$n,m\geq N$ near the finite point of intersection $x_1$. We apply
Lemma~\ref{Lem:MN} to get a contradiction.

Suppose that $y_3(x)$ also has a pole at $\beta$.
We first approximate it by a solution $\hat y_3(x)$ finite at $x=\beta$
such that the graph of the latter solution lies below the graphs of both
solutions $y_1(x)$ and $y_2(x)$. After that we have a situation considered
in the previous paragraph.

(2) In the case when $y_1(x)$ and $y_2(x)$ have an intersection at
the point of infinity at the left bound of their interval of existence and
a finite intersection point to the right of it, the proof is analogous to
that in case (1).

(3) Consider the situation when solutions $y_1(x)$ and $y_2(x)$
have two points of intersection at infinity with the abscissas $\alpha$
and $\beta$, i.e., they have the common interval of existence
$(\alpha,\beta)$. The interval of existence of $y_3(x)$ in this case
includes $(\alpha,\beta)$. According to Theorem~\ref{Thm:intersections}
$y_1(x)$ and $y_2(x)$ have no other finite points of intersection.
Suppose that $y_1(x)>y_2(x)$. A solution $y_0(x)$ with the pole at
$\alpha$ or $\beta$ that have a minimum value (a pole parameter)
between minimum values of $y_1(x)$ and $y_2(x)$ (pole parameters
corresponding to $y_1(x)$ and $y_2(x)$)
does not intersect $y_3(x)$ before it crosses $y_2(x)$. If such
intersection happens then for the triple, $y_0(x)$, $ y_2(x)$ and
$y_3(x)$, we have a situation considered in the previous paragraphs.
Therefore, for the triple $y_0(x)$, $ y_1(x)$ and $y_2(x)$ we
meet a situation considered in the last paragraph of the proof of
Theorem~\ref{Thm:uniqueness}.
\end{proof}
\begin{theorem}
 \label{Thm:dual}
Given any two (different) points, including the points at infinity,
with non-positive abscissas, there are no more than two solutions
of Equation~{\rm(\ref{p1})} whose graphs pass through these points.
\end{theorem}
\begin{proof}
If we have three solutions whose graphs pass through the given points,
then these graphs have no other points of intersection according to
Theorem \ref{Thm:intersections}. Therefore, for $x\in[\alpha,\beta]$,
where $\alpha$ and $\beta$ are abscissas of the given points we have
that $y_1(x)>y_2(x)>y_3(x)$ and we can apply
Projective Lemma~\ref{Lem:MN}.
\end{proof}
\begin{remark}
Theorem~\ref{Thm:dual} is dual to Theorem~\ref{Thm:intersections}.
\end{remark}
\begin{definition}
For $x_0\leq0$ (respectively, $x_0\leq X(0)$), any solution having $x_0$
as the right (left) bound of its interval of existence is called the
{\it left (right) solution for $x_0$}.\\
For any $x_0\leq0$ the {\em left maximum solution for $x_0$} is
the solution from Theorem~{\rm\ref{Thm:uniqueness}}. For any
$x_0\leq X(0)$ the {\em right maximum solution for $x_0$}
is the left maximum solution for the point $X^{-1}(x_0)$.
\end{definition}
\begin{corollary}
Let $x_0\leq X(0)$ and $\{(x_0,a)\}$ be the set of intervals of existence
of the right solutions for $x_0$. Then $X^{-1}(x_0)=\sup\,a\leq0$. The
right maximum solution for $x_0$ is the unique solution with the interval
of existence $(x_0,X^{-1}(x_0))$.
\end{corollary}
\begin{proof}
Since $X(X^{-1}(x_0))=x_0$, the first statement follows from
Proposition~\ref{Prop:existence}. The uniqueness is proved in
Theorem~\ref{Thm:uniqueness}.
\end{proof}
\begin{theorem}
 \label{Thm:monotonic_up}
For $x_0\leq0$ (respectively, $x_0\leq X(0)$), consider left (right)
solutions for $x_0$ with the minimum values $y_l\geq y_{l\mu}$, where
$y_{l\mu}$ is the minimum value of $y_\mu(x)$, the left (right) maximum
solution for $x_0$. Any two such solutions, $y_1(x)$ and $y_2(x)$, have
no finite points of intersection. Moreover, the left (right) bound of
their interval of existence is a continuous strictly monotonically
increasing (decreasing) function of $y_l$ mapping $[y_{l\mu},+\infty)$
onto $[X(x_0),x_0)$ (\,$(x_0,X^{-1}(x_0)]$):
$$
y_{l_2}>y_{l_1}\geq y_{l\mu}\quad\Rightarrow\quad
\mathcal X(x_0,y_{l_2})>\mathcal X(x_0,y_{l_1})\quad
\left(\mathcal\varXi(x_0,y_{l_2})<\mathcal\varXi(x_0,y_{l_1})\right).
$$
\end{theorem}
\begin{proof}
This is an immediate consequence of Projective Lemma~\ref{Lem:MN} for the
triple of solutions: $y_1(x)$, $y_2(x)$, and $y_\mu(x)$. The only point
of intersection that is allowed by this Lemma is the point of intersection
at infinity with abscissa $x_0$. The properties of the monotonic function
are stated in Propositions~\ref{Prop:analyticity} and \ref{Prop:4.3}.
\end{proof}
\begin{theorem}
 \label{Thm:monotonic_down}
For $x_0\leq0$ (respectively, $x_0\leq X(0)$), consider left (right)
solutions for $x_0$ with the minimum values $y_l\leq y_{l\mu}$, where
$y_{l\mu}$ is the minimum value of $y_\mu(x)$, the left (right) maximum
solution for $x_0$. Any two such solutions, $y_1(x)$ and $y_2(x)$, have
one finite point of intersection. The left (right) bound of their interval
of existence is a continuous strictly monotonically decreasing (increasing)
function of $y_l$ mapping $(-\infty,y_{l\mu}]$ onto $[X(x_0),x_0)$
(\,$(x_0,X^{-1}(x_0)]$):
$$
y_{l\mu}\geq y_{l_2}>y_{l_1}\quad\Rightarrow\quad\mathcal X(x_0,y_{l_2})<
\mathcal X(x_0,y_{l_1})\quad
\left(\mathcal\varXi(x_0,y_{l_2})>\mathcal\varXi(x_0,y_{l_1})\right).
$$
\end{theorem}
\begin{proof}
Suppose $y_{l_\mu}\geq y_{l_2}>y_{l_1}$, but solutions $y_1(x)$ and $y_2(x)$
have no finite point of intersection. Let $\alpha$ be an abscissa of the
point of intersection of $y_2(x)$ with $y_\mu(x)$. Then, $y_\mu(x)>y_1(x)$
and $y_2(x)>y_1(x)$ for $x\in[\alpha,x_0)$ (or $x\in(x_0,\alpha]$), the
triple $y_1(x)$, $y_2(x)$, and $y_\mu(x)$ have one point of intersection at
infinity with abscissa $x_0$, and $y_1(x)$ intersects $y_2(x)$ at $\alpha$.
Therefore, we arrive at the contradiction by applying
Projective Lemma~\ref{Lem:MN} to the triple: $y_\mu(x)$, $y_2(x)$, and
$y_1(x)$ on the segment $x\in[\alpha,x_0]$ (or $x\in[x_0,\alpha]$). The
properties of the monotonic function are stated in
Propositions~\ref{Prop:analyticity} and \ref{Prop:4.3}.
\end{proof}
\section{Boundary Value Problems}
 \label{sect:boundary}
\begin{theorem}
 \label{Thm:7.1}
Let $x_1<x_0\leq0$. If\\
\indent
$1$. $x_1<X(x_0)$, then there are no solutions with the interval of
existence $(x_1,x_0)$;\\
\indent
$2$. $x_1=X(x_0)$, then there is one solution with the interval of
existence $(x_1,x_0)$;\\
\indent
$3$. $X(x_0)<x_1<x_0$, then there are two solutions with the interval
of existence $(x_1,x_0)$.
\end{theorem}
\begin{proof}
The first statement follows from Definition~\ref{Def:2.2} and
Theorem~\ref{Thm:1}.
The second statement is proved in Theorem~\ref{Thm:uniqueness}.
The third statement follows immediately from
Theorems~\ref{Thm:monotonic_up} and \ref{Thm:monotonic_down}.
\end{proof}
\begin{remark}
Theorem~\ref{Thm:7.1} can be viewed as a statement about a boundary value
problem with infinite boundary conditions. To study finite boundary value
problems we introduce the function $Z(x_0,y_0,y^0)$ which plays the role of
of the function $X(x_0)$ in Theorem~\ref{Thm:7.1}. The properties
of $Z(x_0,y_0,y^0)$ are studied with the help of an auxiliary
``level-function'' $\mathcal Z(x_0,y_0,y^0;y_l)$.
\end{remark}
\begin{definition}
Let $x_0\leq0,\;y^0\geq y_0>y_l$. Denote by $y_{\pm}(x)$ the solution
with the initial value $y_\pm(x_0)=y_0$ and minimum value
$y_{\pm}(x_{min\pm})=y_l$, where $x_{min-}<x_0$ and $x_{min+}>x_0$.
Define the function,
$\mathcal Z_\pm=\mathcal Z_\pm(x_0,y_0,y^0;y_l)$ as the (unique) solution
of the equation $y_\pm(\mathcal Z_\pm)=y^0$ satisfying inequalities:
$\mathcal Z_-<x_0$ and $\mathcal Z_+>x_0$. We also put
$\mathcal Z_\pm(x_0,y_0,y^0;y_0)=x_0$. (In case $x_0>X(0)$ the
definition of the function $\mathcal Z_+$ is limited to those solutions
whose graphs do not cross the axis of ordinates below $y^0$.)
\end{definition}
\begin{proposition}
 \label{Prop:z-}
The function $\mathcal Z_-(x_0,y_0,y^0;y_l)$ is a smooth function of
all its variables. Moreover, if $y_0>y_l$, then:
\begin{equation}
 \label{ineq:z-}
{\mathcal X}(x_0,y_l)<\mathcal Z_-(x_0,y_0,y^0;y_l)<x_0,
\end{equation}
$$
\underset{y_0\to+\infty}\lim\mathcal Z_-(x_0,y_0,y^0;y_l)=
{\mathcal X}(x_0,y_l),\quad
\underset{y_l\to\pm\infty}\lim\mathcal Z_-(x_0,y_0,y^0;y_l)=x_0.
$$
\begin{eqnarray}
y^1>y^0&\Rightarrow&\mathcal Z_-(x_0,y_0,y^1;y_l)<
\mathcal Z_-(x_0,y_0,y^0;y_l),
\label{ineq:z-mon1}\\
y^0\geq y_{01}>y_0&\Rightarrow&\mathcal Z_-(x_0,y_{01},y^0;y_l)<
\mathcal Z_-(x_0,y_0,y^0;y_l),
\label{ineq:z-mon2}\\
0\geq x_2>x_1&\Rightarrow&0<\mathcal Z_-(x_2,y_0,y^0;y_l)-
\mathcal Z_-(x_1,y_0,y^0;y_l)<x_2-x_1.
\end{eqnarray}
\end{proposition}
\begin{proof}
The smoothness is a consequence of the implicit function theorem and
the smooth behaviour of the solution $y_-(x)$ on initial data $x_0,y_0$
and the minimum value $y_l$ (Remark~\ref{Rem:notation}).
Inequalities~(\ref{ineq:z-}) immediately follow from definitions
of the functions $\mathcal X$ and $\mathcal Z_-$.
To prove the first
limit note that $y^0>y_0\to\infty$ therefore as $y_0\to+\infty$ the
solution corresponding to the function $\mathcal Z_-(x_0,y_0,y^0;y_l)$
is approximating the one that corresponds to the function
${\mathcal X}(x_0,y_l)$ in the sense of Lemma~\ref{Lem:approximation}
(see Proposition~\ref{Prop:4.1}, Definition~\ref{Def:3.1}, and
Proposition~\ref{Prop:3.2}).
The proof of the second limit follows from Inequalities~(\ref{ineq:z-})
and Proposition~\ref{Prop:4.3}.
Inequality~(\ref{ineq:z-mon1}) is evident from the definition of
$\mathcal Z_-$.
Left Inequality~(\ref{ineq:z-mon2}) for $\mathcal Z_-$ follows from
Proposition~\ref{Prop:3.2} and Lemma~\ref{Lem:4.2}. The proof of the right
is literally the same as the one of Proposition~\ref{Prop:3.5}.
\end{proof}
\begin{proposition}
 \label{Prop:z+}
The function $\mathcal Z_+(x_0,y_0,y^0;y_l)$ is a smooth function of
all its variables. Moreover, if $y_0>y_l$, then:
$$
x_0<\mathcal Z_+(x_0,y_0,y^0;y_l)<\varXi(x_0,y_l),
$$
$$
\underset{y_0,y_1\to+\infty}\lim\mathcal Z_+(x_0,y_0,y^0;y_l)=
\varXi(x_0,y_l),\quad
\underset{y_l\to\pm\infty}\lim\mathcal Z_+(x_0,y_0,y^0;y_l)=x_0.
$$
\begin{eqnarray*}
y^1>y^0&\Rightarrow&\mathcal Z_+(x_0,y_0,y^1;y_l)>
\mathcal Z_+(x_0,y_0,y^0;y_l),\\
y^0\geq y_{01}>y_0&\Rightarrow&\mathcal Z_+(x_0,y_{01},y^0;y_l)>
\mathcal Z_+(x_0,y_0,y^0;y_l),\\
X(0,y_l)\geq x_2>x_1&\Rightarrow&0<x_2-x_1<\mathcal Z_+(x_2,y_0,y^0;y_l)-
\mathcal Z_+(x_1,y_0,y^0;y_l).
\end{eqnarray*}
\end{proposition}
\begin {proof}
The proof is very similar to the one of Proposition~\ref{Prop:z-}.
\end{proof}
\begin{definition}
For $x_0\leq0$ and $y_0,\,y^0\geq y_l$ define
$$
\mathcal Z(x_0,y_0,y^0;y_l)=
\left\{
\begin{array}{ccc}
\mathcal Z_-(x_0,y_0,y^0;y_l),&{\rm for}& y^0\geq y_0\geq y_l\\
\mathcal Z_+^{-1}(x_0,y^0,y_0;y_l),&{\rm for}&y_0\geq y^0\geq y_l
\end{array}
\right.,
$$
where the notation $\mathcal Z_+^{-1}(x_0,y^0,y_0;y_l)$ denotes a section
of the fibre $\mathcal Z_+^{-1}(x_0)$ by the hyperplane parallel to the
$x$-axis and passing through the point with the coordinates
$\{0,y^0,y_0,y_l\}$.
\end{definition}
\begin{proposition}
 \label{Prop:z-4}
The function $\mathcal Z(x_0,y_0,y^0;y_l)$ is a smooth function of
all its variables. It has the following properties:
$$
{\mathcal X}(x_0,y_l)<\mathcal Z(x_0,y_0,y^0;y_l)<x_0,
$$
$$
\underset{y_0,y_1\to+\infty}\lim\mathcal Z(x_0,y_0,y^0;y_l)=
{\mathcal X}(x_0,y_l),\quad
\underset{y_l\to\pm\infty}\lim\mathcal Z(x_0,y_0,y^0;y_l)=x_0.
$$
\begin{eqnarray*}
y^1>y^0&\Rightarrow&\mathcal Z(x_0,y_0,y^1;y_l)<
\mathcal Z(x_0,y_0,y^0;y_l),\\
y_{01}>y_0&\Rightarrow&\mathcal Z(x_0,y_{01},y^0;y_l)<
\mathcal Z(x_0,y_0,y^0;y_l),\\
0\geq x_2>x_1&\Rightarrow&0
<\mathcal Z(x_2,y_0,y^0;y_l)-\mathcal Z(x_1,y_0,y^0;y_l)<x_2-x_1.
\end{eqnarray*}
\end{proposition}
\begin{proof}
The properties are the same as the ones for the function
$\mathcal Z_-(x_0,y_0, y^0;y_l)$, without the restriction $y^0\geq y_0$.
Continuity as well as differentiability, with respect to $x_0$ and $y_l$,
at the ``matching hyperplane'' $y_0=y^0$, follows from the identity,
$\mathcal Z_-(x_0,y_0, y_0;y_l)=\mathcal Z_+^{-1}(x_0,y_0,y_0;y_l)$ and
corresponding properties of the functions $\mathcal Z_\pm$.
Differentiability on $y_0$ and $y^0$ at the matching hyperplane is
deduced from the observation that both functions are defined as the
abscissas of the crossing point of the straight line $y=y^0$ by a unique
family of solutions $y(x)$ that depend smoothly on their initial value $y_0$.

Monotonicity properties of the function
$\mathcal Z_+^{-1}(x_0,y^0,y_0;y_l)$
with respect to the $y$-variables do not follow from the corresponding
properties of $\mathcal Z_+(x_0,y^0,y_0;y_l)$ mentioned in
Proposition~\ref{Prop:z+}, however, they can be established in analogous
way on the basis of the results of \S\S~\ref{sect:3} and \ref{sect:5}.
\end{proof}
\begin{definition}
 \label{Def:z-main}
For $x_0\leq0$ and $y_0,\,y^0\in\mathbb R$, denote by $Y$ the set of
solutions $y(x)$ of Equation~(\ref{p1}) with the initial value $y(x_0)=y_0$.
Define the function
$$
Z(x_0,y_0,y^0)=\underset{Y}\inf\left\{z: y(z)=y^0\right\}.
$$
\end{definition}
\begin{proposition}
 \label{Prop:z-unique}
For any given $x_0\leq0$ and $y_0,\,y^0\in\mathbb R$,
there is a unique solution $y(x)$ of Equation~{\rm(\ref{p1})} that
corresponds to the function $Z(x_0,y_0,y^0)$, i.e., $y(x_0)=y_0$ and
$y\left(Z(x_0,y_0,y^0)\right)=y^0$. This solution achieves its (global)
minimum in the segment $\left[Z(x_0,y_0,y^0),x_0\right]$.
\end{proposition}
\begin{proof}
The infimum in Definition~\ref{Def:z-main} is achieved on the solutions
with the minimum to the left of $x_0$. This follows from the monotonic
growth of $\mathcal Z(x_0,y_0,y^0;y_l)$ with respect to the
$x$-variable (Proposition~\ref{Prop:z-4}).
Using the continuity of the function $x_0-\mathcal Z(x_0,y_0,y^0;y_l)$
and the fact that it vanishes as $y_l\to\pm\infty$ we arrive to the
existence of $\hat y_l$ such that
\begin{equation}
 \label{eq:Z_hat_yl}
Z(x_0,y_0,y^0)=\mathcal Z(x_0,y_0,y^0;\hat y_l).
\end{equation}
Equation~(\ref{eq:Z_hat_yl}) verifies existence of the solution $y(x)$.

{}From Theorem~\ref{Thm:dual} it follows that there exist at most two
solutions that correspond to $Z(x_0,y_0,y^0)$. Suppose that there
are two maximal solutions, $y_1(x)$ and $y_2(x)$, passing through
the points $\{x_0,y_0\}$ and $\left\{Z(x_0,y_0,y^0),y^0\right\}$.
Theorem~\ref{Thm:intersections} implies that these solutions
do not have any other points of intersection, therefore $y_1(x)>y_2(x)$
for $x\in\left(Z(x_0,y_0,y^0),x_0\right)$.
Consider a solution $y_3(x)$
passing through the point $\{x_0,y_0\}$ and having the slope between
the slopes of the maximal solutions. The third solution cannot intersect
$y_1(x)$ on the segment $x\in\left[Z(x_0,y_0,y^0),x_0\right]$ as this
contradicts Lemma~\ref{Lem:MN}. Thus, $y_3(x)$ intersects $y_2(x)$
at some $x_1\in\left(Z(x_0,y_0,y^0),x_0\right)$. Solutions $y_3(x)$ and
$y_2(x)$ do not have any other points of intersection. This means that
$y_3(x)<y_2(x)$ for $x\in\left[Z(x_0,y_0,y^0),x_1\right]$. Therefore,
$y_3\left(Z(x_0,y_0,y^0)\right)< y_2\left(Z(x_0,y_0,y^0)\right)=y^0$
and solution $y_3(x)$ should cross the level $y^0$ to the left of
$Z(x_0,y_0,y^0)$. That contradicts the definition of $Z(x_0,y_0,y^0)$
and establishes uniqueness of $y(x)$.

Suppose that the solution $y(x)$ has no minimum in the segment
$\left[Z(x_0,y_0,y^0),x_0\right]$. Denote by $x_{min}$ the global minimum
of $y(x)$ and suppose first that $y'(x)>0$ for
$x\in\left[Z(x_0,y_0,y^0),x_0\right]$, clearly $x_{min}<Z(x_0,y_0,y^0)$.
Consider the segment $[x_{min},x_0]\supset\left[Z(x_0,y_0,y^0),x_0\right]$.
If $\epsilon$ is small enough, then the graphs of all solutions
$y_{\mu}(x):$ $y_{\mu}(x_0)=y(x_0)$, $y'_{\mu}(x_0)\equiv\mu>y'(x_0)$,
and $\mu-\epsilon<y'(x_0)$ would cross the straight line
$x=Z(x_0,y_0,y^0)$ below $y^0$. By virtue of Lemma~\ref{Lem:3.1},
$y_\mu(x)$ would not cross the graph of $y(x)$ at any $x<x_0$ until they
get to their minimum point, which, by continuity, for small $\epsilon$ is
close enough to $x_{min}$. Therefore, the graphs of all these solutions
cross the line $y=y^0$ to the left of the point with the abscissa
$Z(x_0,y_0,y^0)$. That contradicts the definition of the function
$Z(x_0,y_0,y^0)$. Suppose now that $y'(x)<0$ for
$x\in\left[Z(x_0,y_0,y^0),x_0\right]$. Consider the solution $y_0(x)$
with the following initial data $y_0(x_0)=y(x_0)=y_0$ and $y'_0(x_0)=0$.
As follows from Lemma~\ref{Lem:4.1} $y_0(x)<y(x)$ for $x<x_0$. Thus its
graph intersects with the line $y=y^0$ to the left of the point with the
coordinates $\left(Z(x_0,y_0,y^0),y^0\right)$. That contradicts the
definition of $Z(x_0,y_0,y^0)$.
\end{proof}
\begin{definition}
We call the solution from Proposition~\ref{Prop:z-unique}
the {\em maximal solution corresponding to the segment}
$\left[Z(x_0,y_0,y^0),x_0\right]$.
\end{definition}
\begin{proposition}
 \label{Prop:z-main}
The function $Z(x_0,y_0,y^0)$ is a smooth function of all its variables
with the following properties:
\begin{eqnarray*}
X(x_0)<Z(x_0,y_0,y^0)<x_0,&&
\underset{y_0,y^0\to+\infty}\lim Z(x_0,y_0,y^0)=X(x_0),\\
\underset{y_0\to-\infty}\lim Z(x_0,y_0,y^0)=x_0,&&
\underset{y^0\to-\infty}\lim Z(x_0,y_0,y^0)=x_0,
\end{eqnarray*}
\begin{gather*}
\hat y^0>y^0\Rightarrow Z(x_0,y_0,\hat y^0)<Z(x_0,y_0,y^0),\;\;
\hat y_0>y_0\Rightarrow Z(x_0,\hat y_0,y^0)<Z(x_0,y_0,y^0),\\
0\geq\hat x_0>x_0\Rightarrow0<Z(\hat x_0,y_0,y^0)-Z(x_0,y_0,y^0)<\hat x_0-x_0.
\end{gather*}
\end{proposition}
\begin{proof}
The proof of smoothness is analogous to the one for $X(x_0)$ given in
Corollary~\ref{Cor:smoothness}; it is based on the uniqueness of
representation~(\ref{eq:Z_hat_yl}) proved in Proposition~\ref{Prop:z-unique}
and smoothness of the corresponding level function
$\mathcal Z(x_0,y_0,y^0;y_l)$, see Proposition~\ref{Prop:z-4}. The other
properties are proved analogously as the ones in Proposition~\ref{Prop:3.7}
with the help of representation (\ref{eq:Z_hat_yl}).
\end{proof}
\begin{remark}
In Appendix~\ref{App:maximal_solutions} we give estimates for
the function $Z(x_0,y_0,y^0)$ similar to those obtained for the function
$X(x_0)$ in \S~\ref{sect:2}.
\end{remark}
\begin{proposition}
 \label{Prop:monotonic_up}
Let $x_0\leq0$, $y_0,y^0\in\mathbb R$. Consider solutions $y(x)$ of
Equation~{\rm\ref{p1}} with the initial data $y(x_0)=y_0$ and
$y'(x_0)=y_1\leq y_{1m}$, where $y_{1m}$ is the initial slope at $x_0$
of the maximal solution $y_m(x)$ corresponding to the segment
$\left[Z(x_0,y_0,y^0),x_0\right]$.
Graphs of any two such solutions, $y_1(x)$ and $y_2(x)$, have no
(finite/infinite) points of intersection with abscissa
$x\in\left[Z(x_0,y_0,y^0),x_0\right]$ other than $\{x_0,y_0\}$. In
particular, for any fixed $x_1\in\left(Z(x_0,y_0,y^0),x_0\right)$,
$y(x_1)\equiv y(x_1,y_1)$ is a smooth monotonically decreasing function
of $y_1\in[c(x_1),y_{1m}]$ with the image $[y_m(x_1),+\infty]$, where
the function $c(x_1)$ takes values in $(-\infty,y_{1m})$.
\end{proposition}
\begin{proof}
The absence of intersections and monotonicity is a direct consequence of
Projective Lemma~\ref{Lem:MN}. Denote by $Y_1$ the set of initial slopes
$\{y_1\}$ such that corresponding solutions $y(x,y_1)$ cross the straight
line $x=x_1$ with finite values $y(x_1,y_1)$. It is clear that $Y_1$ is
nonempty ($y_{1m}\in Y_1$) and open (by smoothness of $y(x,y_1)$ on $y_1$)
set. It is also bounded. To prove it, consider inequality $y''>6y^2$ for
$y=y(x,y_1)$, $x<0$, and $y_1<0$ and such that $N=y_1^2/2-2y_0^3>0$.
Multiply it on $y'<0$, integrate from $x<x_0$ to $x_0$, and take a square
root to get $\sqrt{4y^3+2N}<-y'$.
After standard rearrangement, integrate the latter inequality from the left
bound $a$ of the interval of existence of $y(x)$ to $x_0$. Estimating the
remaining definite integral, one arrives at the following inequality
$0<x_0-a<\frac2{(2N)^{1/6}}$. This inequality shows that for all
sufficiently large negative values of $y_1$, the solution $y(x,y_1)$ does
not cross the line $x=x_1$. Thus, $c(x_1):=\inf Y_1>-\infty$.
The value $y(x_1,c(x_1))=+\infty$, as otherwise the monotonicity
and smoothness of $y(x_1,y_1)$ with respect to $y_1$ imply a contradiction
with the definition of $c(x_1)$.
\end{proof}
\begin{proposition}
 \label{Prop:monotonic_down}
Let $x_0\leq0$, $y_0,y^0\in\mathbb R$. Consider solutions $y(x)$ of
Equation~{\rm(\ref{p1})} with the initial value $y(x_0)=y_0$ and minimum
value $y_l\leq y_{lm}$, where $y_{lm}$ is the minimum value of the maximal
solution $y_m(x)$ corresponding to the segment
$\left[Z(x_0,y_0,y^0),x_0\right]$. Graphs of any two such solutions,
$y_1(x)$ and $y_2(x)$, have one more (finite) point of intersection in
the segment $x\in\left[Z(x_0,y_0,y^0),x_0\right]$ different from $x_0$.
The function $\mathcal Z(x_0,y_0,y^0;y_l)$ is a strictly monotonically
decreasing function from $x_0$ to $Z(x_0,y_0,y^0)$ as $y_l$ increases
from $-\infty$ to $y_{lm}$.
\end{proposition}
\begin{proof}
Without lost of generality assume that the minimum values of the
solutions are ordered as follows, $y_{lm}>y_{1l}>y_{2l}$.
Denote by $P$ the point of intersection of the graphs of $y_m(x)$ and
$y_1(x)$ with abscissa in $\left(Z(x_0,y_0,y^0),x_0\right)$.
This point exists due to the uniqueness of the maximal solution $y_m(x)$.
If we suppose that the graph of $y_2(x)$ does not intersect the graph of
$y_1(x)$ to the right of $P$, then we arrive at a contradiction with
Lemma~\ref{Lem:MN}. Denote by $Q$ the point of intersection of the graphs
$y_1(x)$ and $y_2(x)$. To prove strict monotonicity of
$\mathcal Z(x_0,y_0,y^0;y_l)$ for $y_l\leq y_{lm}$, let us notice that,
according to the definition, it is the abscissa of the point of intersection
of the graphs of solutions with the level $y=y^0$ located to the left of
their minima. The graphs of $y_1(x)$ and $y_2(x)$ have no points of
intersection other than $\{x_0,y(x_0)\}$ and $Q$. This means that to the
left of $Q$ $y_2(x)>y_1(x)$; therefore the graph of $y_2(x)$ intersects the
level $y=y^0$ at a point with larger abscissa than it does the graph of
$y_1(x)$.
\end{proof}
\begin{theorem}
Let $x_1<x_0\leq0$. The Dirichlet boundary value problem
$y(x_0)=y_0$ and $y(x_1)=y^0$ has\\
\indent
$1$. No solution if $x_1<Z(x_0,y_0,y^0)$;\\
\indent
$2$. A unique solution if $ x_1=Z(x_0,y_0,y^0)$;\\
\indent
$3$. Two solutions if $Z(x_0,y_0,y^0)<x_1<x_0$.
\end{theorem}
\begin{proof}
The first statement follows from Definition~\ref{Def:z-main}
of the function $Z(x_0,y_0,y^0)$ and its boundedness
(Proposition~\ref{Prop:z-main}). The second statement is proved in
Proposition~\ref{Prop:z-unique}.

Consider the third statement. Suppose first that the straight lines $y=y^0$
and $x=x_1$ cross above the maximal solution $y_m(x)$. As follows from
Proposition~\ref{Prop:monotonic_up} there is exactly one solution $y_1(x)$
of the Dirichlet boundary problem with the initial slope
$y'_1(x_0)<y'_m(x_0)$. It satisfies the inequality $y_1(x)>y_m(x)$ for all
$x<x_0$ in the interval of existence of $y_1(x)$.\\
On the other hand, from Proposition~\ref{Prop:monotonic_down} it follows
that there is exactly one solution $y_2(x)$ of the Dirichlet problem in the
class of solutions of Equation~(\ref{p1}) with the (global) minimal value
$y_{2l}<y_{lm}$, where $y_{lm}$ is the minimum value of $y_m(x)$.
It is clear that $y_2(x)$ achieves its minimum value to right of $x_1$,
as the solutions with minimum values less than $y_{lm}$ and (global) minima
located to the left of $x_1$ cross the line $x=x_1$ below $y_m(x_1)<y^0$.
The solution $y_2(x)$ crosses the maximal solution at some point $x_2$:
$x_1<x_2<x_0$.

Consider the case when the lines $y=y^0$ and $x=x_1$ cross below the
maximal solution. According to Proposition~\ref{Prop:monotonic_down}
we still have exactly one solution $y_2(x)$ to the Dirichlet problem
in the class of the functions with the (global) minimum values
$y_{2l}<y_{lm}$ whose minima are located to the {\it right} of $x_1$.
We claim that we also have the second solution, $y_1(x)$, which is the
unique solution of the Dirichlet problem in the class of functions with
the (global) minimum value $y_{1l}<y_{lm}$ whose minima are located to
the {\it left} of $x_1$.
Existence of such a solution is established by considering a family of
solutions of Equation~(\ref{p1}), $y_\mu(x)$, defined by the initial
data: $y_\mu(x_0)=y_0$ and $y'_\mu(x_0)\equiv\mu\geq y'_m(x_0)$.
As follows from Corollary~\ref{Cor:3.1}, the minimum values of $y_\mu(x)$,
satisfy the inequality $y_{l\mu}<y_{lm}\leq y^0$. Therefore, the graphs of
$y_\mu(x)$ cross the line $y=y^0$ to the right of their minima.
As $\mu\to+\infty$ the crossing point approaches $x_0$. This fact follows
from Corollary~\ref{Cor:3.1} and Proposition~\ref{Prop:3.6}. Actually,
the function $y_{l\mu}(\mu)$ is the inverse to $\mu=f(x_0;y_0,y_{l\mu})$
for fixed $x_0$ and $y_0$, thus Corollary~\ref{Cor:3.1} implies that
$y_{l\mu}(\mu)\searrow-\infty$ as $\mu\nearrow+\infty$.
For $\mu\to y'_m(x_0)$, the crossing point
approaches the crossing point of $y_m(x)$ with the line $y=y^0$,
which is located to the left of the point $(x_1,y^0)$. Thus,
existence of $y_1(x)$ follows by standard continuity arguments.
Uniqueness is the consequence of the fact that the graphs of any two
(different) solutions $y_{\mu_1}(x)$ and $y_{\mu_2}(x)$ in the right
half-plane of the largest of their minima have only one point of
intersection, $(x_0,y_0)$, (cf. Lemma~\ref{Lem:3.1}).
In this case $y_k(x)<y_m(x)$, $k=1,2$, for $x\in[x_1,x_0)$.

Finally we have to consider the case when the point $(x_1,y^0)$ lies on
the graph of the maximal solution. This is a limit of the first case
considered in this proof, namely, the ``upper'' solution $y_1(x)$
coincides with $y_m(x)$.
It can also be considered as a limit of the second case considered in
the proof, when again $y_1(x)$ coincides with $y_m(x)$.
Either case leads to existence of two solutions: the first
coincides with $y_m(x)$ and the second is $y_2(x)$ which is the
same in both cases.
Note, that in this case the requirement  $x_1>Z(x_0,y_0,y^0)$ implies
that $y_{lm}<y^0$.
\end{proof}
\section{Numerics}
 \label{sect:numerics}
All calculations for this section, except those presented in
Table~\ref{tab:3}, are done with the help of \textsc{Maple 6}. The results of
Table~\ref{tab:3} are calculated with \textsc{Maple 8}.

We begin with the illustration of Corollary \ref{Cor:3.1} and
Remark~\ref{Rem:symmetry}, namely,
in Table \ref{tab:1} we present some numerical values of the functions
$f(x_0;y_0,y_l)$ and $\Delta(x_0;y_0,y_l)$.
Corresponding calculations are performed twice with the \textsc{Maple}
parameter Digits set to 16 and 18. These settings allow one to calculate
values of the functions $f$ and $\Delta$ up to 6 decimal digits.
Table \ref{tab:1}, together with Remark \ref{Rem:symmetry} allow us to
conjecture that the upper bound of function $\Delta$ is very close to zero,
most probably, it is less than $10^{-1}$. This suggests the following
interesting property, an {\it approximate symmetry}, of solutions of
Equation~(\ref{p1}) defined in Corollary~\ref{Cor:3.1}.
\begin{conjecture}
Consider two solutions: the
first one defined via the initial value $y_0\in\mathbb{R}$ and minimum value
$y_l\leq y_0$, the second one via $-y_l$ and $-y_0$, respectively.
The integer part of the difference between their initial slopes is zero,
i.e., $[\Delta]=0$.
\end{conjecture}
Numerical studies also confirms an analogous conjecture for the
function
$$
\Delta^+(x_0;y_0,y_l)=f^+(x_0;y_0,y_l)-f^+(x_0;-y_l,-y_0).
$$
\begin{table}[ht]
 \caption{}
  \label{tab:1}
\begin{tabular}{|c|c|c||c|c|}
\hline
$x_0$&$y_0$&$y_l$&$f(x_0;y_0,y_l)$&$\Delta(x_0;y_0,y_l)$\\
\hline
\phantom{-}0.0&\phantom{26}0.2&\phantom{-30}0.1
&\phantom{2852}0.270736\ldots&0.004050\ldots\\
\hline
\phantom{-}0.0&\phantom{26}2.0&\phantom{-30}1.0
&\phantom{2852}5.319113\ldots
&0.006387\ldots\\
\hline
-1.3&\phantom{2}46.0&\phantom{-30}2.0
&\phantom{82}624.047258\ldots&0.004173\ldots\\
\hline
-2.0&600.0&-300.0&31176.973126\ldots&0.000401\ldots\\
\hline
-3.0&\phantom{2}10.0&\phantom{20}-5.0&\phantom{842}67.798893\ldots&
0.022622\ldots\\
\hline
-5.0&262.0&\phantom{-30}1.0&\phantom{2}8481.836260\ldots
&0.001160\ldots\\
\hline
\end{tabular}
\end{table}

As proved in Section~\ref{sect:2}, $X(0)$ and $X_{min}(0)$
are finite, however, the results obtained in that section allow us to
find only rough estimates of their numerical values
(see Remark~\ref{Rem:X_min_0}).
In Table~\ref{tab:2}, we present numerical calculations of some important
data characterizing three distinguished solutions of Equation (\ref{p1})
that have a pole at $x_0=0$. In the second and last lines of this table
we give the values of $v(\hat x)=c/7$, where
$c$ is the parameter in the corresponding Laurent expansion (\ref{laurent})
of the solutions at $\hat x=0$ and $x_p$, the right and, respectively, left
bounds of their intervals of existence. The solutions have two zeroes
$z_1$ and $z_2$. In the second column we present numerical data
characterizing $y_{max}(0;x)$, the solution with the maximum interval of
existence, i.e., $x_p=X(0)$. In the third column, the data are given for
$y_{min}(0;x)$, the solution with $x_{min}=X_{min}(0)$. Finally, in the
last column of Table~\ref{tab:2}, we list the data for $y_{s-sym}(x)$,
the symmetric solution that has a pole at $0$. Equation~(\ref{p1}) has
actually two symmetric solutions which were distinguished long ago by
Boutroux \cite{boutroux:I}. They are uniquely characterized by the
following symmetry condition, $y(x)=\varepsilon^2y(\varepsilon x)$, where
$\varepsilon^5=1$, and behaviour at $x=0$: one of them, $y_{r-sym}(x)$ is
regular at $0$ and the other, $y_{s-sym}(x)$, has a pole.
\begin{table}[ht]
 \caption{}
  \label{tab:2}
\begin{tabular}{|c||l|l|l|}
\hline
&\phantom{aaa}$y_{max}(0;x)$&\phantom{aaa}$y_{min}(0;x)$
&\phantom{aaai}$y_{s-sym}(x)$\\
\hline
$v(0)$&\phantom{-}0.110489160\ldots
&\phantom{-}0.125565964\ldots&\phantom{-}0.0\\
\hline
$z_1$&-1.528989716\ldots&-1.537495773\ldots&-1.476591053\ldots\\
\hline
$y'(z_1)$&\phantom{-}1.113243043\ldots&\phantom{-}1.082837664\ldots
&\phantom{-}1.313083166\ldots\\
\hline
$x_{min}$&-2.055505831\ldots&-2.055703500\ldots&-2.044984309\ldots\\
\hline
$y_{min}$&-0.322633511\ldots&-0.307468294\ldots&-0.423460899\ldots\\
\hline
$z_2$&-2.546577118\ldots&-2.539078168\ldots&-2.575998303\ldots\\
\hline
$y'(z_2)$&-1.292743827\ldots&-1.256197484\ldots&-1.527257430\ldots\\
\hline
$x_p$&-3.915285797\ldots&-3.914972029\ldots&-3.902470099\ldots\\
\hline
$v(x_p)$&-0.916786830\ldots&-0.892041655\ldots&-1.091093248\ldots\\
\hline
\end{tabular}
\end{table}
The notation in the first column of Table~\ref{tab:3} has the same meaning
as in the corresponding column of Table~\ref{tab:2}; the only difference is
that in Table~\ref{tab:2} the right pole of the solutions, $x=0$, is known,
while in Table~\ref{tab:3} we add an additional line, labeled $x_{p1}$,
with the information about this pole. In the second column the
data for the symmetric solution $y_{r-sym}(x)$ mentioned above are given.
In the third and fourth columns we list the data for solutions $y_-(x_0;x)$
for $x_0=0$ and $x_0=-1$, respectively. These solutions correspond to
$X_-(x_0)$, i.e., they have a minimum at $x_0$ and the left pole
$x_{p2}=X_-(x_0)$. A comment, probably should be given to the line labeled
$x_{min}$. The zero is the global minimum of the restrictions of
$y_{r-sym}(x)$ and  $y_-(0;x)$ on the non-positive semi-axis. If we consider
this solutions, as we should, on the maximum interval of existence, then
zero is the (local) minimum of $y_-(0;x)$ and an inflection point for
$y_{r-sym}(x)$.
\begin{table}[ht]
 \caption{}
  \label{tab:3}
\begin{tabular}{|c||l|l|l|}
\hline
&\phantom{aa}$y_{r-sym}(x)$&\phantom{aa}$y_-(0;x)$&
\phantom{aaai}$y_-(-1;x)$\\
\hline
$x_{p1}$&\phantom{-}regular for $x\geq0$&\phantom{-}regular for $x\geq0$
&\phantom{-}1.848036525\ldots\\
\hline
$v(x_{p1})$&\phantom{-}N/A&N/A&\phantom{-}0.161869969\ldots\\
\hline
$z_1$&\phantom{-}no zeroes for $x>0$&\phantom{-}no zeroes for $x>0$
&-0.303329968\ldots\\
\hline
$y'(z_1)$&\phantom{-}N/A&N/A&\phantom{-}0.583881922\ldots\\
\hline
$x_{min}$&\phantom{-}0\;\;{\rm for} $x\leq0$
&\phantom{-}0\;\;{\rm for} $x\leq0$&-1\\
\hline
$y_{min}$&\phantom{-}0&-0.124293080\ldots&-0.249902470\ldots\\
\hline
$z_2$&\phantom{-}no zeroes for $x<0$&-0.838060764\ldots&-1.582985950\ldots\\
\hline
$y'(z_2)$&\phantom{-}N/A&-0.398459416\ldots&-0.870257802\ldots\\
\hline
$x_{p2}$&-2.615571209\ldots&-2.677058361\ldots&-3.121759948\ldots\\
\hline
$v(x_{p2})$&-0.371644061\ldots&-0.438744582\ldots&-0.639793560\ldots\\
\hline
\end{tabular}
\end{table}
In particular, from Tables \ref{tab:2} and \ref{tab:3} it follows that:
$$
X_{min}(0)\!=\!-2.055703500\ldots,\,
X_-(0)\!=\!-2.677058361\ldots,\,
X(0)\!=\!-3.915285797\ldots.
$$
It is interesting to compare how these numerical values fit the bounds given
in Section~\ref{sect:2}; namely, in Remark~\ref{Rem:X_min_0} and
Inequalities~(\ref{ineq:eta_est}), (\ref{ineq:X_-eta_est}):
$$
-3.239042305\ldots=-5\left(\frac C4\right)^{4/5}\!<X_-(0)<
X_{min}(0)<-C^{4/5}=-1.963788033\ldots,
$$
and second Inequalities~(\ref{ineq:1thm1}),  (\ref{ineq:2thm1}):
$$
-3.997162138\ldots\!=X_{min}(0)-\frac C{|X_{min}(0)|^{1/4}}<X(0)<\!
-(2C)^{4/5}\!=\!-3.419153556\ldots.
$$
Another observation is that $y_{s-sym}(x)$ appear to have the least
interval of existence comparing to the other solutions in
Table~\ref{tab:2}, however, it has considerably larger spacing between the
zeroes. The solutions presented in Table~\ref{tab:2} have the following
spacing between their zeroes: $1.017587402\ldots$, $1.001582395\ldots$, and
$1.099407250\ldots$, respectively. For a comparison we calculated the
the maximum spacing between the zeroes for solutions  with the pole at
$x_0=0$ and the interval of existence to the left of it, i.e., the function
$\delta_-(0,+\infty)$ in the notation of Appendix~\ref{App:zero-spacing};
$\delta_-(0,+\infty)=1.1808499889180\ldots$. The corresponding solution has
the following data: $v(0)=-0.518045\ldots$, $z_1=-1.3362856\ldots$,
$x_{min}=-1.9417146\ldots$, $y_{min}=-0.741427\ldots$,$z_2=-2.5171356\ldots$
$x_p=-3.7427412\ldots$, $v(x_p)=-1.798000\ldots$.
Numerical calculation shows that one gets more than a twofold increase in
accuracy for the spacing function compared to the accuracy of calculation of
other data, in particular, the zeroes themselves.

It is also important to mention that our numerical studies support the
conjecture of Remark~\ref{Rem:conjecture} that the solution with
$x_{min}=X_{min}(0)$ is unique.

The notation used in the first and fourth columns of Table \ref{tab:4}
has the same meaning as the one in the first column of Table~\ref{tab:2}.
In the second and third columns of Table \ref{tab:4} we present numerical
data for solutions $y_{max}(-1;x)$ and $y_{min}(-1;x)$ which are analogous
to $y_{max}(0;x)$ and $y_{min}(0;x)$, but with the right pole at $x_0=-1$.
\begin{table}[ht]
 \caption{}
  \label{tab:4}
\begin{tabular}{|c||l|l||c||l|}
\hline
&\phantom{aa}$y_{max}(-1;x)$&\phantom{aa}$y_{min}(-1;x)$&
&\phantom{aaai}$\hat y_{min}(x)$\\
\hline
$v(-1)$&-0.064360748\ldots&-0.045841066\ldots&$v(X(0))$&-0.943704177\ldots\\
\hline
$z_1$&-2.374548143\ldots&-2.379585633\ldots&$z_2$&-2.554014196\ldots\\
\hline
$y'(z_1)$&\phantom{-}1.468791556\ldots&\phantom{-}1.441765379\ldots
&$y'(z_2)$&-1.331432779\ldots\\
\hline
$x_{min}$&-2.853609725\ldots&-2.853690013\ldots&$x_{min}$
&-2.055297172\ldots\\
\hline
$y_{min}$&-0.381066130\ldots&-0.369105745\ldots&$y_{min}$
&-0.338834697\ldots\\
\hline
$z_2$&-3.311108506\ldots&-3.306469422\ldots&$z_1$&-1.520595470\ldots\\
\hline
$y'(z_2)$&-1.620851956\ldots&-1.590390526\ldots
&$y'(z_1)$&\phantom{-}1.145676637\ldots\\
\hline
$x_p$&-4.589970403\ldots&-4.589833499\ldots&$x_p$&-0.000570546\ldots\\
\hline
$v(x_p)$&-1.187366438\ldots&-1.161843481\ldots&$v(x_p)$
&\phantom{-}0.093928571\ldots\\
\hline
\end{tabular}
\end{table}

In particular, from Tables~\ref{tab:3} and \ref{tab:4} it follows that:
$
X(-1)=-4.589970403\ldots,
$
$$
X_{min}(-1)=-2.853690013\ldots,\qquad
X_-(-1)=-3.121759948\ldots.
$$
Comparison of these results with the estimations obtained in
Corollary~(\ref{Cor:X_min_asymptotics}) and Theorem~(\ref{Thm:1}) gives:
$$
-3.324707204\ldots\!=-1-C\!<\!X_-(-1)\!<\!X_{min}(-1)\!<\!\!
-1-\frac C{\eta(-1)}\!=\!-2.797524387\ldots,
$$
\begin{eqnarray*}
-4.642303753\ldots&=&X_{min}(-1)-\frac C{|X_{min}(-1)|^{1/4}}<X(-1)\\
&<&-1-\frac{2^{4/5}C}{\eta(-2^{-4/5})}=-4.240073805\ldots,
\end{eqnarray*}
where, $\eta(-1)=1.293282706\ldots$ and $\eta(-2^{-4/5})=1.249215473\ldots$.
Experience with these calculations and numerical verification of
Conjecture~\ref{Con:4} (see below) support the following conjecture made in
Remark~\ref{Rem:conjecture}:
\begin{conjecture}
 \label{Con:2}
For any $x_0$ the solution with the pole at $x_0$ and minimum at
$X_{min}(x_0)$ is unique.
\end{conjecture}
In the last column of Table \ref{tab:4}
we present numerical data for solution $\hat y_{min}(x)$ which has a pole
at $X(0)$, interval of existence to the right of $X(0)$, and
the extremal property: $x_{min}=\Xi_{min}(X(0))=X_-(x_0)$. So, it serves
as an illustration of the content of Section \ref{sect:5} and supports,
together with the check of Conjecture~\ref{Con:4} below,
the conjecture made in Remark~\ref{Rem:conjecture2}:
\begin{conjecture}
 \label{Con:3}
For any $x_0$ the solution with the pole at $x_0$ and minimum at
$\Xi_{min}(x_0)$ is unique.
\end{conjecture}
Comparing the numerical data given in Tables~\ref{tab:2} and \ref{tab:4}
for $y_{min}(x_0;x)$ with that for $y_{max}(x_0;x)$, $x_0=0,-1$,
we see that they are very close, though distinct.

Moreover, in both cases, $x_0=0$ and $x_0=-1$,
\begin{equation}
 \label{ineq:maxmin}
y_{max}(x_0;x)<y_{min}(x_0;x),
\end{equation}
where $x$ belongs to the interval of existence of $y_{min}(x_0;x)$. We
actually, numerically checked Inequality~(\ref{ineq:maxmin}) for
$x_0\in[-3,0]$, so that it is natural to make the following conjecture.
\begin{conjecture}
 \label{Con:4}
Inequality {\rm(\ref{ineq:maxmin})} holds for any $x_0\leq0$.
\end{conjecture}
Moreover, it is easy to see that the leading term of asymptotics
as $x_0\to-\infty$ of both minimum values, $y_m:=y_{max}(x_{min};x_0)$ or
$y_{min}(x_{min};x_0)$ is given by the same formula,
$y_m=v^{max}_{min}\sqrt{|x_0|}+ o\left(\sqrt{|x_0|}\right)$, where
$v^{max}_{min}$ is defined in Remark~\ref{Rem:2.2} (for numerical estimates
it is better to use in this formula $|x_{min}|$ instead of $|x_0|$).

It would be interesting to find the next correction
term to this asymptotics and prove that both minima decrease monotonically
with the decrease of $x_0$, as numerical calculations also confirm. Another
interesting question related with these solutions is to prove or disprove
the following inequality:
\begin{conjecture}
 \label{Con:5}
For $x_1<x_2<0$,
$
0<X(x_2)-X(x_1)<X_{min}(x_2)-X_{min}(x_1).
$
\end{conjecture}
We have established an analogous inequality for the level functions
(see Proposition~\ref{Prop:4.2}, however, it does not imply
Conjecture~\ref{Con:5}.

The data presented in Tables \ref{tab:2}--\ref{tab:4} require precise
calculations, as some of the (different) quantities coincide up to four
digits after the decimal point. An important ingredient of these
calculations is a special representation, the so-called PR\PI-system, of
solutions of Equation (\ref{p1}) in the neighbourhood of a pole. Namely,
let $z(x)$ and $v(x)$ be a solution of the system,
\begin{equation}
\textnormal{PR\PI :}\phantom{aaaaaaaa}
\begin{array}{ll}
\displaystyle\frac{dz}{dx}&\displaystyle =
1-\frac x4z^4-\frac{z^5}4-\frac12v\,z^6\\
&\\
\displaystyle\frac{dv}{dx}&\displaystyle =
\frac{x^2}8\,z+\frac38x\,z^2+
\left(\frac14+x\,v\right)z^3+\frac54v\,z^4+\frac32v^2\,z^5
\end{array}
\end{equation}
with initial data $z(x_0)=0$, $v(x_0)=c/7$.
This is Painlev\'e's {\em regularised polar patch}
(PRP) system \cite{pp:bull}. Note that
\begin{equation}
 \label{eqn:zu}
y(x)=\frac1{z(x)^2},\quad y'(x)=-\,\frac2{z(x)^3}+\,\frac{x z(x)}2
+\frac{z(x)^2}2+v(x)z(x)^3,
\end{equation}
are compatible and the function $y(x)$ solves {\rm\PI}.

To define a solution of Equation (\ref{p1}) in a neighbourhood of a
right pole we choose some parameter $c$ in the initial data for the
PR\PI-system. In neighbourhoods of zeroes of solutions, $y(x)=0$, the
corresponding functions $z(x)$ and $v(x)\to\infty$. Therefore, to go
beyond the zeroes we choose a {\it connection} point in the segment
between the first zero and the right pole and continue our calculation
using our original Equation (\ref{p1}).
To reach the left pole of the solution we have to choose another
connection point after the second zero and continue our
calculations again with the help of the PR\PI-system. Note that in the
first segment $z(x)<0$ whilst in the second it is positive.
Numerical results depend, of course, on the choice of the connection points,
however, the better the precision the lesser the dependence. We believe that
the nine digits after the decimal point for the data presented in
Tables \ref{tab:2}\,-\,\ref{tab:4} coincide with those for the corresponding
true solutions. Our belief is based on calculations performed for
Tables~\ref{tab:2} and \ref{tab:4} with \textsc{Maple 6} by setting the
parameter Digits to $20$, $22$, and $24$, and for Table~\ref{tab:3} with
\textsc{Maple 8} by setting the absolute and relative errors to $10^{-16}$,
$10^{-18}$, and $10^{-20}$, and observation of the stability of the numbers
obtained under variation of the connection points.
\appendix
\section{Estimates for $Z(x_0,y_0,y^0)$}
 \label{App:maximal_solutions}
\begin{definition}
 \label{Def:Cv0}
For $v_0\in\mathbb R$ define the function
$$
C(v_0):=\underset{x>0,\,v_0\sqrt{2x}>-1}\sup C(v_0,x),\;\;
C(v_0,x):=(2x)^{\frac14}\!\int_0^{\sqrt{1+v_0\sqrt{2x}}}
\!\frac{dw}{\sqrt{w^4-3w^2+3+x}}.
$$
\end{definition}
\begin{proposition}
 \label{Prop:Cv0}
$C(v_0)$ is a strictly monotonic continuous function, $0<C(v_0)<C$,
where $C$ is defined in Equation~{\rm(\ref{eq:c})} of \S~{\rm\ref{sect:2}}.
Moreover,
$$
\underset{v_0\to-\infty}\lim C(v_0)=0,\qquad
\underset{v_0\to+\infty}\lim C(v_0)=C.
$$
\end{proposition}
\begin{proof}
The boundedness and strict monotonicity of $C(v_0)$ follow immediately
from Definition~\ref{Def:Cv0}.

If $v_0\leq0$, then the upper limit of the integral
$C(v_0,x)$ is $\leq1$ and $\sqrt{2x}\leq -1/v_0$. Thus, the integral is
uniformly bounded, with respect to $x$. If now $v_0\to-\infty$, then
the factor $(2x)^{1/4}\to0$ and we arrive at the first limit. Note that
$C(v_0,0)=0$, therefore, since $C(v_0,x)$ is a continuous function,
for $v_0\leq0$ the supremum in Definition~\ref{Def:Cv0} is achieved at
some finite value of $x>0$ which is a local maximum of $C(v_0,x)$.

For $v_0>0$, $C(v_0,+\infty):=\underset{x\to+\infty}{\lim}C(v_0,x)=\sqrt[4]2
\int_0^{\sqrt{v_0}\sqrt[4]2}\frac{du}{\sqrt{u^4+1}}$. Now, we make a
change of the variable of integration in $C(v_0,x)$, $w=x^{1/4}u$ and
consider the difference, $C(v_0,x)-C(v_0,+\infty)$.
For all $x$ this difference is strictly larger than the following sum of
two integrals,
$\left(\int_0^{1/\sqrt[4]x}+
\int_{1/\sqrt[4]x}^{\sqrt{v_0}\sqrt[4]2}\right)Ddu$,
where $D$ is the difference of the corresponding integrands.
The first integral is negative, whilst the second one is positive. For
large $x$ the first integral can be estimated as
$\mathcal O\left(1/x^{3/4}\right)$ and the second behaves as
$\mathcal O\left(1/x^{1/2}\right)$. Thus the supremum of $C(v_0,x)$
in Definition~\ref{Def:Cv0} is achieved at some finite local maximum of
$C(v_0,x)$. For each $v_0$ we denote any such maxima as $x_{max}(v_0)$.
It is easy to observe that $\frac\partial{\partial x}C(v_0,x)>0$ for
$v_0\geq0$ and $x\leq3/4$, thus for $v_0\geq0$, $x_{max}(v_0)>3/4$.
The above proof also works for $v_0=+\infty$. In the last case we
denote $x_{max}(+\infty)=x_{max}$. (According to the numerical studies it
is unique and its numerical value is given in Remark~\ref{Rem:2.2}, however,
we do not use this uniqueness in the rest of the proof.)

The numbers $x_{max}(v_0)$ are uniformly bounded, i.e., there exists
$\hat x_{max}$ that for all $v_0$: $x_{max}(v_0)\leq\hat x_{max}$.
Actually suppose that for $v_n\to +\infty$ the corresponding sequence
$x_{max}(v_n)\to+\infty$. Then according to the above considerations
we have that $C(+\infty,x_{max})>C(+\infty,+\infty)$. It follows from this
that for all sufficiently large $n$: $C(v_n,x_{max})>C(v_n,x_{max}(v_n))$,
which contradicts the definition of $x_{max}(v_n)$. If $v_n\to v_0>0$,
the proof is analogous; for  $v_n\to v_0<0$ $\sqrt{2x_{max}(v_0)}<-1/v_0$.
If $v_n\to0$ and $x_{max}(v_n)\to+\infty$, then
$\lim C(v_n,x_{max}(v_n))=0<C(0,x_{max}(0))$. This again leads to a
contradiction with the definition of $x_{max}(v_n)$.

The proof of continuity of $C(v_0)$ is analogous to the above proof of the
uniform boundedness of $x_{max}(v_0)$. In fact, suppose $v_{n_k}\to v_0$ and
$C(v_n,x_{max}(v_n))\not\to C(v_0,x_{max}(v_0))$. Then there is a subsequence
$v_{n_k}$ such that $x_{max}(v_{n_k})\to\hat x$ and
$C(v_{n_k},x_{max}(v_{n_k}))\to C(v_0,\hat x)<C(v_0,x_{max}(v_0))$. The last
inequality means that for all sufficiently large $n_k$:
$C(v_{n_k},x_{max}(v_0))>C(v_{n_k},x_{max}(v_{n_k}))$, which is a
contradiction.

Now the second limit can be proved as follows,
$$
C=C(+\infty,x_{max})=\underset{v_0\to+\infty}{\lim}C(v_0,x_{max})\leq
\underset{v_0\to+\infty}{\underline{\lim}}C(v_0,x_{max}(v_0)).
$$
On the other hand, for all $v_0\in\mathbb R$ and $x\geq0$: $C(v_0,x)<C$,
that implies
$
\underset{v_0\to+\infty}{\overline{\lim}}C(v_0,x_{max}(v_0))\leq C.
$
\end{proof}
\begin{remark}
Most probably the functions $x_{max}(v_0)$ from the above proof are
single-valued. This leads to the smoothness of $C(v_0)$ in the spirit
of Corollary~\ref{Cor:smoothness}.
\end{remark}
\begin{definition}
 \label{Def:X_{min}(x_0,y_0)}
For $x_0\leq0$ and $y_0\in\mathbb R$
$$
X_{min}(x_0,y_0):=\underset{y_l}{\inf}\{x_{min}: \mathrm{the\;minimum\;of\;
solution}\;y(x;x_0,y_0,y_l)\}.
$$
\end{definition}
\begin{proposition}
 \label{Prop:X_{min}(x_0,y_0)}
For $x_0<0$ and $y_0\in\mathbb R$, define $\eta_{y_0}(x_0)$ as the unique
positive solution of the equation $\eta^5-|x_0|\eta=C(y_0/\eta^2)$. Then
\begin{equation}
 \label{ineq::X(x_0,y_0)}
\eta_{y_0}^4(x_0)=|x_0|+
\frac{C\Big(\frac{y_0}{\eta_{y_0}^2(x_0)}\Big)}{\eta_{y_0}(x_0)}<
|X_{min}(x_0,y_0)|<
|x_0|+\frac{C\Big(\frac{y_0}{\sqrt{|x_0|}}\Big)}{|x_0|^{1/4}},
\end{equation}
\end{proposition}
\begin{proof}
The proof is similar to the one for the analogous estimates of
$X_{min}(x_0)$ in Corollary~\ref{Cor:X_min_asymptotics} with the help of
Proposition~\ref{Prop:Cv0}.
\end{proof}
\begin{remark}
 \label{Rem:eta-algebraic}
Left Inequality~(\ref{ineq::X(x_0,y_0)}) remains valid if instead of the
function $\eta_{y_0}(x_0)$ satisfying the transcendental equation
one substitutes the function $\hat\eta_{y_0}(x_0)$ which is defined as
the unique positive solution of the algebraic equation,
$\hat\eta^5-|x_0|\hat\eta=C(v_0)$, where
$$
v_0=\left\{\begin{aligned}
&\frac{y_0}{\sqrt{|x_0|}},&y_0\leq0\\
&\frac{y_0}{\sqrt{|x_0|+\kappa(x_0,y_0)}},&y_0\geq0
\end{aligned}\right.
\qquad\mathrm{and}\qquad
\kappa(x_0,y_0)=\frac{C\Big(\frac{y_0}{\sqrt{|x_0|}}\Big)}{|x_0|^{1/4}}.
$$
This makes the inequality less accurate.
\end{remark}
\begin{remark}
It is easy to see in a similar way as for $X_{min}(x_0)$, that there is at
least one solution $y(x;x_0,y_0,y_l)$ whose minimum equals
$X_{min}(x_0,y_0)$. However, we can only conjecture the uniqueness of this
solution.
\end{remark}
\begin{remark}
Concerning the behaviour of $X_{min}(x_0,y_0)$ for $x_0$: $-1<x_0\leq0$
we can make an analogous comment as in Remark~\ref{Rem:X_min_0} with
the change $C$ on $C(y_0/\eta_{y_0}^2(x_0))$.
\end{remark}
\begin{proposition}
 \label{Prop:Z-above}
$$
|Z(x_0,y_0,y^0)|<|X_{min}(x_0,y_0)|+\frac{C(v^0)}{|X_{min}(x_0,y_0)|^{1/4}},
\qquad v^0=\frac{y^0}{\sqrt{|X_{min}(x_0,y_0)|}}.
$$
If $y^0\geq\tilde y_{lm}$, where $\tilde y_{lm}$ is the minimum value of
the solution $y(x;x_0,y_0,\tilde y_{lm})$ with the minimum at
$X_{min}(x_0,y_0)$, then $|X_{min}(x_0,y_0)|\leq|Z(x_0,y_0,y^0)|$, where
the equality is only possible if $y^0=\tilde y_{lm}$.
\end{proposition}
\begin{proof}
Denote $y_{lm}$ the minimum value of the maximal solution. Consider
the solution $y(x;X_{min}(x_0,y_0),y_{lm},y_{lm})$. It follows from
Lemma~\ref{Lem:4.2} that the graph of  $y(x;X_{min}(x_0,y_0),y_{lm},y_{lm})$
crosses the straight line $y=y^0$ at some point $P$ which is located to the
left of the point where the maximal solution crosses this line. This  means
that $|Z(x_0,y_0,y^0)|-|X_{min}(x_0,y_0)|<|P_x|-|X_{min}(x_0,y_0)|$, where
$P_x$ is the abscissa of $P$. The estimation of the last difference is
similar to the one for $X_-(x_0)$ in Lemma~\ref{Lem:2.2} and yields the
announced result.

The second statement is evident as the graph of $y(x;x_0,y_0,y_{lm})$ crosses
the straight line $y=y^0$ to the left of $X_{min}(x_0,y_0)$.
\end{proof}
\begin{proposition}
 \label{Prop:Z-below}
Let $y_{m0}=\min\{y_0,y^0\}$ and $\eta_{y_{m0}}(x_0)$ be the unique
positive solution of the equation $\eta^5-|x_0|\eta=2C(y_{m0}/\eta^2)$.
Then
\begin{equation}
 \label{ineq:belowZ(x_0,y_0,y^0)}
\eta_{y_{m0}}^4(x_0)=|x_0|+
\frac{2C\Big(\frac{y_{m0}}{\eta_{y_{m0}}^2(x_0)}\Big)}{\eta_{y_{m0}}(x_0)}<
|Z(x_0,y_0,y^0)|
\end{equation}
\end{proposition}
\begin{proof}
The proof is very similar to the one for the second
Inequality~(\ref{ineq:1thm1}) for $X(x_0)$ in Theorem~\ref{Thm:1}) with the
help of monotonicity of the function $C(v_0)$.
\end{proof}
\begin{remark}
As in Remark~\ref{Rem:eta-algebraic} we can, with a loss of accuracy,
substitute the function $\eta_{y_{m0}}(x_0)$ satisfying the transcendental
equation with the function $\hat\eta_{\hat v_0}(x_0)$ which is defined as
the unique positive solution of the algebraic equation,
$\hat\eta^5-|x_0|\hat\eta=2C(v_0)$, where now
$$
\hat v_0=\left\{\begin{aligned}
&\frac{y_{m0}}{\sqrt{|X_{min}(x_0,y_0)|}},&y_0\leq0\\
&\frac{y_{m0}}{\sqrt{|X_{min}(x_0,y_0)|+\hat\kappa(x_0,y_0)}},&y_0\geq0
\end{aligned}\right.
\quad\mathrm{and}\quad
\hat\kappa(x_0,y_0)=\frac{C\Big(\!
\frac{y_{m0}}{\sqrt{|X_{min}(x_0,y_0)|}}\!\Big)}{|X_{min}(x_0,y_0)|^{1/4}}.
$$
\end{remark}
\section{Spacing of Zeroes}
 \label{App:zero-spacing}
\begin{definition}
 \label{Def:deltas}
Let  $x_0\leq0$, $y_l<0$, and $y_0\geq0$.
Denote $z_2^-<z_1^-<x_0$ the zeroes of the solution $y(x;x_0,y_0,y_l)$ and
$x_0<z_2^+<z_1^+$ the zeroes of the solution $y_+(x;x_0,y_0,y_l)$ with the
initial slope $y_1<0$ at $x_0\leq X(0)$.
Define
$$
\delta_\pm(x_0,y_0,y_l)=z_1^\pm-z_2^\pm,\qquad
\delta_\pm(x_0,y_0)=\underset{y_l<0}{\sup}\;\delta_\pm(x_0,y_0,y_l).
$$
\end{definition}
\begin{proposition}
\label{Prop:delta-existence}
The functions $\delta_\pm(x_0,y_0)$ are finite for all $x_0$ and $y_0$,
in their domains of definition. For given $x_0$ and $y_0$ there exists
a value of $y_l=\hat y_l$ such that the solution
$y(x;x_0,y_0,\hat y_l)$ has the zero spacing $\delta_\pm(x_0,y_0)$.
In the case $y_0>0$ the number $\hat y_l$ for $\delta_-(x_0,y_0)$ is less
than the minimum value of the maximal solution corresponding to the segment
$\left[Z(x_0,y_0,0),x_0\right]$ (see \S~{\rm\ref{sect:boundary}},
Proposition~{\rm\ref{Prop:z-unique}}), and $\hat y_l$ for $\delta_+(x_0,y_0)$
is less than the minimum value of the maximal solution corresponding to the
segment $\left[x_0,Z^{-1}(x_0,0,y_0)\right]$, where
$\hat x_0=Z^{-1}(x_0,0,y_0)$ means that  $x_0=Z(\hat x_0,0,y_0)$.
\end{proposition}
\begin{proof}
It is clear that $0<\delta_-(x_0,y_0)<x_0-Z(x_0,y_0,0)$ and
$0<\delta_+(x_0,y_0)<Z^{-1}(x_0,0,y_0)-x_0$.

The functions $\delta_\pm(x_0,y_0,y_l)$ are smooth since they are
differences of two smooth functions, $z_1^\pm(x_0,y_0,y_l)$ and
$z_2^\pm(x_0,y_0,y_l)$. That the functions $z_1^\pm$ and $z_2^\pm$ are
smooth follows from the implicit function theorem, due to the fact that
$y'(z_k^\pm;x_0,y_0,y_l)\neq0$ for $k=1,2$.

The functions $\delta_\pm(x_0,y_0,y_l)$ are bounded as $y_l\to-\infty$; for
$\delta_-(x_0,y_0,y_l)$ it follows from Proposition~\ref{Prop:4.3}, for
$\delta_+(x_0,y_0,y_l)$ an analogous estimate as in
Proposition~\ref{Prop:4.3} clearly holds. Thus, we get existence of finite
value of $y_l=\hat y_l$, which, of course, depends on $x_0$, $y_0$, and on
the sign $\pm$, such that
$\delta_\pm(x_0,y_0)=\delta_\pm(x_0,y_0,\hat y_l)$.

Now we prove that $\hat y_l$ for $\delta_-(x_0,y_0,y_l)$ is less than
the minimum value of the solution corresponding to $Z(x_0,y_0,0)$.
The proof of the corresponding statement for $\delta_+(x_0,y_0,y_l)$ is
absolutely analogous. Let $y_{lm}$ be the minimum value of the maximal
solution, i.e., $Z(x_0,y_0,0)=\mathcal Z_-(x_0,y_0,0,y_{lm})$. The graphs
of solutions $y(x;x_0,y_0,y_l)$ with $y_l:\,y_{lm}<y_l<0$ have only one
point of intersection, $(x_0,y_0)$. Therefore, their first zero, $z_1^-$,
is monotonically increasing and the second, $z_2^-$, is monotonically
decreasing as $y_l\searrow y_{lm}$. These mean that for such values of
$y_l$ the function $\delta_-(x_0,y_0,y_l)$ is monotonically increasing.
Thus $\hat y_l\leq y_{lm}$. For $y_l<y_{lm}$, $z_1^-$ continues its monotonic
increase, while $z_2^-$ begins a monotonic increase too. However since
$y_{lm}$ is the minimum of $z_2^-$, for some $\epsilon>0$ we have
$\Delta z_2^-=\mathcal O((\Delta y_l)^2)$, while
$\Delta z_1^-=\mathcal O(\Delta y_l)$, where $\Delta y_l=y_{lm}-y_l<\epsilon$.
Therefore, $\Delta z_1^--\Delta z_2^->0$ in some lower neighbourhood of
$y_{lm}$, which means that $\hat y_l<y_{lm}$.
\end{proof}
\begin{proposition}
The functions  $\delta_\pm(x_0,y_0)$ are continuous functions of
two variables.
\end{proposition}
\begin{proof}
Let $\hat y_l$ be defined as in Proposition~\ref{Prop:delta-existence},
i.e., $\delta_\pm(x_0,y_0)=\delta_\pm(x_0,y_0,\hat y_l)$. Consider a
sequence of points $(x_n,y_n)\to(x_0,y_0)$ and denote by $\hat y_{ln}$
a corresponding sequence of the levels, such that
$\delta_\pm(x_n,y_n)=\delta_\pm(x_n,y_n,\hat y_{ln})$. The sequence
$\hat y_{ln}$ is bounded, since for any unbounded subsequence
$\hat y_{ln_k}$ the corresponding values of
$\delta_\pm(x_{n_k},y_{n_k},\hat y_{ln_k})\to0$
(cf. Proposition~\ref{Prop:4.3}).

Suppose that there is a subsequence
$\delta_\pm(x_{n_k},y_{n_k},\hat y_{ln_k})$ that does not converge
to $\delta_\pm(x_0,y_0)$. Changing the notation, if necessary, we can
assume that it converges to some $\delta<\delta_\pm(x_0,y_0)$.
Consider another subsequence $\delta_\pm(x_{n_k},y_{n_k},\hat y_l)$.
It converges to $\delta_\pm(x_0,y_0)$. Therefore, for all $n>N$,
$\delta_\pm(x_{n_k},y_{n_k},\hat y_l)>
\delta_\pm(x_{n_k},y_{n_k},\hat y_{ln_k})$, which contradicts the
definition of the numbers $\hat y_{ln_k}$.
\end{proof}
\begin{proposition}
 \label{Prop:delta-}
For $\hat x_0<x_0\leq0$,
$$
0<\delta_-(x_0,y_0,y_l)-\delta_-(\hat x_0,y_0,y_l)<x_0-\hat x_0,\qquad
0<\delta_-(x_0,y_0)-\delta_-(\hat x_0,y_0)<x_0-\hat x_0.
$$
For $\hat y_0>y_0\geq0>y_l$,
\begin{equation}
 \label{ineq:delta-level-monotonic}
0<\delta_-(x_0,y_0,y_l)-\delta_-(x_0,\hat y_0,y_l)<
\frac{\hat y_0-y_0}{f(x_0;y_0,y_l)},
\end{equation}
where the function $f(x_0;y_0,y_l)$ is defined in
Corollary~{\rm\ref{Cor:3.1}}. Moreover,
\begin{equation}
 \label{ineq:delta-monotonic}
0<\delta_-(x_0,y_0)-\delta_-(x_0,\hat y_0)<
\frac{x_0-X_{min}(x_0,y_0)}{y_0-y_{lm}(x_0,y_0)}(\hat y_0-y_0),
\end{equation}
where the function $y_{lm}(x_0,y_0)$ is the minimum value of the maximal
solution corresponding to the segment $[Z(x_0,y_0,0),x_0]$ and
$X_{min}(x_0,y_0)$ is given in Definition~{\rm\ref{Def:X_{min}(x_0,y_0)}}.
\end{proposition}
\begin{proof}
Consider solutions $y_1(x)=y(x;x_0,y_0,y_l)$ and
$y_2(x)=y(x;\hat x_0,y_0,y_l)$.
Denote by $x_1(y)$ and $x_2(y)$, respectively, their inverse functions on
the segments bounded by minima and $x_0$. The difference
$x_1(y)-x_2(y)$ is monotonically decreasing with $y\searrow y_l$, see
the proof of Proposition~\ref{Prop:3.5}.
Monotonic decrease of the difference continues with $y\nearrow+\infty$
if we consider the inverse functions on the segments bounded by
the left poles and the minima. Thus we have the following inequalities:
$$
0<\mathcal Z_-(x_0,y_0,y^0,y_l)-\mathcal Z_-(\hat x_0,y_0,y^0,y_l)<
z_2^--\hat z_2^-<z_1^--\hat z_1^-<x_0-\hat x_0,
$$
where $y^0\geq y_0$ is arbitrary and the ``hat-notation'' is used to
distinguish the zeroes of solution $y_2(x)$. This proves both inequalities
for the difference $\delta_-(x_0,y_0,y_l)-\delta_-(\hat x_0,y_0,y_l)$.
Now choosing in the last difference $y_l=\hat{\hat{y_l}}$ such that
$\delta_-(\hat x_0,y_0,\hat{\hat{y_l}})=\delta_-(\hat x_0,y_0)$
(see Proposition~\ref{Prop:delta-existence}) and
increasing then, if necessary, $\delta_-(x_0,y_0,\hat{\hat{y_l}})$ to
$\delta_-(x_0,y_0)$ we arrive at the left inequality for
$\delta_-(x_0,y_0)-\delta_-(\hat x_0,y_0)$. The right inequality is
obtained by choosing $y_l=\hat y_l$ such that
$\delta_-(x_0,y_0,\hat y_l)=\delta_-(x_0,y_0)$ and increasing,
if necessary, $\delta_-(\hat x_0,y_0,\hat y_l)$.

Consider now solutions $y_1(x)=y(x;x_0,y_0,y_l)$ and
$y_3(x)=y(x;x_0,\hat y_0,y_l)$. Denote $\hat x_0$ the right point of
intersection of the graph of $y_3(x)$ with the straight line $y=y_0$. Then
we can apply already proved inequalities to the zero spacing
functions of $y_1(x)$ and $y_3(x)$. Using the concavity of graphs of
the solutions and monotonicity of the function $f(x_0;y_0,y_l)$ one writes,
$$
f(x_0;y_0,y_l)(x_0-\hat x_0)<f(\hat x_0;y_0,y_l)(x_0-\hat x_0)
<\hat y_0-y_0<f(x_0;\hat y_0,y_l)(x_0-\hat x_0).
$$
This completes the proof of Inequality~(\ref{ineq:delta-level-monotonic}).
Inequality~(\ref{ineq:delta-monotonic}) is a consequence of
Inequality~(\ref{ineq:delta-level-monotonic}) and the estimate
$f(x_0;y_0,y_l)(x_0-X_{min}(x_0, y_0))>y_0-y_{lm}(x_0,y_0)$,
which follows from the concavity and Proposition~\ref{Prop:delta-existence}.
\end{proof}
\begin{remark}
 \label{Rem:delta1}
It is interesting to note that the function  $|Z(x_0,y_0,+\infty)|$
introduced in \S~\ref{sect:boundary} (cf. Definition~\ref{Def:z-main}) that
measures a distance to the farthest (to the left) pole of the solutions with
initial value $y_0$ at $x_0$ is monotonically increasing with the increase
of $y_0$ (see Proposition~\ref{Prop:z-main}), while the behaviour of
$\delta_-(x_0,y_0)$ is directly opposite. It is natural to expect that the
solutions with larger intervals of existence have a larger spacing between
their zeroes, and it is actually true. However, it is true
only if we consider zeroes of solutions whose intervals of existence
contain the point $x_0$ rather than having it as their right bound.
\end{remark}
\begin{proposition}
 \label{Prop:delta+}
For $\hat x_0<x_0\leq\mathcal X(0,y_l)$,
\begin{equation}
 \label{ineq:delta+x_level-monotonic}
0<\delta_+(x_0,y_0,y_l)-\delta_+(\hat x_0,y_0,y_l)<
(L_\Xi(y_l)-1)(x_0-\hat x_0),
\end{equation}
where
$
L_\Xi(y_l):=\underset{x_0\leq\mathcal X(0,y_l)}{\sup}\;
\frac\partial{\partial x_0}\varXi(x_0,y_l)
$
satisfies the following inequalities
$$
1<L_\Xi(y_l)\leq L_\Xi
$$ with
$L_\Xi$ defined in Definition~{\rm\ref{Def:Lipschitz}}
of \S~{\rm\ref{sect:5}}.\\
For $\hat x_0<x_0\leq X(0)$ and $\hat y_0>y_0\geq0>y_l$:
\begin{eqnarray}
 \label{ineq:delta+x-monotonic}
0<\delta_+(x_0,y_0)-\delta_+(\hat x_0,y_0)<(L_\Xi-1)(x_0-\hat x_0);&&
\end{eqnarray}
\begin{equation}
 \label{ineq:delta+y_level-monotonic}
0<\delta_+(x_0,\hat y_0,y_l)-\delta_+(x_0,y_0,y_l)<
\frac{L_\Xi(y_l)-1}{|f^+(x_0;y_0,y_l)|}(\hat y_0-y_0),
\end{equation}
where the function $f^+(x_0;y_0,y_l)$ is defined in
Remark~{\rm\ref{Rem:symmetry}};
\begin{equation}
 \label{ineq:delta+y-monotonic}
0<\delta_+(x_0,\hat y_0)-\delta_+(x_0,\hat y_0)<
(L_\Xi-1)\frac{\Xi_{min}(x_0)-x_0}{y_0-y_{lm}(x_0,y_0)}(\hat y_0-y_0),
\end{equation}
where the function $y_{lm}(x_0,y_0)$ is the minimum value of the maximal
solution corresponding to the segment $[x_0,Z^{-1}(x_0,0,y_0)]$.
\end{proposition}
\begin{proof}
Consider solutions $y_1(x)=y_+(x;x_0,y_0,y_l)$ and
$y_2(x)=y_+(x;\hat x_0,y_0,y_l)$.
Denote by $x_1(y)$ and $x_2(y)$, respectively, their inverse functions on
the segments bounded by minima and $x_0$. The difference
$x_1(y)-x_2(y)$ is monotonically increasing with $y\searrow y_l$, see
the proof of Proposition~\ref{Prop:3.5} (cf. Proposition~\ref{Prop:5.5}).
Monotonic increase of the difference continues with $y\nearrow+\infty$
if we consider the inverse functions on the segments bounded by
the left poles and the minima. Thus we have the following inequalities:
$$
0<x_0-\hat x_0<z_2^+-\hat z_2^+<
\mathcal Z_+(x_0,y_0,+\infty,y_l)-\mathcal Z_+(\hat x_0,y_0,+\infty,y_l),
$$
where the ``hat-notation'' is used to distinguish the zeroes of solution
$y_2(x)$. Now, we notice that
$
\mathcal Z_+(x_0,y_0,+\infty,y_l)-\mathcal Z_+(\hat x_0,y_0,+\infty,y_l)=
\varXi(\tilde x_0,y_l)-\varXi(\tilde{\hat x}_0,y_l),
$
for some $\tilde x_0$ and $\tilde{\hat x}_0$, such that
$\tilde x_0-\tilde{\hat x}_0<x_0-\hat x_0$, by virtue of the monotonicity.
This proves both Inequalities~(\ref{ineq:delta+x_level-monotonic}) and
(\ref{ineq:delta+x-monotonic}), since according to
Proposition~\ref{Prop:XiLipschitz} $L_\Xi$ is finite.

Consider solutions $y_1(x)=y_+(x;x_0,y_0,y_l)$ and
$y_3(x)=y(x;x_0,\hat y_0,y_l)$. Denote by $\tilde x_0$ the left point of
intersection of the graph of $y_1(x)$ with the straight line $y=\hat y_0$.
Then we can apply already proved
Inequalities~(\ref{ineq:delta+x_level-monotonic}) for the zero spacing
functions of $y_1(x)$ and $y_3(x)$ with $y_0$ changed to $\hat y_0$.
Now, taking into account that
$\delta_+(\tilde x_0,\hat y_0,y_l)=\delta_+(x_0,y_0,y_l)$
and, by concavity of the graph of $y_1(x)$,
$|f^+(x_0;y_0,y_l)|(x_0-\tilde x_0)<\hat y_0-y_0$ we arrive at
Inequalities~(\ref{ineq:delta+y_level-monotonic}).
Inequalities~(\ref{ineq:delta+y-monotonic}) are the consequence of the
ones for the level functions~(\ref{ineq:delta+y_level-monotonic}) and the
estimate $|f^+(x_0;y_0,y_l)|(\Xi_{min}(x_0)-x_0)>y_0-y_{lm}(x_0,y_0)$,
which follows from the concavity and Proposition~\ref{Prop:delta-existence}.
\end{proof}
\begin{corollary}
The functions $\delta_\pm(x_0,0)$ are smooth and have the following
properties:
$$
\delta_-(x_0,0)=\underset{y_0\geq0}{\sup}\;\delta_-(x_0,y_0),\qquad
\delta_+(x_0,0)=\underset{y_0\geq0}{\inf}\;\delta_+(x_0,y_0).
$$
For all $x_0<0$
$$
\delta_+(x_0,0)\underset{x_0\to-\infty}=\frac{2C(0)}{|x_0|^{1/4}}+
o\left(\frac1{|x_0|^{3/2}}\right)>
\delta_-(x_0,0)\underset{x_0\to-\infty}=\frac{2C(0)}{|x_0|^{1/4}}+
\mathcal{O}\left(\frac1{|x_0|^{3/2}}\right),
$$
where $C(0)$ is given in Definition~{\rm\ref{Def:Cv0}} for $v_0=0$ and
$-\frac{C^2(0)}{|x_0|^{3/2}}<\mathcal{O}\left(\frac1{|x_0|^{3/2}}\right)<0$.
\end{corollary}
\begin{proof}
The extremal properties of $\delta_\pm(x_0,0)$ follows from
Propositions~\ref{Prop:delta-} and \ref{Prop:delta+}.
The smoothness is the consequence of the
representations
\begin{equation}
 \label{eq:delta_0}
\delta_\pm(x_0,0)=|Z^{\mp1}(x_0,0,0)-|x_0||,
\end{equation}
where
$Z^{+1}(x_0,0,0)\equiv Z(x_0,0,0)$ and $Z^{-1}(x_0,0,0)$ is the the
section of the fiber $Z^{-1}(x_0)$ with the straight line passing through
the origin in direction of the vector $(1,0,0)$. Now the uniqueness of
the maximal solutions (see Proposition~\ref{Prop:z-unique}) allows one to
prove smoothness in the spirit of Corollary~\ref{Cor:smoothness}.
Of course, smoothness of $\delta_-(x_0,0)$ follows directly from
representation~(\ref{eq:delta_0}) and Proposition~\ref{Prop:z-main}.
The estimates are special cases of the estimates proved in
Propositions~\ref{Prop:Z-above} and \ref{Prop:Z-below}.
\end{proof}
\begin{remark}
Numerical calculations shows that the integral (see Definition~\ref{Def:Cv0})
whose maximum equals $C(0)$ has only one extremum at
$$
x_{max}=2.04124321493909675\ldots,\;\;
v_{min}^{max}=-\frac1{\sqrt{2x_{max}}}=-0.49492298670007907\ldots
$$
the accuracy of $18$ digits in $x_{max}$ allows one to determine $36$
digits of $C(0)$,
$$
C(0)=0.69663587640019346382935052393992493\ldots
$$
The minimum value $y_{lm}$ of the maximal solution corresponding to the segment
$[Z(x_0,0,0),x_0]$ has the following asymptotics,
$$
y_{lm}=\underset{x_0\to-\infty}=v_{min}^{max}\sqrt{|x_0|}+
o\left(\sqrt{|x_0|}\right).
$$
For a better numerical approximation in the last formula is preferable to
use instead of $x_0$ the corresponding minimum of the maximal solution.
\end{remark}
\begin{remark}
If the solution with the zero spacing equal to $\delta_\pm(x_0,y_0)$
(cf. Proposition~\ref{Prop:delta-existence}) is unique, then smoothness
of $\delta_\pm(x_0,y_0)$ with respect to both variables can be immediately
established by the analogous proof as the one for smoothness of $X(x_0)$
in \S~\ref{sect:main} Corollary~\ref{Cor:smoothness}.
\end{remark}
\begin{remark}
We define functions $\delta_\pm(x_0,y_0)$ in the natural domain $y_0\geq0$,
however it is possible to define their continuation into the lower
half-plane, $y_0<0$. To do so we consider solutions with initial
values at the point $(x_0,y_0)$ with the zeroes, $z_2<x_0<z_1$ for
$\delta_-(x_0,y_0)$, and  $z_1<x_0<z_2$ for $\delta_+(x_0,y_0)$ and the
supremum in Definition~\ref{Def:deltas} is taken over non-negative slopes
for $\delta_-(x_0,y_0)$ and non-positive for $\delta_+(x_0,y_0)$. Then
in some neighbourhood of $y_0=0$ function  $\delta_-(x_0,y_0)$ continues
its monotonic increase, respectively function $\delta_+(x_0,y_0)$ its
decrease, with the decrease of $y_0$. This
follows by the same arguments as in the proof of monotonicity of these
functions for $y_0\geq0$ in Propositions~\ref{Prop:delta-} and
\ref{Prop:delta+}. The proof works for all values of $y_0$ such that
the minima of solutions corresponding to $\delta_-(x_0,y_0)$ are located
to the left of $x_0$ and, respectively, to the right of $x_0$ for
$\delta_+(x_0,y_0)$. Since we know that
$\underset{y_0\to-\infty}\lim\delta_\pm(x_0,y_0)=0$, both minima have to
reach abscissa $x_0$ at some finite values (maybe different) $y_-(x_0)$ and
$y_+(x_0)$. It is a further interesting question whether actually
$y_-(x_0)=y_+(x_0)=y_0(x_0)$ and for $y_0\leq y_0(x_0)$,
$\delta_-(x_0,y_0)=\delta_+(x_0,y_0)$?
\end{remark}

\end{document}